\documentclass[10pt]{article} % For LaTeX2e
\usepackage[accepted]{tmlr}
\usepackage{amsmath,amsfonts,bm}

\def\eqref#1{equation~\ref{#1}}
\def\1{\bm{1}}

\DeclareMathAlphabet{\mathsfit}{\encodingdefault}{\sfdefault}{m}{sl}
\SetMathAlphabet{\mathsfit}{bold}{\encodingdefault}{\sfdefault}{bx}{n}

\usepackage{graphicx,graphics,epsfig,epstopdf}
\usepackage{hyperref}
\usepackage{url}
\usepackage{natbib}
\usepackage[utf8]{inputenc} % allow utf-8 input
\usepackage[T1]{fontenc}    % use 8-bit T1 fonts
\usepackage{url}            % simple URL typesetting
\usepackage{booktabs}       % professional-quality tables
\usepackage{amsfonts}       % blackboard math symbols
\usepackage{amssymb,amsmath,amsthm}
\usepackage{nicefrac}       % compact symbols for 1/2, etc.
\usepackage{microtype}      % microtypography
\usepackage{tikz}
\usepackage{caption}
\usepackage{wrapfig}
\usepackage[linesnumbered,ruled,vlined,algo2e,algoruled,boxed]{algorithm2e}
\usepackage{subcaption}
\usepackage{soul}
\usepackage{mathtools}
\usepackage{accents}
\usepackage{setspace}
\usepackage{multirow,makecell}
\usepackage{threeparttable}
\usepackage{optidef}

\SetKwInput{kwAdd}{Add}
\SetKwInput{kwSet}{Set}
\SetKwInput{kwReset}{Reset}
\SetKwInput{kwUpdate}{Update}
\SetKwInput{kwClear}{Clear}
\SetKwInput{kwSave}{Save}
\SetKwInput{kwBuild}{Build}
\SetKwInput{kwIni}{Initialization}
\usetikzlibrary{matrix,chains,positioning,decorations.pathreplacing,arrows}

\newtheorem{lemma}{Lemma}
\newtheorem{remark}{Remark}
\newtheorem{theorem}{Theorem}
\newtheorem{assumption}{Assumption}

\DeclareMathOperator{\rank}{rank}
\newcommand{\col}{{\rm col\;}}

\newcommand{\bI}{{\mathbf I}}

\usepackage{hyperref}

\DeclareMathAlphabet\mathbfcal{OMS}{cmsy}{b}{n}

\title{Online Control-Informed Learning}

\author{\name Zihao Liang \email liang331@purdue.edu \\
      \addr School of Aeronautics and Astronautics\\
  Purdue University
      \AND
      \name Tianyu Zhou \email zhou1043@purdue.edu \\
      \addr School of Aeronautics and Astronautics\\
  Purdue University
      \AND
      \name Zehui Lu \email lu846@purdue.edu\\
      \addr School of Aeronautics and Astronautics\\
  Purdue University
    \AND
      \name Shaoshuai Mou \email  mous@purdue.edu\\
      \addr School of Aeronautics and Astronautics\\
  Purdue University}

\def\month{03}  % Insert correct month for camera-ready version
\def\year{2025} % Insert correct year for camera-ready version
\def\openreview{\url{https://openreview.net/forum?id=LDzvZEVl5H}} % Insert correct link to OpenReview for camera-ready version

\begin{document}

\maketitle

\begin{abstract}
This paper proposes an Online Control-Informed Learning (OCIL) framework, which employs the well-established optimal control and state estimation techniques in the field of control to solve a broad class of learning tasks in an online fashion. This novel integration effectively handles practical issues in machine learning such as noisy measurement data, online learning, and data efficiency. By considering any robot as a tunable optimal control system, we propose an online parameter estimator based on extended Kalman filter (EKF) to incrementally tune the system in an online fashion, enabling it to complete designated learning or control tasks. The proposed method also improves the robustness in learning by effectively managing noise in the data. Theoretical analysis is provided to demonstrate the convergence of OCIL. Three learning modes of OCIL, i.e. Online Imitation Learning, Online System Identification, and Policy Tuning On-the-fly, are investigated via experiments, which validate their effectiveness. 
\end{abstract}

\section{Introduction}

Informed Machine Learning (IML) \citep{von2021informed} represents an emerging approach that integrates prior knowledge into the machine learning (ML) process. While classic classification tasks in unsupervised, semi-supervised, or supervised ML primarily focus on extracting patterns from labeled or unlabeled data \citep{lecun2015deep}, IML leverages prior knowledge such as physical laws, expert knowledge, or existing models to uncover underlying connections within data \citep{karniadakis2021physics}. This integration enables models to produce more reliable and interpretable predictions, especially when dealing with noisy data. This approach is especially advantageous in the domains where theoretical understanding is well-established and thus can guide ML. One notable example of IML is physics-informed machine learning \citep{wu2018physics,karniadakis2021physics,kashinath2021physics}, which is particularly valuable for solving partial differential equations for computational fluid dynamics.

Control-informed learning (CIL) is a subset of IML tailored for system control, autonomy, and robotics. 
This approach merges standard control theory with ML techniques to enhance the capabilities of autonomous systems. The integration leverages the complementary strengths of control and learning. Control theory provides model structures and optimization guidance that enable efficient and reliable algorithms for handling complex tasks. Meanwhile, ML improves these models by learning from data, a capacity that some conventional control methods lack \citep{jin2020pontryagin,jin2021safe}. This paper aims to tackle learning tasks in autonomous systems that are governed by optimal control (OC) systems. An optimal control system usually consists of dynamics, a control policy, and an objective function. From a unified perspective, learning these components can be understood as addressing a common problem with unknown parameters in different parts of the system and using different loss functions. 
For example, in learning dynamics, the task involves parameterizing a differential equation, with the loss function measuring the prediction error between the model’s output and the target data. In learning policies, the unknown parameters are within the feedback policy and the control objective itself serves as a loss function. When learning control objective functions, the objective is parameterized, and the loss measures the discrepancy between the reproduced trajectory and observed demonstrations.

To tackle these problems, many works in the field of so-called Learning for Dynamics and Control aim to leverage the integration of learning and control but often treat them as separate or sequential tasks. For example, control theories are used for algorithm development and convergence analysis of online unconstrained or constrained optimizations \citep{casti2023control,10189107, 10771962}; model-based reinforcement learning \citep{heess2015learning,gu2016continuous}, improves sample efficiency by using dynamics models; Koopman-operator control \citep{proctor2018generalizing,abraham2019active,hao2024distributed}, employs learning to transform nonlinear systems into linear observable space, simplifying control design.
In contrast, CIL integrates these processes, allowing learning algorithms to incorporate control insights directly. The integration enables ML and control techniques to perform simultaneously, reducing computational complexity, and improving practical applicability. 
CIL differentiates itself by utilizing Pontryagin's maximum principle, a foundational concept in OC theory. This principle defines the optimality conditions for the state and input trajectories of an OC system. CIL employs these conditions to provide gradients for machine learning \citep{jin2020pontryagin,jin2021safe,bottcher2022ai}. CIL integrates these gradients directly into its learning process, ensuring that machine learning outcomes are efficient while remaining consistent with established control theories and physical models. This approach enhances both the reliability and accuracy of the results.

\subsection{Related Work}

This section presents existing research on learning various components of an autonomous system and explores related learning frameworks that address these problems from a unified perspective.

\textbf{Learning dynamics.} 
To learn a nonlinear system with possibly noisy measurement, Markov decision-process-based methods are widely used, such as linear regression \citep{haruno2001mosaic}, observation-transition modeling \citep{finn2016unsupervised}, latent space modeling \citep{watter2015embed}, (deep) neural networks (NN) \citep{li2018maximum,li2018optimal,han2019mean,zhang2019you,benning2019deep,liu2020selection,beintema2023deep,pillonetto2025deep}, Gaussian processes \citep{deisenroth2011pilco}, and transition graphs \citep{zhang2018composable}.
Despite their widespread use, these methods often must balance data efficiency with prediction accuracy.
To improve both metrics, physics-informed learning approaches \cite{lutter2019deep,xu2020second,saemundsson2020variational,sharma2023review} incorporate physical laws into learning models.
Koopman operator theory offers a method for lifting states to an infinite-dimensional linear observable space \citep{mauroy2020koopman,liang2023data,hao2023policy,liu2024koopa}.

\textbf{Learning objective functions.} 
Objective learning is typically referred to as inverse reinforcement learning (IRL) in the ML community and inverse optimal control (IOC) in the system control community.
These methods aim to deduce a control objective function with observed optimal demonstrations. \citep{brown2019extrapolating}
The objective function is generally represented as a weighted sum of features \citep{abbeel2004apprenticeship,ratliff2006maximum,ziebart2008maximum,arora2021survey}.
Approaches to find these unknown weights include feature matching \citep{abbeel2004apprenticeship}, maximum entropy \citep{ziebart2008maximum}, maximum margin \citep{ratliff2006maximum}, and approximate variational reward imitation learning \citep{chan2021scalable}.
As for learning nonlinear parameter mapping of objective functions, prior and system-dependent knowledge is required to further extend the methods above.
On the other hand, with system dynamics, IOC aims for efficient learning approaches \citep{mombaur2010human}.
For example, some methods \citep{keshavarz2011imputing,jin2019inverse,jin2021inverse,jin2021distributed,liang2022iterative,liang2023data} directly calculate unknown weights by minimizing the violation of optimality conditions by the observed demonstration data, which avoids repeatedly solving OC problems. 

\textbf{Learning control policies.} Learning policies are generally termed reinforcement learning (RL) and OC in the ML and control communities, respectively. In the RL community, there are mainly two streams of research, namely model-free and model-based RL. Model-free RL learns policies by directly interacting with the environment, without using a model of it \citep{mnih2013playing,mnih2015human,oh2016control}. To improve data complexity, model-based RL learns a dynamics model before policy learning \citep{schneider1997exploiting,abbeel2006using,deisenroth2011pilco,levine2014learning,gu2016continuous}. For OC, the first strategy is based on dynamical programming, such as the linear quadratic regulator (LQR) \citep{scokaert1998constrained}, which solves the OC problem with linear dynamics and quadratic cost, the linear quadratic Gaussian \citep{athans1971role}, which combines LQR with a Kalman filter to solve OC problem with linear system affected by Gaussian noise, the iterative linear quadratic regulator (iLQR) \citep{li2004iterative}, which linearizes the dynamics and quadratizes the value function, and differential dynamical programming, which quadratizes the dynamics and value function. Another strategy relies on Pontryagin's maximum/minimal principle (PMP) \citep{pontryagin2018mathematical}, such as shooting methods \citep{bock1984multiple} and collocation methods \citep{patterson2014gpops}. These open-loop methods are further improved by closed-loop methods such as model predictive control (MPC) \citep{schwenzer2021review}, which repeatedly solves an OC problem over a finite horizon to generate control inputs.
Recently, \cite{jin2020pontryagin} proposed a framework for learning an optimal policy based on differentiating Pontryagin's Maximum Principle.

Many research studies also focus on incremental policy tuning. One of the most popular tracks is transfer learning, which exploits the generalization of existing knowledge so that it can be transferred across different domains \citep{taylor2009transfer}.
Recently, transfer learning has been implemented to speed up the learning process in RL \citep{taylor2009transfer}.
Another popular method is behavior cloning \citep{torabi2018behavioral,czarnecki2019distilling,sasaki2021behavioral}.
In the control community, tuning OC systems initially refers to neighboring extremal optimal control (NEOC) \citep{bryson1975applied,ghaemi2009neighboring}. There are other popular methods including adaptive control \citep{ioannou2012robust,bertsekas2022lessons,luo2023adaptive,guo2023composite}, which adjusts its parameters in real-time to maintain optimal performance, even in the presence of uncertainties or changes in system dynamics, and Bayesian optimization for controller tuning, \citep{khosravi2021performance,sorourifar2021data,berkenkamp2023bayesian}.

To sum up, there are numerous existing methods focused on individual tasks. These approaches are effective when only one component of the system is unknown. However, in many real-world scenarios, multiple components may be unavailable or uncertain simultaneously. For instance, in autonomous driving, the dynamics of the vehicle may be unknown due to changes in road conditions or vehicle wear and tear. Simultaneously, the control policy may also be unavailable due to a lack of predefined rules or data. 
In such cases, existing methods often fall short, as they are not designed to handle the joint learning of multiple interdependent components, limiting their applicability in more complex or incomplete systems.

\textbf{Unified learning frameworks.} Several studies have explored unified learning frameworks to tackle learning challenges in autonomous systems. These approaches integrate an implicit planner directly within the policy \citep{okada2017path,pereira2018mpc,amos2018differentiable,srinivas2018universal}. The main challenge in these methods is learning the OC system, which is very similar to the goal of this work.
\citep{okada2017path,pereira2018mpc} learn a path-integral OC system \citep{kappen2005path}, which is a special class of OC systems. \citep{srinivas2018universal} learns an OC system in a latent space.
These methods rely on an ``unrolling'' strategy to make differentiation easier. 
Essentially, they treat solving an OC problem as an ``unrolled'' computational graph created by applying gradient descent repeatedly. This allows automatic differentiation tools \citep{abadi2016tensorflow} to be used directly. This approach faces a few challenges: (i) it requires storing all intermediate steps, making it memory-intensive, and (ii) the accuracy of the gradients depends on how many steps are included in the graph, leading to a trade-off between computational cost and accuracy.
To tackle these issues, \cite{amos2018differentiable} proposed a differentiable MPC framework. In the forward pass, it uses an LQR approximation of the OC system, and in the backward pass, gradients are computed by differentiating this LQR approximation. This framework has a major challenge: differentiating LQR requires solving a large linear equation, involving the inversion of a matrix with size proportional to the time horizon, making it very costly for long-horizon systems. 
To address the challenges of the framework mentioned above, \cite{jin2020pontryagin} proposed Pontryagin's differential programming (PDP). PDP avoids unrolled computational graphs by only storing the resulting trajectory without concern about how it is solved. Instead of relying on intermediate LQR approximations, it directly differentiates through Pontryagin's Maximum Principle (PMP) to obtain exact gradients. Furthermore, its backward pass uses an auxiliary control system to obtain the gradient, reducing memory and computational complexity.  First, it lacks the ability for online learning, as it relies on gradient descent to update unknown parameters in the OC system, requiring significant computation time to reach convergence. This drawback is particularly problematic in applications like autonomous driving, where quick adaptation to new scenarios is essential for safety and performance. Second, PDP does not account for noisy measurement data, limiting its effectiveness in real-world situations where sensor data is often unreliable or noisy.

% Online control-informed learning (OCIL) represents a subset of CIL that continually updates its model as new data streams. This is beneficial for robotic applications, where the environment can change unpredictably. One advantage of online methods is their ability to learn and adapt on the fly. Unlike offline methods, which require a static dataset for training, online learning algorithms adjust their parameters or models in response to new information without training from scratch. This continual learning capability allows autonomous systems to improve their policies in real time. Additionally, OCIL can handle dynamic environments by continually updating its parameters. This feature is critical in applications such as autonomous driving, where adapting to new scenarios quickly is crucial for safety and performance. This concept intersects with many  ML methods such as transfer learning \citep{pan2009survey}, continual learning \citep{aljundi2019task}, and learning on-the-fly \citep{ornik2019control}. 

% \zihao{highlight limitations of current frameworks}
% \zihao{summarize limitation of limited method and unified frameworks}

\subsection{Contributions}
This paper introduces an online learning framework called Online Control-Informed Learning (OCIL). This framework is designed to be data efficient for various learning and control tasks while providing robustness against noisy data. In this paper, we consider an autonomous system as an OC system, which is parameterized by tunable parameters within different components of the system, including dynamics, policy, and objective function. By tuning the OC system in an online fashion, the proposed OCIL tackles three learning tasks in robotics, namely Online Imitation Learning, Online System Identification, and Policy Tuning On-the-fly. The proposed OCIL consists of two main components, both of which are inspired by control theory. Specifically, the framework first proposes an online parameter estimator based on the classic online state estimation techniques in control theory.
The estimator continually updates the parameter estimates in an online fashion as new data becomes available, aiming to minimize a cumulative loss defined for a specific task. To do so, the gradient information for the loss with respect to the tunable parameter is required. Therefore, OCIL employs a gradient generator (GG) based on Pontryagin Differential Programming in OC theory to calculate the exact gradient.

\textbf{Notations.} 
$\lVert \cdot \lVert$ denotes the Euclidean norm. Given a matrix $A\in\mathbb{R}^{n\times m}$, let $A'$ denotes its transpose.
% $A^\dagger$ denotes its Moore-Penrose pseudoinverse.
For positive integers $n$ and $m$, let $\bI_{n}$ be the ${n \times n}$ identity matrix; $\mathbf{0}_n \in\mathbb{R}^n$ denotes a vector with all value $0$; $\mathbf{0}_{n\times m}$ denotes a $n\times m$ matrix with all value $0$. Let $\text{col}\{\boldsymbol{v}_1,\hdots,\boldsymbol{v}_a\}$ denote a column stack of elements $\boldsymbol{v}_1,\hdots,\boldsymbol{v}_a$, which may be scalars, vectors or matrices, i.e. $\text{col}\{\boldsymbol{v}_1,\hdots,\boldsymbol{v}_a\}\triangleq[\boldsymbol{v}_1^\prime \hdots \boldsymbol{v}_a^\prime]$.

\section{Problem Formulation}\label{ProF}

Consider the following class of OC systems $\Sigma(\boldsymbol{\theta}^*)$, where $\boldsymbol{\theta}^*\in\mathbb{R}^p$ denotes the unknown and constant parameter. The behavior of $\Sigma(\boldsymbol{\theta}^*)$ is determined by minimizing a control objective function:
% \begin{equation} \label{eq:optiOC}
% \begin{aligned}
% \{\boldsymbol{x}_{1:T}(\boldsymbol{\theta}^*),\boldsymbol{u}_{0:T-1}(\boldsymbol{\theta}^*)\}\in&\arg\min_{\substack{\boldsymbol{x}_{0:T}\\\boldsymbol{u}_{0:T-1}}} \quad J(\boldsymbol{\theta}^*)= \textstyle\sum\nolimits_{t=0}^{T-1}c(\boldsymbol{x}_t,\boldsymbol{u}_t,\boldsymbol{\theta}^*)+h(\boldsymbol{x}_T,\boldsymbol{\theta}^*),\\
% &\textrm{subject to} \quad \boldsymbol x_{t+1} = \boldsymbol{f}(\boldsymbol{x}_t, \boldsymbol{u}_t,\boldsymbol{\theta}^*), \quad \text{with}\ \boldsymbol{x}_0\ \text{given}. \\
% \end{aligned}
% \end{equation}
\begin{argmini!}|s|
{ \boldsymbol{x}_{1:T},\boldsymbol{u}_{0:T-1} }{ J(\boldsymbol{x}_{0:T},\boldsymbol{u}_{0:T-1},\boldsymbol{\theta}^*)= \textstyle\sum\nolimits_{t=0}^{T-1}c(\boldsymbol{x}_t,\boldsymbol{u}_t,\boldsymbol{\theta}^*)+h(\boldsymbol{x}_T,\boldsymbol{\theta}^*) \label{eq:tru_obj}}
{\label{eq:optiOC}}{ \{\boldsymbol{x}_{1:T}(\boldsymbol{\theta}^*),\boldsymbol{u}_{0:T-1}(\boldsymbol{\theta}^*)\} = }
\addConstraint{ \boldsymbol x_{t+1} = \boldsymbol{f}(\boldsymbol{x}_t, \boldsymbol{u}_t,\boldsymbol{\theta}^*), \quad \text{with}\ \boldsymbol{x}_0\ \text{given}. \label{eq:tru_dyn}}
\end{argmini!}
where $t=0,1,2,\cdots, T$ is the time index with $T$ being the final time; $\boldsymbol{x}_t\in\mathbb{R}^n$ and $\boldsymbol{u}_t \in\mathbb{R}^m$ denote the system state and control input, respectively;  $\boldsymbol{x}_{0:T}(\boldsymbol{\theta}^*)\triangleq\col\{\boldsymbol{x}_0(\boldsymbol{\theta}^*),\cdots,\boldsymbol{x}_T(\boldsymbol{\theta}^*)\}$ and $\boldsymbol{u}_{0:T-1}(\boldsymbol{\theta}^*)\triangleq\col\{\boldsymbol{u}_0(\boldsymbol{\theta}^*),\cdots,\boldsymbol{u}_{T-1}(\boldsymbol{\theta}^*)\}$ denote the states and inputs trajectory given parameter $\boldsymbol{\theta}^*$, respectively; $\boldsymbol{x}_t^*(\boldsymbol{\theta}^*)$ and $\boldsymbol{u}_t^*(\boldsymbol{\theta}^*)$ denote the state and input given $\boldsymbol{\theta}^*$ at time $t$ respectively; $\boldsymbol{f}:\mathbb{R}^n\times\mathbb{R}^m\times\mathbb{R}^p\rightarrow\mathbb{R}^n$ denotes a twice-differentiable time-invariant system dynamics; $c:\mathbb{R}^n\times\mathbb{R}^m\times\mathbb{R}^p\mapsto\mathbb{R}$ and $h:\mathbb{R}^n\times\mathbb{R}^p\mapsto\mathbb{R}$ denote running cost the final cost, respectively, both of which are assumed to be twice-differentiable.

\begin{remark}
Including the parameter $\theta^*$ in the system dynamics allows for the representation of both partially known and completely unknown dynamics. For partially known dynamics, it is parameterized via a known physical dynamic model with unknown physical parameters. For example, this could be a quadrotor dynamics with known structure and unknown inertia and mass \citep{wang2014robust,jin2020pontryagin,revach2022kalman}. In the case of completely unknown dynamics, parameterization is done by neural networks. In this case, the neural network captures the evolution of the state, where the parameter $\boldsymbol{\theta}^*$ represents the weights and biases of the neural network \citep{kumpati1990identification,lewis1998neural,nelles2020nonlinear}.
\end{remark}

% For a choice of $\boldsymbol{\theta}$, the trajectory of optimal control system (\ref{eq:OCsys}) can be determined by solving the optimal control problem:
% \begin{equation} \label{eq:optiOC}
% \begin{aligned}
% \{\boldsymbol{x}_{0:T}(\boldsymbol{\theta}),\boldsymbol{u}_{0:T-1}(\boldsymbol{\theta})\}\in&\arg\min_{\substack{\boldsymbol{x}_{0:T}\\\boldsymbol{u}_{0:T-1}}} \quad J(\boldsymbol{\theta})\\
% &\textrm{subject to} \quad \boldsymbol x_{t+1} =\boldsymbol f(\boldsymbol x_t, \boldsymbol u_{t},\boldsymbol{\theta}), \ \forall t \text{ with given } \boldsymbol{x}_0. \\
% \end{aligned}
% \end{equation}
% \begin{equation} \label{eq:optiOC}
% \{\boldsymbol{x}_{0:T}(\boldsymbol{\theta}),\boldsymbol{u}_{0:T-1}(\boldsymbol{\theta})\}\in\arg\min_{{\boldsymbol{x}_{0:T},\boldsymbol{u}_{0:T-1}}} \  J(\boldsymbol{\theta}) \ \ \text{s.t.} \ \  \boldsymbol x_{t+1} =\boldsymbol f(\boldsymbol x_t, \boldsymbol u_{t},\boldsymbol{\theta}), \ \forall t \text{ with given } \boldsymbol{x}_0.
% \end{equation}
% where $\boldsymbol{x}_{0:T}(\boldsymbol{\theta})\triangleq\col\{\boldsymbol{x}_0(\boldsymbol{\theta}),\cdots,\boldsymbol{x}_T(\boldsymbol{\theta})\}$ and $\boldsymbol{u}_{0:T-1}(\boldsymbol{\theta})\triangleq\col\{\boldsymbol{u}_0(\boldsymbol{\theta}),\cdots,\boldsymbol{u}_{T-1}(\boldsymbol{\theta})\}$ denote the states and inputs trajectory given parameter $\boldsymbol{\theta}$, respectively.
For notation simplicity, we define the unknown trajectory of the optimal control system $\Sigma(\boldsymbol{\theta}^*)$ as
\begin{equation}
    \boldsymbol{\xi}({\boldsymbol{\theta}^*})\triangleq
\text{col}\{ \boldsymbol{x}_{0:T}(\boldsymbol{\theta}^*),
\boldsymbol{u}_{0:T-1}(\boldsymbol{\theta}^*) \}
\in\mathbb{R}^{(T+1)n+Tm}
\end{equation}
% $\boldsymbol{\xi}({\boldsymbol{\theta}})\triangleq
% \text{col}\{ \boldsymbol{x}_{0:T}(\boldsymbol{\theta}),
% \boldsymbol{u}_{0:T-1}(\boldsymbol{\theta}) \}
% \in\mathbb{R}^{(T+1)n+Tm}$.

The goal of this paper is to estimate $\boldsymbol{\theta}^*$. Define $\hat{\boldsymbol{\theta}}\in\mathbb{R}^p$ as an arbitrary estimation of $\boldsymbol{\theta}^*$.
Then for estimation purposes, a copy, $\Sigma(\hat{\boldsymbol{\theta}})$, of the autonomous system $\Sigma(\boldsymbol{\theta}^*)$ can be proposed by replacing $\boldsymbol{\theta}^*$ with $\hat{\boldsymbol{\theta}}$ in \mbox{(\ref{eq:optiOC})}, i.e.,
\begin{argmini!}|s|
{ \boldsymbol{x}_{1:T},\boldsymbol{u}_{0:T-1} }{ J(\boldsymbol{x}_{0:T},\boldsymbol{u}_{0:T-1},\hat{\boldsymbol{\theta}})= \textstyle\sum\nolimits_{t=0}^{T-1}c(\boldsymbol{x}_t,\boldsymbol{u}_t,\hat{\boldsymbol{\theta}})+h(\boldsymbol{x}_T,\hat{\boldsymbol{\theta}}) \label{eq:est_obj}}
{\label{eq:est_auto_sys}}{ \{\boldsymbol{x}_{1:T}(\hat{\boldsymbol{\theta}}), \boldsymbol{u}_{0:T-1}(\hat{\boldsymbol{\theta}})\} = }
\addConstraint{ \boldsymbol{x}_{t+1} = \boldsymbol{f}(\boldsymbol{x}_t, \boldsymbol{u}_t, \hat{\boldsymbol{\theta}}), \  \text{with}\ \boldsymbol{x}_0\ \text{given}. \label{eq:est_dyn}}
\end{argmini!}

At each time $t$, a noisy measurement $\boldsymbol{O}_t\in\mathbb{R}^r$ is observed, where
\begin{equation} \label{eq:measFunc}
\boldsymbol{O}_t=\boldsymbol{h}(\boldsymbol{\xi}_t(\boldsymbol{\theta}^*))+\boldsymbol{v}_t.
\end{equation}
Here, $\boldsymbol{h}:\mathbb{R}^{n+m}\mapsto\mathbb{R}^r$ denotes a twice-differentiable measurement function; $\boldsymbol{\xi}_t({\boldsymbol{\theta}}^*)=\{\boldsymbol{x}_t^*(\boldsymbol{\theta}^*), \boldsymbol{u}_t^*(\boldsymbol{\theta}^*)\}$; $\boldsymbol{v}_t \sim \mathcal{N}(\boldsymbol{0}_{r},\boldsymbol{R}_t)$ denotes the measurement noise which is a multivariate Gaussian, with $\boldsymbol{R}_t\in\mathbb{R}^{r\times r}$ being the covariance matrices of the measurement noise.

With the measurement equation (\ref{eq:measFunc}) defined, this paper considers a signed residual function:
\begin{equation}
    \boldsymbol{l}(\boldsymbol{\xi}_t(\hat{\boldsymbol{\theta}}),\boldsymbol{O}_t)=\boldsymbol{O}_t - \boldsymbol{h}(\boldsymbol{\xi}_t(\hat{\boldsymbol{\theta}}))\in\mathbb{R}^r.
\end{equation}
Then, the performance of the entire trajectory can be evaluated by a cumulative loss which is assumed to be twice-differentiable:
\begin{equation} \label{eq:cumuLoss}
    L(\boldsymbol{\xi}(\hat{\boldsymbol{\theta}})) = \textstyle\sum_{t=0}^{T} \|\boldsymbol{l}(\boldsymbol{\xi}_t(\hat{\boldsymbol{\theta}}),\boldsymbol{O}_t)\|^2.
\end{equation}

The \textit{problem of interest} is to develop an online method to update the estimation $\hat{\boldsymbol{\theta}}_t\in\mathbb{R}^p$ of $\boldsymbol{\theta}^*$ at every time $t$, such that its trajectory $\boldsymbol{\xi}(\hat{\boldsymbol{\theta}}_t)$ from (\ref{eq:optiOC}) minimizes a task-specific cumulative loss $ L(\boldsymbol{\xi}(\hat{\boldsymbol{\theta}}))$.
% \begin{equation} \label{eq:generalOpti0}
% \min_{\boldsymbol{\theta}} \quad  L(\boldsymbol{\xi}(\boldsymbol{\theta})) \quad \textrm{subject to} \quad \boldsymbol{\xi}(\boldsymbol{\theta}) \textrm{ be the trajectory of \ref{eq:optiOC}}.
% \end{equation}

To achieve a specific learning or control task, one needs to select the most suitable measurement $\boldsymbol{O}_t$.  Below, we will present three modes of the proposed OCIL framework. It is worth noting that in different applications, adjustments to the configuration of system $\Sigma(\hat{\boldsymbol{\theta}})$ are required according to the task.

\textbf{Online SysID:} For a SysID problem, the goal is to identify the dynamics model of a physical system from the state-input trajectory $\boldsymbol{\xi}^o=\{\boldsymbol{x}^o_{0:T},\boldsymbol{u}_{0:T-1}\}$, where the superscript $o$ denotes the observed trajectory. The trajectory is often generated by persistent excitation of the system without considering any control objectives \mbox{\citep{keesman2011system}}. Therefore, we can set $J(\boldsymbol{x}_{0:T},\boldsymbol{u}_{0:T-1},\hat{\boldsymbol{\theta}})=0$:
\begin{equation} 
    \Sigma(\hat{\boldsymbol{\theta}}): 
\begin{aligned} \label{sys:sysID}
\text{dynamics:}\quad&& &\boldsymbol{x}_{t+1} = \boldsymbol{f}(\boldsymbol{x}_t, \boldsymbol{u}_t,\hat{\boldsymbol{\theta}}), \quad \text{with}\ \boldsymbol{x}_0\ \text{given},\\
\text{objective:}\quad&& &J(\boldsymbol{x}_{0:T},\boldsymbol{u}_{0:T-1},\hat{\boldsymbol{\theta}}) = 0.
\end{aligned}
\end{equation}
% \begin{equation} \label{sys:sysID}
%     \Sigma(\boldsymbol{\theta}):  \quad \text{dynamics:}\quad\boldsymbol{x}_{t+1} = \boldsymbol{f}(\boldsymbol{x}_t, \boldsymbol{u}_t,\boldsymbol{\theta}), \text{with}\ \boldsymbol{x}_0\ \text{given}; \quad  \text{objective:}\quad J(\boldsymbol{\theta})= 0.
% \end{equation}
To identify the model dynamics, namely finding the $\boldsymbol{\theta}^*$ in the dynamics $\boldsymbol{f}(\boldsymbol{x}_t, \boldsymbol{u}_t,\boldsymbol{\theta}^*)$, one could design the signed residual function to represent the discrepancy between the observed trajectory and the trajectory produced by $\hat{\boldsymbol{\theta}}$, i.e. $\boldsymbol{l}(\boldsymbol{\xi}_t(\hat{\boldsymbol{\theta}}),\boldsymbol{\xi}_t^o)=\boldsymbol{\xi}_t^o-\boldsymbol{\xi}_t(\hat{\boldsymbol{\theta}})$, where $\boldsymbol{\xi}_t^o$ is a slice of $\boldsymbol{\xi}^o$ at time $t$. In the SysID mode, the measurement $\boldsymbol{O}_t$ received at time $t$ is a slice of the trajectory of a physical system $\boldsymbol{\xi}_t^o$.

\textbf{Online Imitation Learning:} The objective function and the model dynamics are parameterized by an unknown $\boldsymbol{\theta}^*$. The OC system follows (\ref{eq:est_auto_sys}).
% \begin{equation}
%     \Sigma(\boldsymbol{\theta}): 
% \begin{aligned} 
%     \text{dynamics:}\quad&& \boldsymbol{x}_{t+1} &= \boldsymbol{f}(\boldsymbol{x}_t, \boldsymbol{u}_t), \quad \text{with}\ \boldsymbol{x}_0\ \text{given},\\
%     \text{objective:}\quad&& J(\boldsymbol{\theta}) &= \textstyle\sum\nolimits_{t=0}^{T-1}c(\boldsymbol{x}_t,\boldsymbol{u}_t,\boldsymbol{\theta})+h(\boldsymbol{x}_T,\boldsymbol{\theta}),
% \end{aligned}
% \end{equation}
Suppose one can observe the measurement of the expert demonstration $\boldsymbol{y}_t^*$ at each time $t$. Then, the signed residual function can be designed as $\boldsymbol{l}(\boldsymbol{\xi}_t(\hat{\boldsymbol{\theta}}),\boldsymbol{y}_t^*)=\boldsymbol{y}_t^*-\boldsymbol{g}(\boldsymbol{x}_t(\hat{\boldsymbol{\theta}}),\boldsymbol{u}_t(\hat{\boldsymbol{\theta}}))$. In this case, the measurement $\boldsymbol{O}_t$ received at time $t$ is the expert demonstration $\boldsymbol{y}^*_t$. The optimal demonstration can vary between being continuous or sparse, depending on practical application scenarios.

\textbf{Tuning Policy On-the-fly:} For an autonomous system, one would like to obtain a control policy such that the trajectory minimizes certain task loss. This mode considers a feedback controller which is parameterized by $\hat{\boldsymbol{\theta}}$, i.e. $\boldsymbol{u}_t=\boldsymbol{\mu}(\boldsymbol{x}_t,\hat{\boldsymbol{\theta}})$. Then the OC system is written as follows:
\begin{equation}
    \Sigma(\hat{\boldsymbol{\theta}}): 
\begin{aligned} \label{sys:Control}
    \text{dynamics:}\quad&& &\boldsymbol{x}_{t+1} = \boldsymbol{f}(\boldsymbol{x}_t, \boldsymbol{\mu}(\boldsymbol{x}_t,\hat{\boldsymbol{\theta}})), \quad \text{with}\ \boldsymbol{x}_0\ \text{given},\\
    \text{objective:}\quad&& &J(\boldsymbol{x}_{0:T},\boldsymbol{u}_{0:T-1},\hat{\boldsymbol{\theta}}) = \textstyle\sum\nolimits_{t=0}^{T-1}c(\boldsymbol{x}_t,\boldsymbol{\mu}(\boldsymbol{x}_t,\hat{\boldsymbol{\theta}}))+h(\boldsymbol{x}_T).
\end{aligned}
\end{equation}
Then we can design the signed residual function such that it represents trajectory tracking. For instance, the signed residual function could be $\boldsymbol{l}(\boldsymbol{\xi}_t(\hat{\boldsymbol{\theta}}),\boldsymbol{\xi}_t^d)=\boldsymbol{\xi}_t^d-\boldsymbol{\xi}_t(\hat{\boldsymbol{\theta}})$, where $\boldsymbol{\xi}_t^d$ is a slice of desired trajectory to track at time $t$.

\section{Main Results}

% This section introduces an online parameter estimator based on the extended Kalman filter (EKF). Going forward, we will show the challenge of obtaining the Kalman gain. To tackle this challenge, we employ a gradient generator (GG) based on the Pontryagin differential programming (PDP). Then the proposed OCIL framework will be introduced and supported with theoretical analysis.

The proposed OCIL consists of two main components, both of which are inspired by control theory. Specifically, OCIL first proposes an online parameter estimator based on the extended Kalman filter (EKF).
Going forward, we will show the challenge of obtaining the Kalman gain. To tackle this challenge, the gradient information for the loss with respect to the tunable parameter is required. Therefore, OCIL employs a gradient generator (GG) based on Pontryagin Differential Programming to calculate the exact gradient. Then the proposed OCIL framework will be introduced and supported with theoretical analysis.

\subsection{Online Parameter Estimator}

To minimize the cumulative task loss $L(\boldsymbol{\xi}(\hat{\boldsymbol{\theta}}))$ with measurement $\boldsymbol{O}_t$, which is unavailable until time $t$, the optimization problem that needs to be solved in an online fashion is:
\begin{equation} \label{eq:generalOpti}
\min_{\boldsymbol{\theta}} \quad \textstyle\sum_{t=0}^{T} \|\boldsymbol{l}(\boldsymbol{\xi}_t(\hat{\boldsymbol{\theta}}), \boldsymbol{O}_t)\|^2\quad\textrm{subject to} \quad \boldsymbol{\xi}(\hat{\boldsymbol{\theta}}) \textrm{ is the trajectory of } (\ref{eq:est_auto_sys}).
\end{equation}
% A method is desired to incrementally find $\boldsymbol{\theta}$ in real time because $\boldsymbol{O}_t$ is unavailable until time $t$.
The optimization problem (\ref{eq:generalOpti}) is essentially a least squares problem, although under constraints. One of the most famous methods to solve the least squares problems incrementally is the EKF \citep{1996Bertsekas,ribeiro2004kalman}.
The EKF was proposed to incrementally estimate the state of a system using measured output available at each time step. In our problem setting, instead of estimating the state of a system, our goal is to estimate the parameter $\boldsymbol{\theta}^*$ by utilizing the measurement $\boldsymbol{O}_t$ that is available at each time $t$. Therefore, by considering the parameter $\boldsymbol{\theta}^*$ as the state to be estimated, one can introduce a new dynamical system:
\begin{equation} \label{eq:dynTheta}
\text{dynamics: } \boldsymbol{\theta}_{t+1} =  \boldsymbol{\theta}_{t}, \text{ with }\boldsymbol{\theta}_0 = \boldsymbol{\theta}^*, \quad
\text{measurement: } \boldsymbol{O}_t=\boldsymbol{h}(\boldsymbol{\xi}_t(\boldsymbol{\theta}_t))+\boldsymbol{v}_t,
\end{equation}
 The online estimation of $\boldsymbol{\theta}^*$ via EKF can be done as follows \citep{ribeiro2004kalman}:
\begin{subequations} \label{eq:main_EKF}
\begin{align}
&\hat{\boldsymbol{\theta}}^{-}_{t} := \hat{\boldsymbol{\theta}}_{t-1}, \ \boldsymbol{P}^{-}_{t} := \boldsymbol{P}_{t-1} \label{eq:EKFpredict} \\
&\boldsymbol{K}_{t} := \boldsymbol{P}^{-}_t\boldsymbol{L}_{t}^\prime(\boldsymbol{L}_t\boldsymbol{P}^{-}_{t}\boldsymbol{L}^\prime_t+\boldsymbol{R}_t)^{-1}, \ 
\boldsymbol{P}_{t} := (\boldsymbol{I}_p-\boldsymbol{K}_t\boldsymbol{L}_t)\boldsymbol{P}^{-}_{t}, \ 
\hat{\boldsymbol{\theta}}_{t} := \hat{\boldsymbol{\theta}}^{-}_t + \boldsymbol{K}_{t}(\boldsymbol{O}_t - \boldsymbol{h}(\boldsymbol{\xi}_t(\hat{\boldsymbol{\theta}}^{-}_t))) , \label{eq:EKFupdate}
\end{align}
\end{subequations}
\begin{equation}
    \boldsymbol{L}_t \triangleq \frac{d \boldsymbol{l}(\boldsymbol{\xi}_t(\boldsymbol{\theta}_t), \boldsymbol{O}_t)}{d \boldsymbol{\theta}_t}|_{\boldsymbol{\theta}_t = \hat{\boldsymbol{\theta}}^{-}_{t}}\in\mathbb{R}^{r\times p}
\end{equation}
where (\ref{eq:EKFpredict}) predicts the dynamics; (\ref{eq:EKFupdate}) updates the parameter estimate.
Here, the superscript $^-$ means the term is not yet updated by measurement residual; $\boldsymbol{P}_t\in\mathbb{R}^{p\times p}$ is a positive-definite matrix that denotes the covariance of the estimate; $\boldsymbol{K}_t\in\mathbb{R}^{p\times r}$ denotes the Kalman gain. 
Throughout the estimation process, all of the terms are known except $\boldsymbol{L}_t$. It is challenging to obtain this term as the signed residual function $\boldsymbol{l}(\boldsymbol{\xi}_t(\hat{\boldsymbol{\theta}}), \boldsymbol{O}_t)$ is not an explicit function of $\boldsymbol{\theta}$. In the next subsection, we will present a \textit{gradient generator} which computes the exact value for $\boldsymbol{L}_t$.

% Predict:
% \begin{align}
% \hat{\boldsymbol{\theta}}^{-}_{t} &= \hat{\boldsymbol{\theta}}_{t-1}, \label{eq:EKFpredict} \\
% \boldsymbol{P}^{-}_{t}&=\boldsymbol{P}_{t-1}+\boldsymbol{Q}_{t-1},
% \end{align}

% Update:
% \begin{align}
% \boldsymbol{K}_{t}&=\boldsymbol{P}^{-}_t\boldsymbol{L}_{t}^\prime(\boldsymbol{L}_t\boldsymbol{P}^{-}_{t}\boldsymbol{L}^\prime_t+R_t)^{-1}, \label{eq:Kalmangain}\\ 
% \boldsymbol{P}_{t}&=(\boldsymbol{I}_p-\boldsymbol{K}_t\boldsymbol{L}_t)\boldsymbol{P}^{-}_{t}, \label{eq:Kalmancov}\\
%     \begin{split}
% \hat{\boldsymbol{\theta}}_{t}&=\hat{\boldsymbol{\theta}}^{-}_t + \boldsymbol{K}_{t}l(\boldsymbol{\xi}_t(\hat{\boldsymbol{\theta}}^{-}_t), \boldsymbol{O}_t), \label{eq:EKFupdate}
%      \end{split}
% \end{align}
% % with
% \begin{equation}
% \boldsymbol{L}_t \triangleq \frac{\partial l(\boldsymbol{\xi}_t(\boldsymbol{\theta}), \boldsymbol{O}_t)}{\partial \boldsymbol{\theta}}|_{\boldsymbol{\theta} = \hat{\boldsymbol{\theta}}^{-}_{t}}\ \in\mathbb{R}^{1\times p},
% \end{equation}

% In the convergence analysis in \ref{sec:convergenveAnalysis}, we will show how the covariance matrices $\boldsymbol{Q}_t$ and $\boldsymbol{R}_t$ affect the convergence of the proposed estimator.

\subsection{Gradient Generator}

In this section, for notation simplicity, the parameter estimate $\hat{\boldsymbol{\theta}}_t^-$ is simplified to $\boldsymbol{\theta}$; $\frac{d \boldsymbol{l}(\boldsymbol{\xi}_t(\boldsymbol{\theta}_t), \boldsymbol{O}_t)}{d \boldsymbol{\theta}_t}|_{\boldsymbol{\theta}_t = \hat{\boldsymbol{\theta}}^{-}_{t}}$ is written as $\frac{d \boldsymbol{l}(\boldsymbol{\xi}_t(\boldsymbol{\theta}))}{d \boldsymbol{\theta}}$. To obtain the gradient $\frac{d \boldsymbol{l}(\boldsymbol{\xi}_t(\boldsymbol{\theta}))}{d \boldsymbol{\theta}}$, one can employ the chain rule by definition,
% i.e. $\scriptsize{\displaystyle\frac{d \boldsymbol{l}(\boldsymbol{\xi}_t(\boldsymbol{\theta}))}{d \boldsymbol{\theta}}=\frac{\partial \boldsymbol{l}(\boldsymbol{\xi}_t(\boldsymbol{\theta}))}{\partial \boldsymbol{\xi}_t({\boldsymbol{\theta}})}\frac{\partial  \boldsymbol{\xi}_t({\boldsymbol{\theta}})}{\partial \boldsymbol{\theta}}}$,
\begin{equation} \label{eq:chainRule}
   \frac{d \boldsymbol{l}(\boldsymbol{\xi}_t(\boldsymbol{\theta}))}{d \boldsymbol{\theta}}=\frac{\partial \boldsymbol{l}(\boldsymbol{\xi}_t(\boldsymbol{\theta}))}{\partial \boldsymbol{\xi}_t({\boldsymbol{\theta}})}\frac{\partial  \boldsymbol{\xi}_t({\boldsymbol{\theta}})}{\partial \boldsymbol{\theta}},
\end{equation}
where $\frac{\partial \boldsymbol{l}(\boldsymbol{\xi}_t(\boldsymbol{\theta}))}{\partial \boldsymbol{\xi}_t({\boldsymbol{\theta}})}$ is known since the signed residual function is pre-designed. The challenge that remains is to find the partial derivative $\frac{\partial  \boldsymbol{\xi}_t({\boldsymbol{\theta}})}{\partial \boldsymbol{\theta}}$, i.e. an analytical relation between trajectory $\boldsymbol{\xi}_t$ and the tunable parameter $\boldsymbol{\theta}$. To tackle this challenge, the gradient generator in \cite{jin2020pontryagin} is used to obtain the exact value of $\frac{\partial  \boldsymbol{\xi}_t({\boldsymbol{\theta}})}{\partial \boldsymbol{\theta}}$.
% which is summarized in Algorithm~\ref{algorithm:GG}.

Given the OC system (\ref{eq:est_auto_sys}), one can obtain the Hamiltonian equation 
\begin{equation}
    H_{t} = c(\boldsymbol{x}_{t},\boldsymbol{u}_{t},\boldsymbol{\theta})+\boldsymbol{f}(\boldsymbol{x}_{t},\boldsymbol{u}_{t},\boldsymbol{\theta})^\prime\boldsymbol{\lambda}_{t+1}
\end{equation}
for all $t=0,\cdots,T-1$, where $\boldsymbol{\lambda}_{t}\in\mathbb{R}^n$ denotes the Lagrangian multiplier associated with the equality constraint of model dynamics. With the definition of $\boldsymbol{\xi}({\boldsymbol{\theta}})$, one has $\frac{\partial \boldsymbol{\xi}({\boldsymbol{\theta}})}{\partial \boldsymbol{\theta}}=\text{col}\{\frac{\partial\boldsymbol{x}_{1:T}(\boldsymbol{\theta})}{\partial\boldsymbol{\theta}},\frac{\partial\boldsymbol{u}_{0:T-1}(\boldsymbol{\theta})}{\partial\boldsymbol{\theta}}\}$.
% \begin{equation}
% \frac{\partial \boldsymbol{\xi}({\boldsymbol{\theta}})}{\partial \boldsymbol{\theta}}=\text{col}\{ \frac{\partial\boldsymbol{x}_{0:T}(\boldsymbol{\theta})}{\partial\boldsymbol{\theta}}, \frac{\partial\boldsymbol{u}_{0:T-1}(\boldsymbol{\theta})}{\partial\boldsymbol{\theta}} \}.
% \end{equation}
% \begin{equation}
% \frac{\partial \boldsymbol{\xi}({\boldsymbol{\theta}})}{\partial \boldsymbol{\theta}}=\begin{bmatrix}
%     \frac{\partial\boldsymbol{x}_{0:T}(\boldsymbol{\theta})}{\partial\boldsymbol{\theta}}\\
%     \frac{\partial\boldsymbol{u}_{0:T-1}(\boldsymbol{\theta})}{\partial\boldsymbol{\theta}}
% \end{bmatrix}.
% \end{equation}
By defining 
% $\boldsymbol{X}_{t}\triangleq$$\frac{\partial\boldsymbol{x}_{t}({\boldsymbol{\theta}})}{\partial\boldsymbol{\theta}}$$\in\mathbb{R}^{n\times p}$, $\boldsymbol{U}_{t}\triangleq$$\frac{\partial\boldsymbol{u}_{t}({\boldsymbol{\theta}})}{\partial\boldsymbol{\theta}}$$\in\mathbb{R}^{m\times p}$,
\begin{equation}
X_{t}\triangleq\frac{\partial\boldsymbol{x}_{t}({\boldsymbol{\theta}})}{\partial\boldsymbol{\theta}}\in\mathbb{R}^{n\times p}, \quad U_{t}\triangleq\frac{\partial\boldsymbol{u}_{t}({\boldsymbol{\theta}})}{\partial\boldsymbol{\theta}}\in\mathbb{R}^{m\times p},
\end{equation}
one can utilize the following lemma from \cite{jin2020pontryagin} to obtain the partial derivatives $\displaystyle \frac{\partial  \boldsymbol{\xi}_t({\boldsymbol{\theta}})}{\partial \boldsymbol{\theta}}$:

\begin{lemma} \label{lemma:GG}
\cite{jin2020pontryagin} Define the Jacobian and Hessian matrices related to $\boldsymbol{\xi}(\boldsymbol{\theta})$ as:
\begin{equation} \label{eq:auxTerms}
% \small
\begin{aligned}
&\boldsymbol{F}_{t}=\frac{\partial\boldsymbol{f}}{\partial\boldsymbol{x}_t},\  \boldsymbol{G}_{t}=\frac{\partial\boldsymbol{f}}{\partial\boldsymbol{u}_t},\  \boldsymbol{E}_{t}=\frac{\partial\boldsymbol{f}}{\partial\boldsymbol{\theta}},\boldsymbol{H}_{t}^{xx}=\frac{\partial^2H_{t}}{\partial\boldsymbol{x}_{t}\partial\boldsymbol{x}_{t}},\  \boldsymbol{H}_{t}^{xu}=\frac{\partial^2H_{t}}{\partial\boldsymbol{x}_{t}\partial\boldsymbol{u}_{t}}=(\boldsymbol{H}_{t}^{ux})^\prime,\\
& \boldsymbol{H}_{t}^{uu}=\frac{\partial^2H_{t}}{\partial\boldsymbol{u}_{t}\partial\boldsymbol{u}_{t}},
\boldsymbol{H}_{t}^{x\theta}=\frac{\partial^2H_{t}}{\partial\boldsymbol{x}_{t}\partial\boldsymbol{\theta}},\  \boldsymbol{H}_{t}^{u\theta}=\frac{\partial^2H_{t}}{\partial\boldsymbol{u}_{t}\partial\boldsymbol{\theta}},\ \boldsymbol{H}_{T}^{xx}=\frac{\partial^2h}{\partial\boldsymbol{x}_{T}\partial\boldsymbol{x}_{T}},\  \boldsymbol{H}_{T}^{x\theta}=\frac{\partial^2h}{\partial\boldsymbol{x}_{T}\partial\boldsymbol{\theta}}.
\end{aligned}
\end{equation}

If $\boldsymbol{H}_{t}^{uu}$ is invertible for all $t=0,\cdots,T-1$, the following recursions from $t=T$ to $t=0$ hold:
% \begin{subequations}
% \begin{align}
% V_{t}&=C_{t}+A_{t}^\prime(I+V_{t+1}B_{t})^{-1}V_{t+1}A_{t}, \label{eq:Pt} \\
% W_{t}&=A_{t}^\prime(I+V_{t+1}B_{t})^{-1}(W_{t+1}+V_{t+1}M_{t})+N_{t}, \label{eq:Wt}
% \end{align}
% \end{subequations}
\begin{equation} \label{eq:backwardPass}
% \small 
\begin{aligned}
&\boldsymbol{V}_{t}=\boldsymbol{C}_{t}+\boldsymbol{A}_{t}^\prime(\boldsymbol{I}+\boldsymbol{V}_{t+1}\boldsymbol{B}_{t})^{-1}\boldsymbol{V}_{t+1}\boldsymbol{A}_{t}, \\&\boldsymbol{W}_{t}=\boldsymbol{A}_{t}^\prime(\boldsymbol{I}+\boldsymbol{V}_{t+1}\boldsymbol{B}_{t})^{-1}(\boldsymbol{W}_{t+1}+\boldsymbol{V}_{t+1}\boldsymbol{M}_{t})+\boldsymbol{N}_{t},
\end{aligned}
\end{equation}
%     \begin{equation} \label{eq:Pt}
%     V_{t}=C_{t}+A_{t}^\prime(I+V_{t+1}B_{t})^{-1}V_{t+1}A_{t},
%     \end{equation}
%     \begin{equation} \label{eq:Wt}
% W_{t}=A_{t}^\prime(I+V_{t+1}B_{t})^{-1}(W_{t+1}+V_{t+1}M_{t})+N_{t}, 
%     \end{equation}
with $\boldsymbol{V}_{T}=\boldsymbol{H}^{xx}_{T}$ and $\boldsymbol{W}_{T}=\boldsymbol{H}^{x\theta}_{T}$. Here, $\boldsymbol{A}_{t}=\boldsymbol{F}_{t}-\boldsymbol{G}_{t}(\boldsymbol{H}^{uu}_{t})^{-1}\boldsymbol{H}_{t}^{ux}$, $\boldsymbol{B}_{t}=\boldsymbol{G}_{t}(\boldsymbol{H}^{uu}_{t})^{-1}\boldsymbol{G}_{t}^\prime$, $\boldsymbol{M}_{t}=\boldsymbol{E}_{t}-\boldsymbol{G}_{t}(H^{uu}_{t})^\prime \boldsymbol{H}^{u\theta}_{t}$, $\boldsymbol{C}_{t}=\boldsymbol{H}^{xx}_{t}-\boldsymbol{H}^{xu}_{t}(\boldsymbol{H}^{uu}_{t})^{-1}\boldsymbol{H}^{ux}_{t}$, $\boldsymbol{N}_{t}=\boldsymbol{H}^{x\theta}_{t}-\boldsymbol{H}^{xu}_{t}(\boldsymbol{H}^{uu}_{t})^\prime \boldsymbol{H}^{u\theta}_{t}$ are all known given (\ref{eq:auxTerms}). Then, the partial derivative $\displaystyle\frac{\partial \boldsymbol{\xi}({\boldsymbol{\theta}})}{\partial \boldsymbol{\theta}}$ can be obtained by recursively solving the following equations from $t=0$ to $T-1$ with $\boldsymbol{X}_{0}(\boldsymbol{\theta})=\boldsymbol{0}$:
\begin{equation} \label{eq:recursiveGG}
% \small 
\begin{aligned}
&\boldsymbol{U}_{t}=-(\boldsymbol{H}^{uu}_{t})^{\text{-}1}(\boldsymbol{H}^{ux}_{t}\boldsymbol{X}_{t}+\boldsymbol{H}^{u\theta}_{t}+\boldsymbol{G}_{t}^\prime(\boldsymbol{I}+\boldsymbol{V}_{t+1}\boldsymbol{B}_{t})^{\text{-}1}(\boldsymbol{V}_{t+1}\boldsymbol{A}_{t}\boldsymbol{X}_{t}+\boldsymbol{V}_{t+1}\boldsymbol{M}_{t}+\boldsymbol{W}_{t+1})), \\
&\boldsymbol{X}_{t+1}=\boldsymbol{F}_{t}\boldsymbol{X}_{t}+\boldsymbol{G}_{t}\boldsymbol{U}_{t}+\boldsymbol{E}_{t}.
\end{aligned}
\end{equation}
\end{lemma}

The terms in (\ref{eq:auxTerms}) are based on the trajectory $\boldsymbol{\xi}({\boldsymbol{\theta}})$ and the associated Lagrangian multiplier $\boldsymbol{\lambda}_{0:T-1}$. According to the discrete-time Pontryagin Maximum Principle \citep{jin2020pontryagin}, the trajectory of the Lagrangian multiplier can be obtained by 
% $\boldsymbol{\lambda}_{T}=\frac{\partial h}{\partial \boldsymbol{x}_{T}}$, $
% \boldsymbol{\lambda}_{t}\triangleq\frac{\partial H_{t}}{\partial \boldsymbol{x}_{t}}=\frac{\partial c}{\partial \boldsymbol{x}_{t}}+\frac{\partial h}{\partial \boldsymbol{x}_{t}}\boldsymbol{\lambda}_{t+1}$, for $t=T-1,\cdots,1$.
\begin{equation} \label{eq:lagrangianMulti}
% \small
    \boldsymbol{\lambda}_{T}=\frac{\partial h}{\partial \boldsymbol{x}_{T}},\quad
    \boldsymbol{\lambda}_{t}\triangleq\frac{\partial H_{t}}{\partial \boldsymbol{x}_{t}}=\frac{\partial c}{\partial \boldsymbol{x}_{t}}+\frac{\partial h}{\partial \boldsymbol{x}_{t}}\boldsymbol{\lambda}_{t+1}, \text{ for } t=T-1,\cdots,1.
\end{equation}
\begin{remark}
    Lemma~\ref{lemma:GG} proposes a recursive way to obtain the exact gradient of the trajectory $\boldsymbol{\xi}(\boldsymbol{\theta})$ with respect to the parameter $\boldsymbol{\theta}$, i.e. $\frac{\partial \boldsymbol{\xi}({\boldsymbol{\theta}})}{\partial \boldsymbol{\theta}}$.
\end{remark}

\subsection{OCIL Framework}

\begin{figure}[!htbp]
\center
\includegraphics[width=0.95\columnwidth]{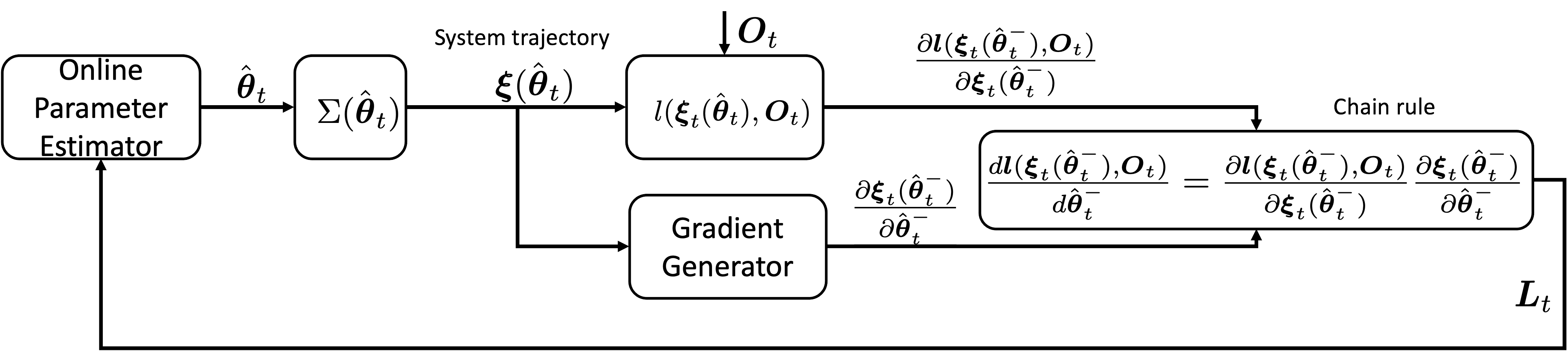}
\caption{Framework of Online Control-Informed Learning.}
\label{fig:framework}
\end{figure}

With the online parameter estimator and the gradient generator, we propose the Online Control-Informed Learning framework in Fig. \ref{fig:framework}. The framework is summarized in Algorithm \ref{algorithm:OPDP}.

As shown in Fig. \ref{fig:framework}, at each time step, the predefined OC system $\Sigma(\hat{\boldsymbol{\theta}}_t)$ generates a system trajectory $\boldsymbol{\xi}(\hat{\boldsymbol{\theta}}_t)$ by performing optimal control with given $\boldsymbol{x}_0$ and $\hat{\boldsymbol{\theta}}_t$. The trajectory $\boldsymbol{\xi}(\hat{\boldsymbol{\theta}}_t)$ is then fed into the signed residual function $\boldsymbol{l}(\boldsymbol{\xi}_t(\hat{\boldsymbol{\theta}}_t),\boldsymbol{O}_t)$ and the gradient generator. Along with the information $\boldsymbol{O}_t$ obtained at time $t$, the signed residual function generates $\frac{\partial \boldsymbol{l}(\boldsymbol{\xi}_t(\hat{\boldsymbol{\theta}}_t^-),\boldsymbol{O}_t)}{\partial \boldsymbol{\xi}_t(\hat{\boldsymbol{\theta}}_t^-)}$, while $\frac{\partial \boldsymbol{\xi}_t(\hat{\boldsymbol{\theta}}_t^-)}{\partial \hat{\boldsymbol{\theta}}_t^-}$ is generated by the gradient generator in Algorithm \ref{algorithm:GG}. The chain rule is then performed to obtain the Jacobian matrix $\boldsymbol{L}_t$, which is then passed into the online parameter estimator for the estimation of $\boldsymbol{\theta}^*$.

% \vspace{-1em}

% \vspace{-1em}
\begin{minipage}{0.45\textwidth}
\begin{algorithm2e}[H]
\setstretch{1.06}
\small \caption{\small Gradient Generator (GG)} \label{algorithm:GG}
\KwIn{Trajectory $\boldsymbol{\xi}(\hat{\boldsymbol{\theta}}_t^-)$ from $\Sigma(\hat{\boldsymbol{\theta}}_t^-)$}
Compute the coefficient matrices in (\ref{eq:auxTerms}) ;\\
Set $\boldsymbol{V}_T=\boldsymbol{H}^{xx}_T$ and $\boldsymbol{W}_T=\boldsymbol{H}^{x\theta}_T$;\\
\For{$t\leftarrow T$ to $0$ by $\Delta t$}{
Update $\boldsymbol{V}_t$ and $\boldsymbol{W}_t$ using (\ref{eq:backwardPass})
}
Set $\boldsymbol{X}_0(\hat{\boldsymbol{\theta}}_t^-)=\boldsymbol{0}$;\\
\For {$t\leftarrow 0$ to $T$ by $\Delta t$}{
Update $\boldsymbol{X}_t(\hat{\boldsymbol{\theta}}_t^-)$ and $\boldsymbol{U}_t(\hat{\boldsymbol{\theta}}_t^-)$ using (\ref{eq:recursiveGG})
}
\KwOut{$\frac{\partial \boldsymbol{\xi}(\hat{\boldsymbol{\theta}}_t^-)}{\partial \hat{\boldsymbol{\theta}}_t^-}=\{\boldsymbol{X}_{1:T}(\hat{\boldsymbol{\theta}}_t^-),\boldsymbol{U}_{0:T-1}(\hat{\boldsymbol{\theta}}_t^-)\}$}  
\end{algorithm2e}
\end{minipage}
\hfill
\begin{minipage}{0.5\textwidth}
\begin{algorithm2e}[H]
    % \setstretch{1}
	\small \caption{\small Online Control-Informed Learning}            \label{algorithm:OPDP}
	\SetKwInput{initialize}{Initialize}
        \SetKwInput{sysloss}{System and Residual}
        \sysloss{$\Sigma(\hat{\boldsymbol{\theta}})$ and $\boldsymbol{l}(\boldsymbol{\xi}_t(\hat{\boldsymbol{\theta}}),\boldsymbol{O}_t)$}
	\initialize{$\hat{\boldsymbol{\theta}}_0, \boldsymbol{P}_0$}
        \For{$t=t_0,t_1,\cdots$}{
	    Obtain new information $\boldsymbol{O}_t$;\\
            Solve $\boldsymbol{\xi}(\hat{\boldsymbol{\theta}}_t)$ from current OC system $\Sigma(\hat{\boldsymbol{\theta}}_t)$;\\
            Obtain $\frac{\partial \boldsymbol{\xi}_t(\hat{\boldsymbol{\theta}}^-_t)}{\partial \hat{\boldsymbol{\theta}}^-_t}$ with GG in Algorithm \ref{algorithm:GG};\\
            Obtain $\frac{\partial \boldsymbol{l}(\boldsymbol{\xi}_t(\hat{\boldsymbol{\theta}}^-_t),\boldsymbol{O}_t)}{\partial \boldsymbol{\xi}_t(\hat{\boldsymbol{\theta}}^-_t)}$ from $\boldsymbol{l}(\boldsymbol{\xi}_t(\hat{\boldsymbol{\theta}}^-_t),\boldsymbol{O}_t)$;\\
            % Obtain ${\scriptsize{\displaystyle\frac{d \boldsymbol{l}(\boldsymbol{\xi}_t(\boldsymbol{\theta}_t),\boldsymbol{O}_t)}{d \boldsymbol{\theta}_t}}}$ via chain rule \ref{eq:chainRule};\\
            Obtain $\boldsymbol{L}_t$ via the chain rule (\ref{eq:chainRule});\\
            % Obtain $\boldsymbol{L}_t$ via the chain rule;\\
            Update $\hat{\boldsymbol{\theta}}_t$ using the estimator (\ref{eq:main_EKF});
            }	
\end{algorithm2e}
\end{minipage}

\subsection{Convergence Analysis} \label{sec:convergenveAnalysis}

This subsection presents the convergence analysis of the online parameter estimator. The analysis employs a candidate Lyapunov function and introduces how the measurement covariance matrix $\boldsymbol{R}_t$ affects the convergence of the cumulative loss $L(\boldsymbol{\xi}(\hat{\boldsymbol{\theta}}))$. In this section, for brevity, the signed residual function $\boldsymbol{l}(\boldsymbol{\xi}_t(\hat{\boldsymbol{\theta}}), \boldsymbol{O}_t)$ is written as $\boldsymbol{l}(\boldsymbol{\xi}_t(\hat{\boldsymbol{\theta}}))$.
Suppose for a specific task, the optimal cumulative loss $L(\boldsymbol{\xi}(\boldsymbol{\theta}^*))=0$. Then, we define the estimation error as $\tilde{\boldsymbol{\theta}}_t = \boldsymbol{\theta}^*-\hat{\boldsymbol{\theta}}_t$
% , and prediction error as $\tilde{\boldsymbol{\theta}}^-_t = \boldsymbol{\theta}_t-\hat{\boldsymbol{\theta}}^-_t$
. Furthermore, we define
% \begin{subequations}
% \begin{align}
% \text{Estimation error: }\tilde{\boldsymbol{\theta}}_t &= \boldsymbol{\theta}_t-\hat{\boldsymbol{\theta}}_t, \label{eq:estimationError} \\
% \text{Prediction error: }\tilde{\boldsymbol{\theta}}^-_t &= \boldsymbol{\theta}_t-\hat{\boldsymbol{\theta}}^-_t, \label{eq:predictionError} \\
% \text{Measurement error: }\boldsymbol{e}_t&=\boldsymbol{l}(\boldsymbol{\theta}_t)-\boldsymbol{l}(\hat{\boldsymbol{\theta}}^{-}_t). \label{eq:observationError}
% \end{align}
% \end{subequations}
% \begin{subequations}
% \begin{align}
% &\text{Estimation error: }\tilde{\boldsymbol{\theta}}_t = \boldsymbol{\theta}_t-\hat{\boldsymbol{\theta}}_t, \quad
% \text{Prediction error: }\tilde{\boldsymbol{\theta}}^-_t = \boldsymbol{\theta}_t-\hat{\boldsymbol{\theta}}^-_t, \label{eq:thetaErrors} \\
% &\text{Measurement error: }\boldsymbol{e}_t=\boldsymbol{l}(\boldsymbol{\theta}_t)-\boldsymbol{l}(\hat{\boldsymbol{\theta}}^{-}_t). \label{eq:observationError}
% \end{align}
% \end{subequations}
\begin{equation}\label{eq:observationError}
\begin{aligned}
    &\text{Measurement error: }\boldsymbol{e}_t=\boldsymbol{l}(\boldsymbol{\xi}(\boldsymbol{\theta}^*))-\boldsymbol{l}(\boldsymbol{\xi}(\hat{\boldsymbol{\theta}}^{-}_t))\\
    &\text{Prediction error: }\tilde{\boldsymbol{\theta}}^-_t = \boldsymbol{\theta}^*-\hat{\boldsymbol{\theta}}^-_t.
    \end{aligned}
\end{equation}
To perform the convergence analysis, a candidate Lyapunov function is employed:
\begin{equation} \label{eq:Lyapunov}
V_t=\tilde{\boldsymbol{\theta}}^\prime_t\boldsymbol{P}^{-1}_{t}\tilde{\boldsymbol{\theta}}_t.
    \end{equation}
The goal here is to determine conditions for which the candidate Lyapunov function $\{V_t\}_{t=1,2,\hdots}$ is a decreasing sequence, i.e. $V_{t+1}-V_{t}\le0, \forall t$. 
% One standard way to approximate the measurement error defined in \ref{eq:observationError} is the first-order approximation via the Taylor expansion, i.e. $\boldsymbol{e}_t\approx\boldsymbol{L}_t\tilde{\boldsymbol{\theta}}^-_t
% $. However, this approximation is available only if $\hat{\boldsymbol{\theta}}^-_{t}$ belongs to a neighborhood of $\boldsymbol{\theta}_{t}$.To obtain an exact equality rather than approximation, 
For rigorous analysis of the candidate Lyapunov function, as proposed in \cite{Boutayeb1997}, unknown diagonal matrices $\mathbfcal{F}_t\in\mathbb{R}^{r\times r}$ and $\mathbfcal{G}_t\in\mathbb{R}^{p\times p}$ are introduced to model the measurement and prediction error defined in (\ref{eq:observationError}):
\begin{equation} \label{eq:errorapproximate}
\mathbfcal{F}_t\boldsymbol{e}_t=\boldsymbol{L}_t\tilde{\boldsymbol{\theta}}^-_t, \quad \tilde{\boldsymbol{\theta}}^-_t=\mathbfcal{G}_t\tilde{\boldsymbol{\theta}}_{t-1}.
\end{equation}
% \begin{equation} \label{eq:errorapproximate}
% \begin{aligned}
% &\mathbfcal{F}_t\boldsymbol{e}_t=\boldsymbol{L}_t\tilde{\boldsymbol{\theta}}^-_t\\
% &\tilde{\boldsymbol{\theta}}^-_t=\mathbfcal{G}_t\tilde{\boldsymbol{\theta}}_{t-1}.
% \end{aligned}
% \end{equation}
To ensure convergence of the proposed estimator, the following assumptions need to be made. 
 
\begin{assumption} \label{assumption:observable}
The derivative $\boldsymbol{L}_t=\frac{d \boldsymbol{l}(\boldsymbol{\xi}_t(\hat{\boldsymbol{\theta}}^{-}_t))}{d \hat{\boldsymbol{\theta}}^{-}_t}$ is of full rank for every $\hat{\boldsymbol{\theta}}^{-}_t$.
\end{assumption}

\begin{remark}
 The discrete-time dynamical system (\ref{eq:dynTheta}) satisfies the observability rank condition, i.e., for every $\hat{\boldsymbol{\theta}}^{-}_t$, $\rank (\emph{col}\{\frac{d \boldsymbol{l}(\boldsymbol{\xi}_t(\hat{\boldsymbol{\theta}}^{-}_t))}{d \hat{\boldsymbol{\theta}}^{-}_t}, \frac{d \boldsymbol{l}(\boldsymbol{\xi}_t(\hat{\boldsymbol{\theta}}^{-}_t))}{d \hat{\boldsymbol{\theta}}^{-}_t}\boldsymbol{I}_p, \cdots, \frac{d \boldsymbol{l}(\boldsymbol{\xi}_t(\hat{\boldsymbol{\theta}}^{-}_t))}{d \hat{\boldsymbol{\theta}}^{-}_t}\boldsymbol{I}_p^{p-1} \}) = p$ \citep{Song1992}.
% \begin{equation}
% \rank (\emph{col}\{ \frac{d \boldsymbol{l}(\boldsymbol{\theta})}{d \boldsymbol{\theta}}, \frac{d \boldsymbol{l}(\boldsymbol{\theta})}{d \boldsymbol{\theta}}\boldsymbol{I}_p, \cdots, \frac{d \boldsymbol{l}(\boldsymbol{\theta})}{d \boldsymbol{\theta}}\boldsymbol{I}_p^{p-1} \}) = p.
% \end{equation} 
    % \begin{equation}
    %     \rank 
    %     \begin{bmatrix}
    %     \frac{\partial l(\boldsymbol{\theta})}{\partial \boldsymbol{\theta}} \\
    %     \frac{\partial l(\boldsymbol{\theta})}{\partial \boldsymbol{\theta}}\boldsymbol{I}_p\\
    %     \vdots\\
    %     \frac{\partial l(\boldsymbol{\theta})}{\partial \boldsymbol{\theta}}\boldsymbol{I}_p^{p-1}
    % \end{bmatrix}
    % =p.
    % \end{equation} 
That means if Assumption \ref{assumption:observable} is satisfied for every $\hat{\boldsymbol{\theta}}^{-}_t$, the system (\ref{eq:dynTheta}) is observable for every $\hat{\boldsymbol{\theta}}^{-}_t$. The observability condition assures that $\boldsymbol{P}_t$ is a bounded matrix from above and below \citep{Song1992,Boutayeb1999}. 
\end{remark}
As common in the EKF analysis, we adopt the following assumption:
\begin{assumption} \label{assumption:Lbounded}
    $\boldsymbol{L}_t$ is a uniformly bounded matrix.
\end{assumption}

We have the following lemma to show how the measurement covariance matrix $\boldsymbol{R}_t$ affects the convergence of the tunable parameter. The proof can be found in Appendix \ref{Appendix:ProofLemma2}.

\begin{lemma} \label{lemma:convergenceKalman}
    Let Assumptions \ref{assumption:observable} and \ref{assumption:Lbounded} hold. If the following inequalities are satisfied:
    \begin{equation} \label{eq:alphaR}
    (\mathbfcal{F}_t-\boldsymbol{I}_{s})^2\le\boldsymbol{R}_t(\boldsymbol{L}_t\boldsymbol{P}^{-}_{t}\boldsymbol{L}^\prime_t+\boldsymbol{R}_t)^{-1},
    \end{equation}
    \begin{equation} \label{eq:Q}
        \mathbfcal{G}_t^{\prime}\boldsymbol{P}_t^{-1}\mathbfcal{G}_t-\boldsymbol{P}_t^{-1}\leq 0,
    \end{equation}
   Then the proposed estimator (\ref{eq:main_EKF}), when used as an observer for the system (\ref{eq:dynTheta}), ensures local asymptotic convergence, i.e. $\lim_{t\rightarrow\infty}\tilde{\boldsymbol{\theta}}_t=\boldsymbol{0}$.

\end{lemma}
% Please refer to section \ref{Appendix:ProofLemma3} in the appendix for the proof of Lemma \ref{lemma:convergenceKalman}.

\begin{remark}
Lemma \ref{lemma:convergenceKalman} provides sufficient conditions for the convergence of $\hat{\boldsymbol{\theta}}_t$. As the diagonal matrices $\mathbfcal{F}_t$ and $\mathbfcal{G}_t$ are unknown, one can design the matrix $\boldsymbol{R}_t$ to satisfy inequalities (\ref{eq:alphaR}). For example, one can set the matrix $\boldsymbol{R}_t$ to be sufficiently large, i.e. much larger than $\boldsymbol{L}_t\boldsymbol{P}^{-}_{t}\boldsymbol{L}^\prime_t$, so that (\ref{eq:alphaR}) is satisfied, which means the parameter estimator can tolerate arbitrary large initial prediction error. It is worth to note that as long as (\ref{eq:alphaR}) and (\ref{eq:Q}) are satisfied, $\hat{\boldsymbol{\theta}}_t$ converges to $\boldsymbol{\theta}^*$ and consequently $\mathbfcal{F}_t$ and $\mathbfcal{G}_t$ become identity matrix. In the case when there is no measurement noise, i.e.  $\boldsymbol{R}_t=\boldsymbol{0}_{s\times s}$, $\mathbfcal{F}_t$ and $\mathbfcal{G}_t$ can only be identity matrices to satisfy the inequalities (\ref{eq:alphaR}) and (\ref{eq:Q}), indicating the convergence of $\hat{\boldsymbol{\theta}}_t$ to $\boldsymbol{\theta}^*$.
% According to \ref{eq:errorapproximate}, $\mathbfcal{F}_t=\boldsymbol{I}_s$ indicates that $\hat{\boldsymbol{\theta}}_t$ converges to $\boldsymbol{\theta}_t$. 
    % In the case when there is no noise, i.e. $\boldsymbol{Q}_t=\boldsymbol{0}_{p\times p}$ and $\boldsymbol{R}_t=\boldsymbol{0}_{s\times s}$, the candidate Lyapunov function \ref{eq:Lyapunov} constantly equals to zero. As the matrix $\boldsymbol{P}^{-1}_{t}$ is not a zero matrix, it indicates that $\tilde{\boldsymbol{\theta}}_t=\boldsymbol{0}$, which means $\hat{\boldsymbol{\theta}}_t$ converges to $\boldsymbol{\theta}_t$.
\end{remark}

\begin{remark} \label{remark:initialGuess}
    Equation (\ref{eq:alphaR}) and (\ref{eq:Q}) indicate one of the limitations of the estimator, which is the selection of initial guess. If the initial guess of $\boldsymbol{\theta}^*$ results in $\mathbfcal{F}_0$ and $\mathbfcal{G}_0$ that do not satisfy (\ref{eq:alphaR}) and (\ref{eq:Q}), the value of the Lyapunov function (\ref{eq:Lyapunov}) becomes larger, which leads to even larger $\mathbfcal{F}_t$ and $\mathbfcal{G}_t$, causing the estimation to diverge.
\end{remark}

% \begin{remark}
%     As the initial $\alpha_t$ is unknown, one intuitive idea is to set $R_t$ to be sufficiently large so that \ref{eq:alphaR} is satisfied, with a reasonable assumption that $\alpha_t$ is bounded. The sufficient condition \ref{eq:alphaR} also provides qualitative analysis when there exists measurement noise.
% \end{remark}

We have the following main theorem shows how the measurement covariance matrix $\boldsymbol{R}_t$ affects the convergence of cumulative loss $L(\boldsymbol{\xi}(\hat{\boldsymbol{\theta}}))$ by utilizing the inequalities introduced in Lemma \ref{lemma:convergenceKalman}. The proof can be found in Appendix \ref{Appendix:ProofTheorem}. 

\begin{theorem} \label{theorem}
    Let Assumptions \ref{assumption:observable} and \ref{assumption:Lbounded} hold. If the inequalities in Lemma~\ref{lemma:convergenceKalman} are met,
    % \begin{equation}
    %     (\mathbfcal{F}_t-\boldsymbol{I}_{s})^2\le\boldsymbol{R}_t(\boldsymbol{L}_t\boldsymbol{P}^{-}_{t}\boldsymbol{L}^\prime_t+\boldsymbol{R}_t)^{-1}, \quad
    %     \mathbfcal{G}_t^{\prime}(\boldsymbol{P}_t+\boldsymbol{Q}_t)^{-1}\mathbfcal{G}_t-\boldsymbol{P}_t^{-1}\leq 0,
    % \end{equation}
    then estimating $\boldsymbol{\theta}^*$ with the proposed estimator (\ref{eq:main_EKF}) employing the gradient generator in (\ref{eq:backwardPass})-(\ref{eq:recursiveGG}) ensures local asymptotic convergence of the cumulative loss $L$ in (\ref{eq:cumuLoss}) to 0 , i.e. $\lim_{t\rightarrow\infty}L(\boldsymbol{\xi}(\hat{\boldsymbol{\theta}}))=0$.
\end{theorem}
% Please refer to section \ref{Appendix:ProofTheorem} in the appendix for the proof of Theorem \ref{theorem}.

% \zihao{remove this part} To evaluate the effectiveness of online algorithms, the regret analysis is often adopted \citep{li2021onlineoptimization,zinkevich2003online,hazan2007logarithmic}. In our case, the regret at time $T$ is defined as:
% \begin{equation} \label{eq:regret}
%      \text{Regret}_T=\textstyle\sum_{t=1}^T(\|\boldsymbol{l}(\hat{\boldsymbol{\theta}}_t)\|^2-\|\boldsymbol{l}(\boldsymbol{\theta}^*_t)\|^2),
% \end{equation}
% which represents the accumulative performance discrepancy between $\hat{\boldsymbol{\theta}}_t$ and $\boldsymbol{\theta}^*_t$. Assume $\|\boldsymbol{l}(\boldsymbol{\theta})\|^2$ is convex with respect to $\boldsymbol{\theta}$. As Lemma~ \ref{lemma:convergenceKalman} states, the sufficient conditions (\ref{eq:alphaR}) and (\ref{eq:Q}) are satisfied for all $t$, then $\hat{\boldsymbol{\theta}}_t$ will converge to $\boldsymbol{\theta}^*_t$ asymptotically. Furthermore, according to (\ref{eq:errorapproximate}), $\boldsymbol{e}_t$ will converge to zero, which indicates that $\|\boldsymbol{l}(\hat{\boldsymbol{\theta}}_t)\|^2$ will converge to $\|\boldsymbol{l}(\boldsymbol{\theta}^*_t)\|^2$ asymptotically. Therefore, we can say that the regret in (\ref{eq:regret}) will converge to a bounded value as $t$ goes to infinity. Please refer to \ref{Appendix:regretAnalysis} in Appendix for details. 

\section{Applications to Different Online Learning Modes and Experiments} \label{sec:applications}
This section demonstrates the capability of the proposed OCIL framework with its three modes by three applications, Online Imitation Learning, Online System Identification, and Learning Policy on-the-fly.
This section includes a performance comparison with some state-of-the-art frameworks for three environments that are summarized in Table \ref{table:simEnv}.
Let $\boldsymbol{O}_t^* =\boldsymbol{h}(\boldsymbol{\xi}_t(\boldsymbol{\theta}^*))$ denotes the measurement without noise. The measurement noise is subject to a multivariate Gaussian distribution $\mathcal{N}(\boldsymbol{O}_t^*, \sigma^2\boldsymbol{I}_r)$.

To highlight the flexibility of OCIL, each experiment includes two phases: 1) online phase, where OCIL keeps learning the unknown parameter while new data comes in before the final time $T$; 2) offline phase, where OCIL keeps learning the parameter given the learned parameter at time $T$ and the entire trajectory obtained from time $t=0$ to time $T$.
For each environment and task, a terminal time $T \in \mathbb{Z}$ is defined to represent a desired time duration where the system shall finish the task.

To unify the data visualization of both online and offline phases, the horizontal axis represents the number of data points, where a vertical red line corresponds to the final time $T$, i.e. the end of the online phase.
The number of data points reflects the number of iterations multiplied by the total number of time steps for each iteration. 
The solid blue curves indicate the online portion of OCIL, whereas the dashed blue curves indicate the offline portion.
For every environment and every method, 5 trials are performed given random initial conditions due to the high computational cost for other methods.
The computational performance and analysis for OCIL are shown in Section \ref{Appendix:computationPerformacnce} of the Appendix.

% \vspace{-1em}
\begin{table}[ht]
\caption{Experiment Environments}
\begin{center}
	\begin{tabular}{c| c| c}
		\hline
		\textbf{Systems}  & Dynamics parameter $\boldsymbol{\theta}_{dyn}$   & Objective parameter $\boldsymbol{\theta}_{obj}$ \\
		\hline
		 Cartpole &    cart mass, pole mass and length& \multirowcell{3}{$c(\boldsymbol{x},\boldsymbol{u})=\boldsymbol{\theta}_{obj}\|\boldsymbol{x}-\boldsymbol{x}_g\|^2+\|\boldsymbol{u}\|^2$\\$h(\boldsymbol{x})=\boldsymbol{\theta}_{obj}\|\boldsymbol{x}-\boldsymbol{x}_g\|^2$}  \\
		% \hline
		6-DoF Quadrotor & mass, wing length, inertia matrix\\
		% \hline
		 6-DoF Rocket & mass, rocket length, inertia matrix \\
		\hline
	\end{tabular}
	\label{table:simEnv}
	\end{center}
\end{table}
% \vspace{-1em}

\textbf{Online Imitation Learning.}
The control objective is parameterized as a weighted distance to the goal. Set the signed residual function of imitation learning $l(\boldsymbol{\xi}_t(\hat{\boldsymbol{\theta}}),\boldsymbol{y}_t^*)=\boldsymbol{y}_t^*-\boldsymbol{g}(\boldsymbol{x}_t(\hat{\boldsymbol{\theta}}),\boldsymbol{u}_t(\hat{\boldsymbol{\theta}}))$. The optimal cumulative loss is zero, i.e.  $ L(\boldsymbol{\xi}(\boldsymbol{\theta}^*))=0$, with full knowledge of the parameter. Four existing methods are used for comparisons: (i) inverse KKT \citep{englert2017inverse} (ii) neural policy cloning \citep{bojarski2016end} and (iii) PDP \citep{jin2020pontryagin}.
These methods don't handle measurement noise well because of their limitations, so we performed the experiments without including measurement noise for these methods. For OCIL, $\sigma =0.1$ for all of the systems.

\begin{figure}[ht]
     \centering
     \begin{subfigure}[b]{0.3\textwidth}
         \centering
         \includegraphics[width=\textwidth]{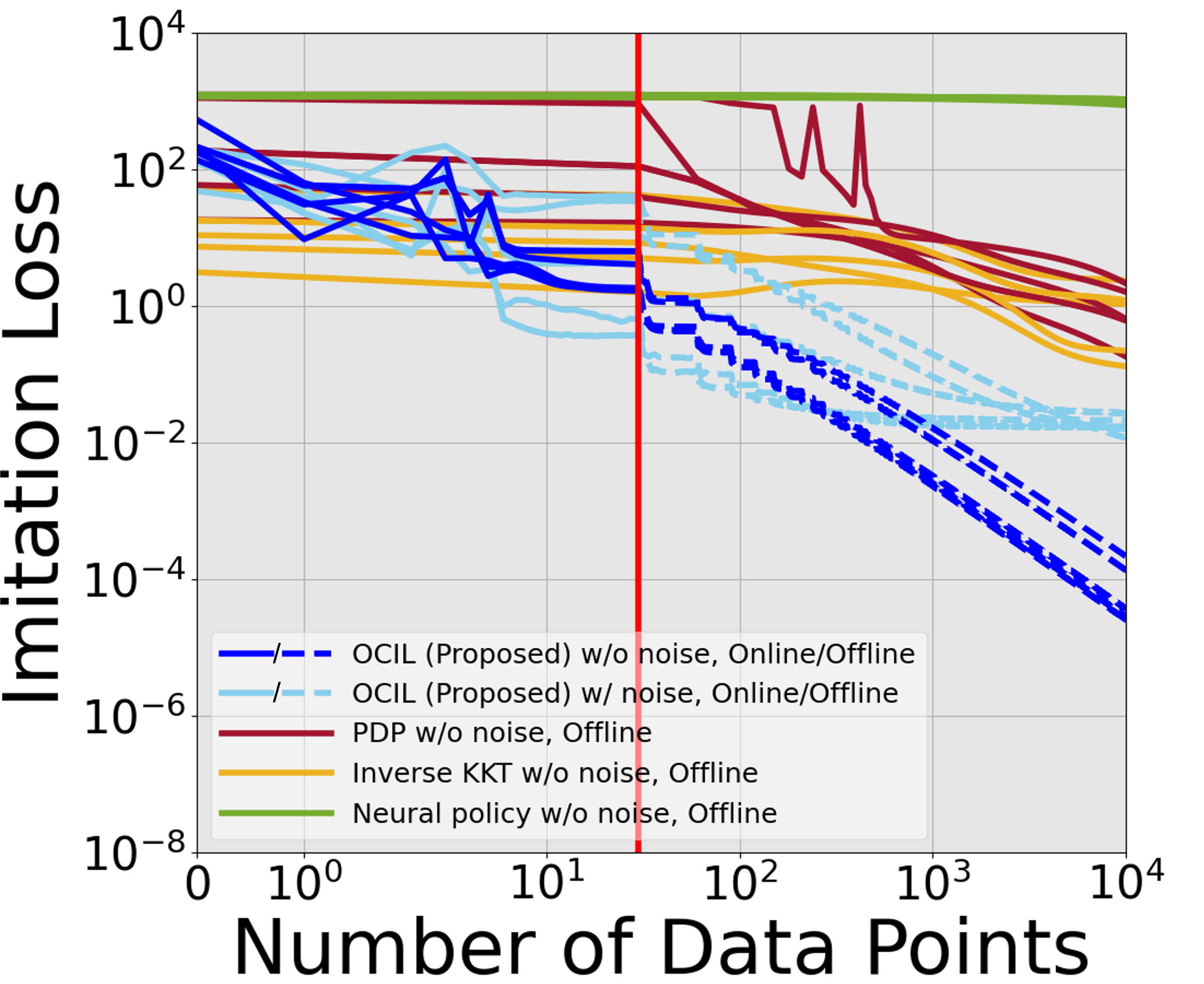}
         \caption{Cartpole}
         \label{fig:IL_Cartpole}
     \end{subfigure}
     \begin{subfigure}[b]{0.3\textwidth}
         \centering
         \includegraphics[width=\textwidth]{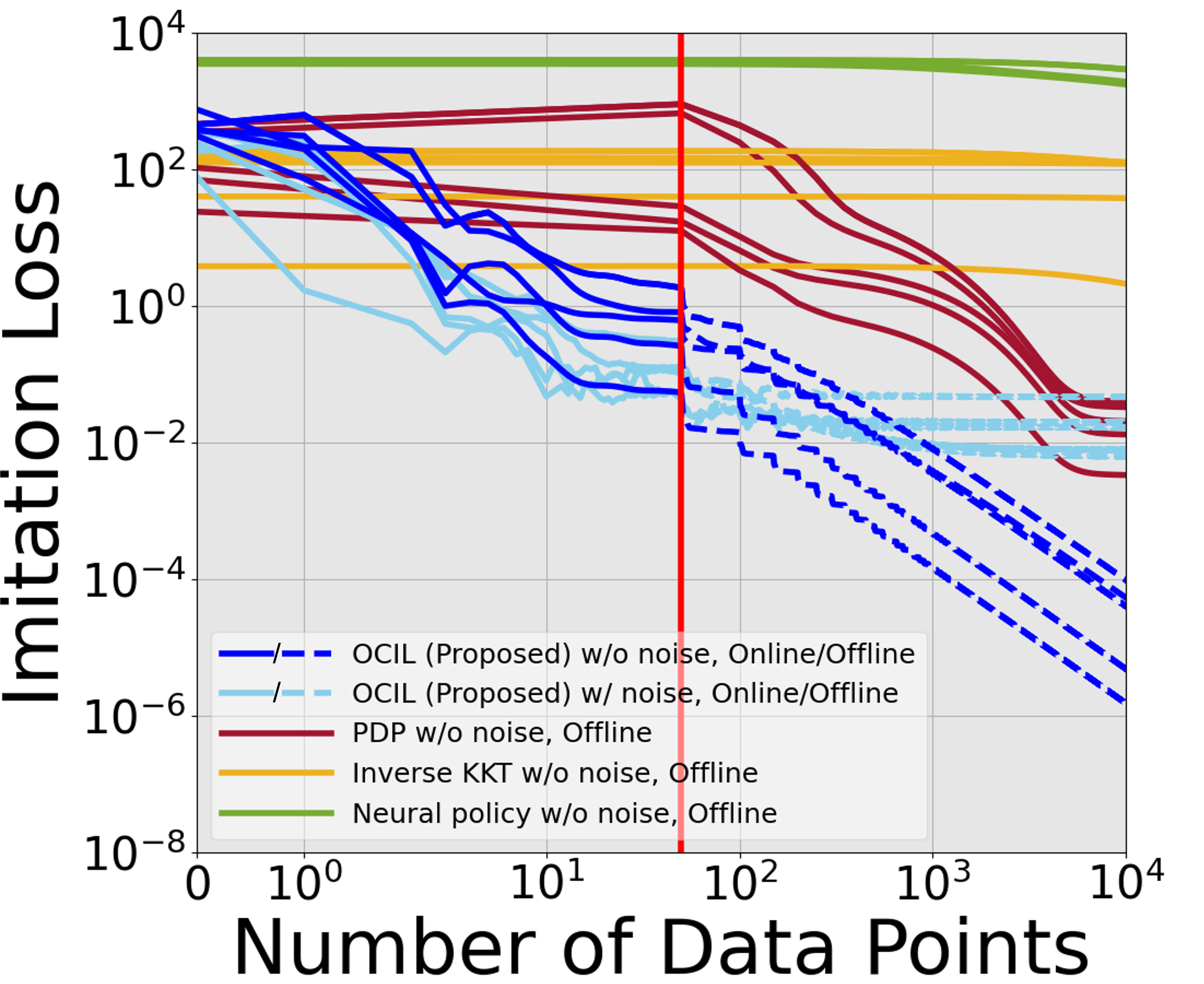}
         \caption{Quadrotor}
         \label{fig:IL_Quadrotor}
     \end{subfigure}
     \begin{subfigure}[b]{0.3\textwidth}
         \centering
         \includegraphics[width=\textwidth]{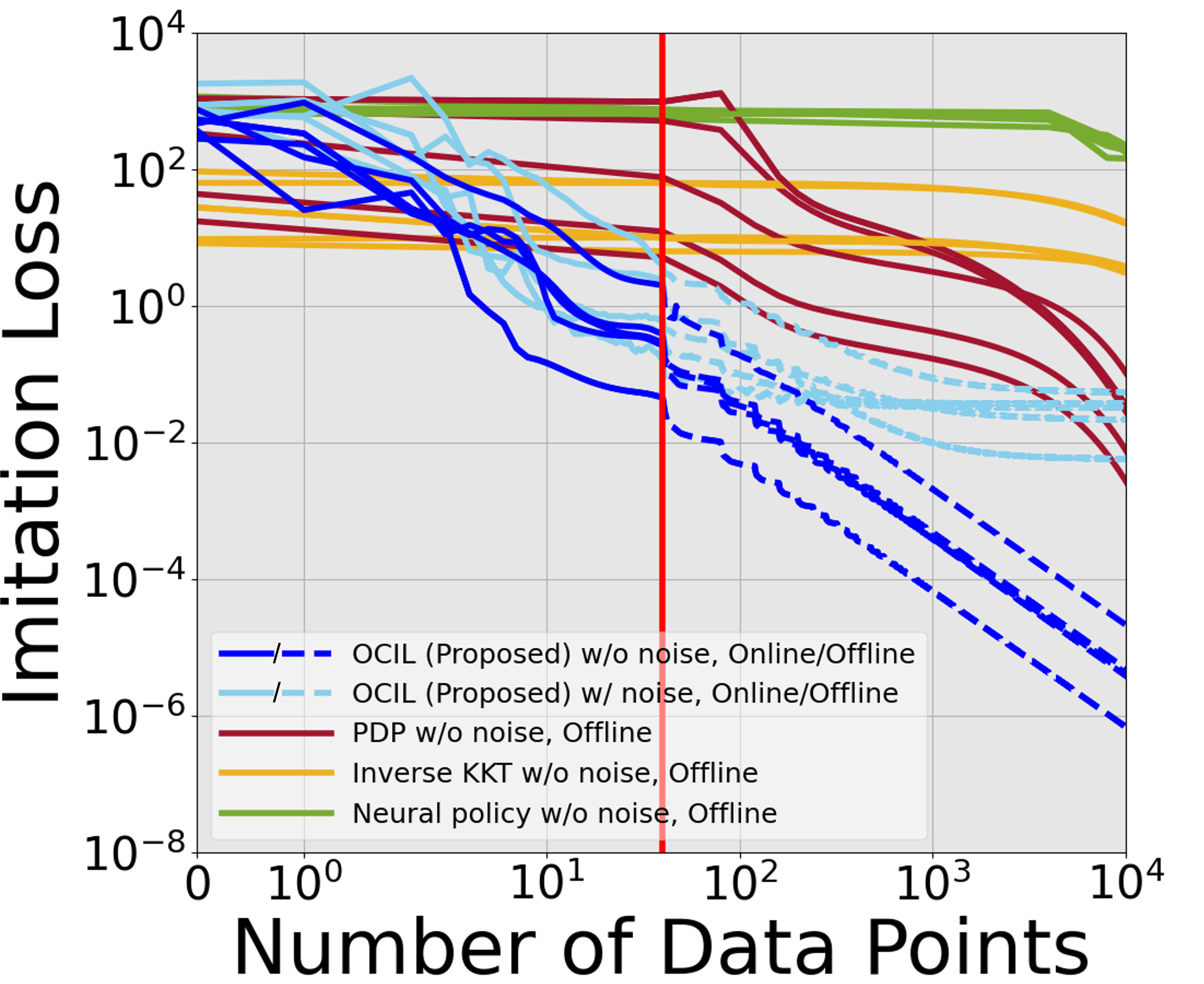}
         \caption{Rocket}
         \label{fig:IL_Rocket}
     \end{subfigure}
     %    \begin{subfigure}[b]{0.24\textwidth}
     %     \centering
     %     \includegraphics[width=\textwidth]{IR_cartpole.png}
     %     \caption{Neural}
     %     \label{fig:IL_Neural}
     % \end{subfigure}
     \caption{Imitation loss v.s. number of data points}
        \label{fig:IL}
\end{figure}

Fig.~\ref{fig:IL_Cartpole}-\ref{fig:IL_Rocket} summarize the comparison result, where OCIL converges faster and obtains lower loss than the other offline methods, in both online and offline phases.
The initial loss for each method is different because the learning representation (parameterization) is different. Thus, it is hard to guarantee that an initial neural network has the same loss as another initial parameter vector.
Nevertheless, the initial representation of each method is adjusted such that OCIL does not take advantage of good initialization.
Fig.~\ref{fig:IL_Cartpole}-\ref{fig:IL_Rocket} validate the effectiveness of OCIL's both online and offline performance, even with measurement noise.

% \vspace{-1em}
% \vspace{-1em}

\begin{figure}[ht]
\centering
\includegraphics[width=0.30\textwidth]{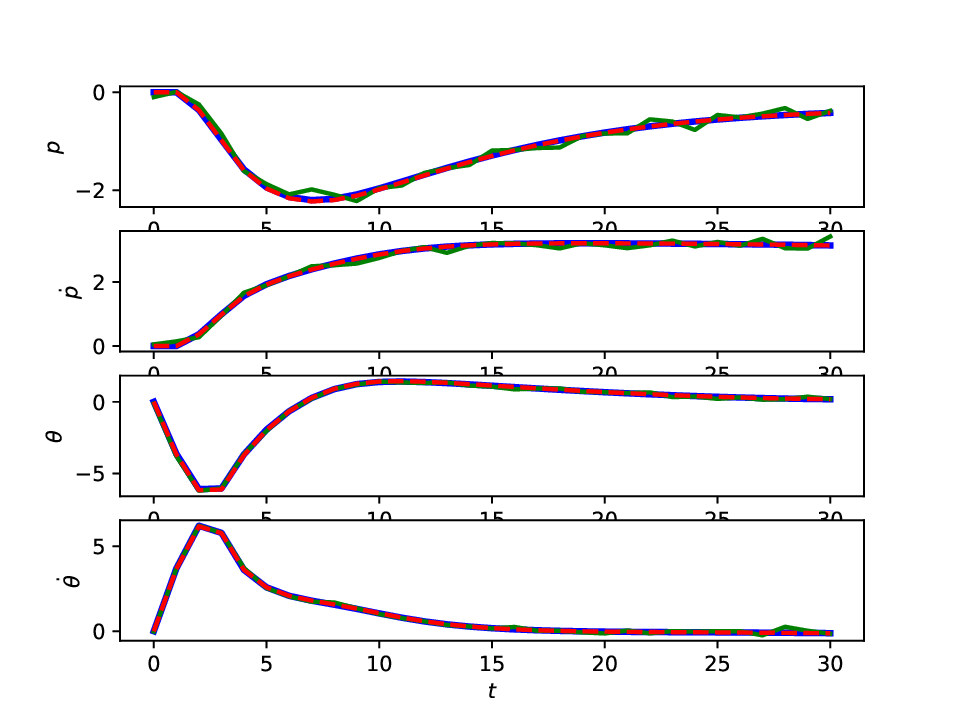}
\centering\includegraphics[width=0.30\linewidth]{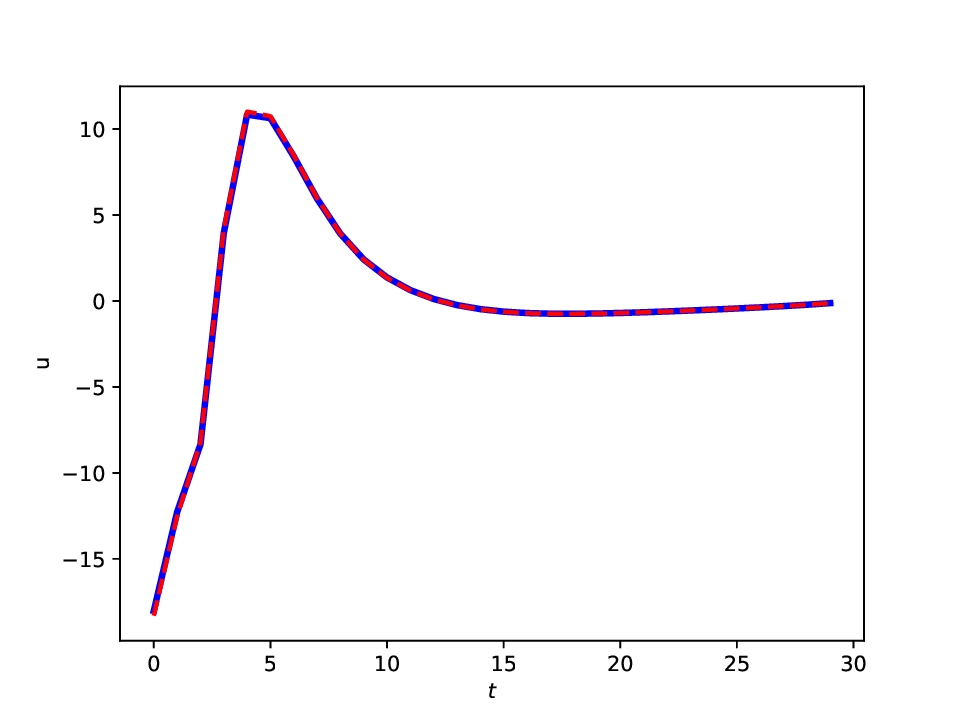}
\caption{Trajectories of the cartpole in online imitation learning. Blue solid lines: learned trajectory. Green solid lines: observed noisy trajectory. Red dashed lines: ground truth.}
\label{fig/IL_cartpole_state}
\end{figure}

\begin{figure*}[ht]
\centering\includegraphics[width=0.32\linewidth]{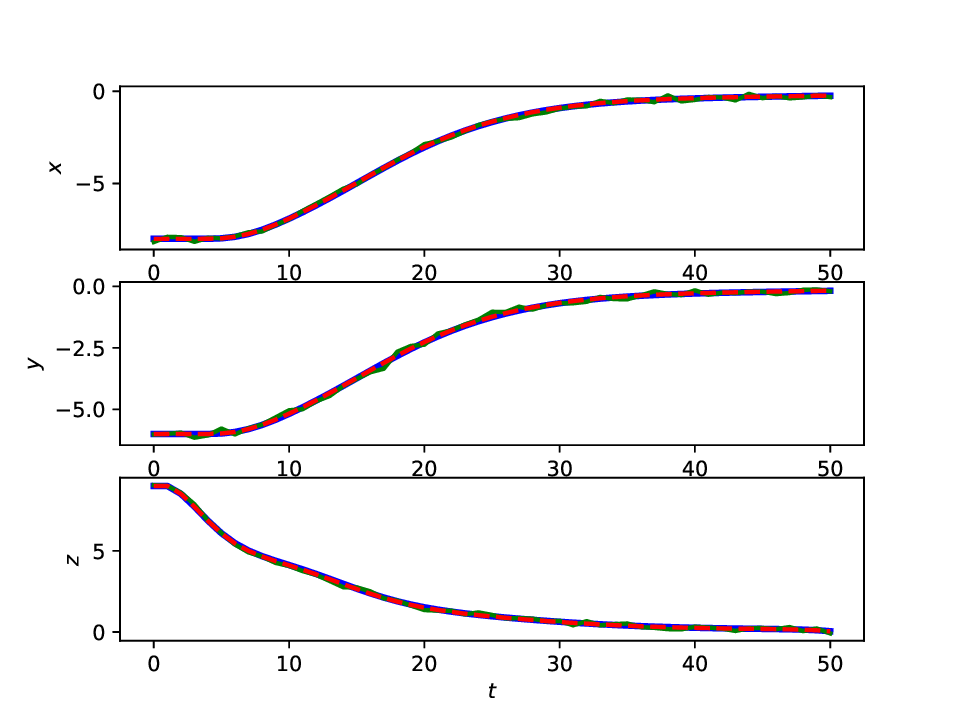}\label{fig:IL_uav_state_1}
\centering\includegraphics[width=0.32\linewidth]{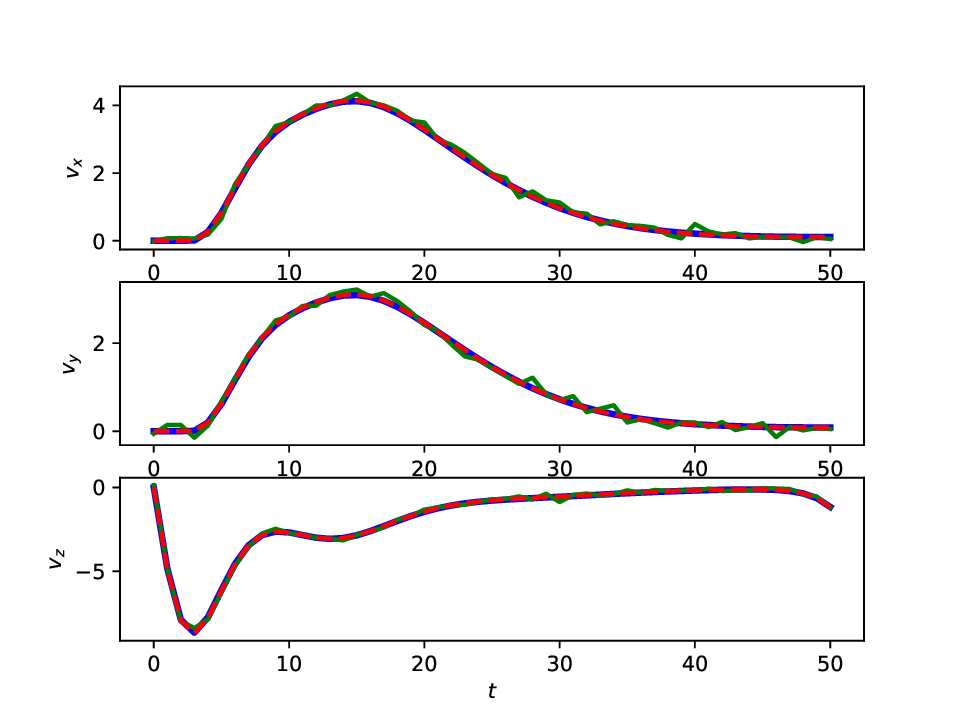}
% \hfill
\centering\includegraphics[width=0.32\linewidth]{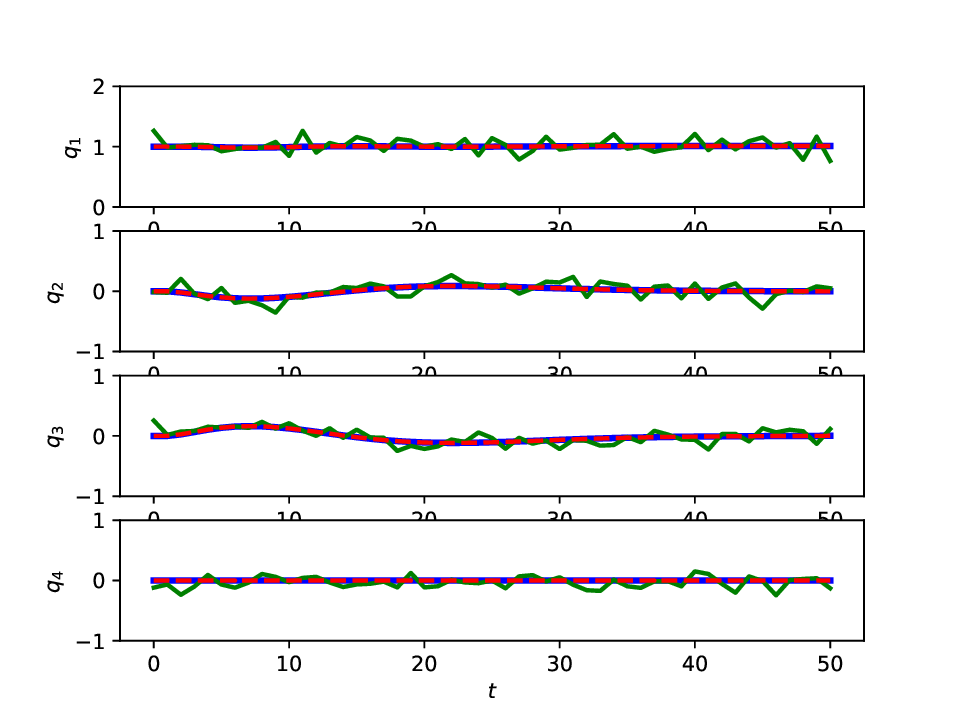}
\centering\includegraphics[width=0.32\linewidth]{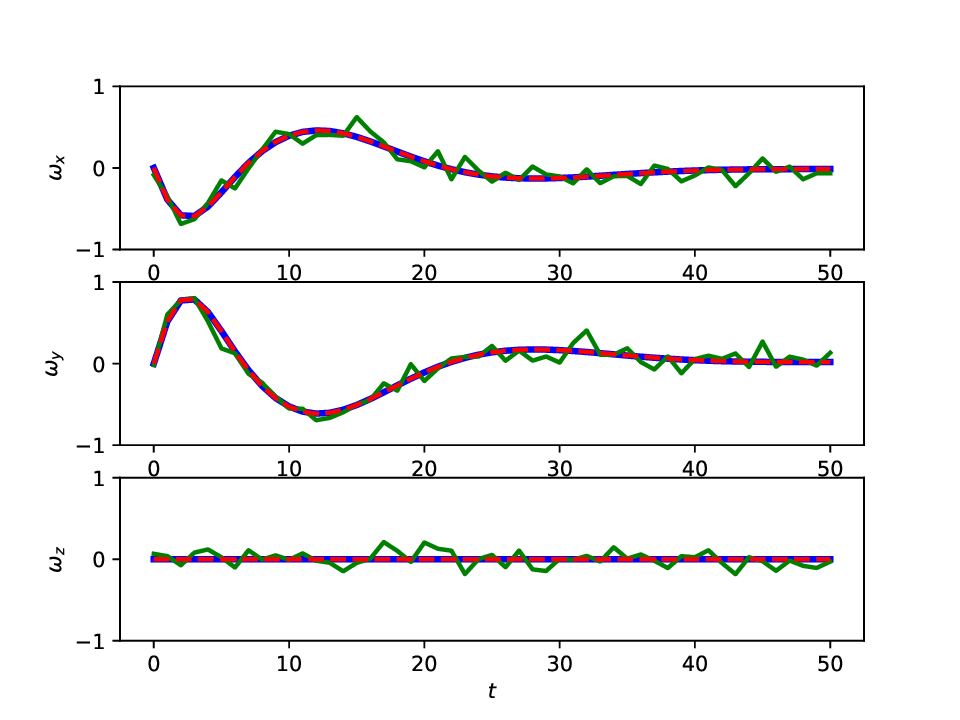}
\centering\includegraphics[width=0.32\linewidth]{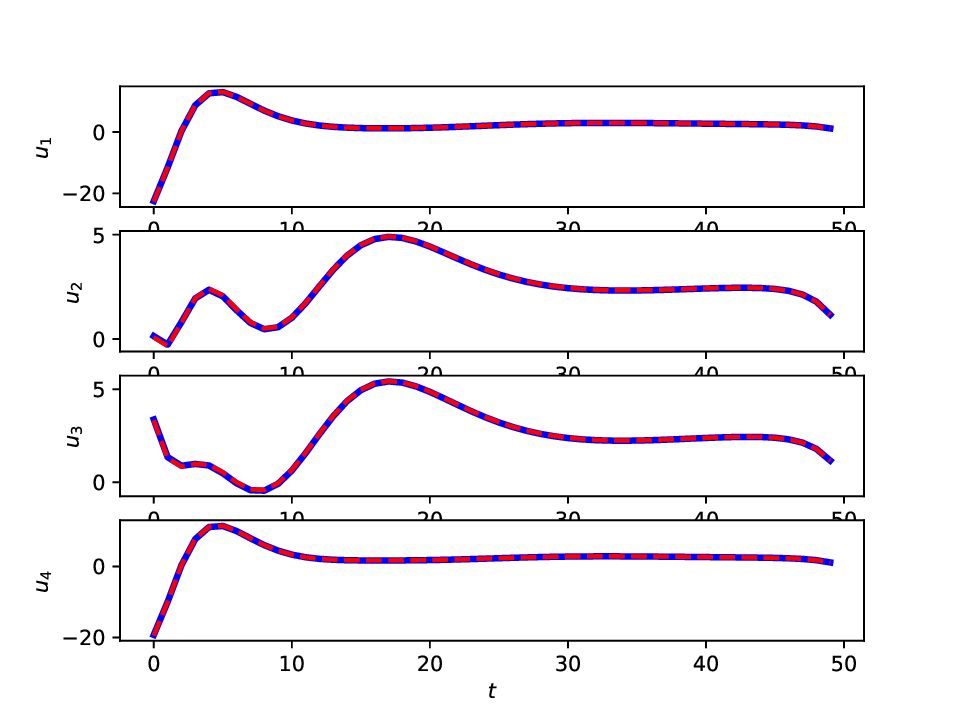}
\caption{Trajectories of the quadrotor in online imitation learning. Blue solid lines: learned trajectory. Green solid lines: observed noisy trajectory. Red dashed lines: ground truth.} \label{fig:IL_uav_state}
\end{figure*} 

\textbf{Online System Identification}
% OCIL is used to solve the online system identification problem with environments in Table \ref{table:simEnv}.
% In this mode the system in \ref{sys:sysID} is used.
The signed residual function is set to be $l(\boldsymbol{\xi}_t(\hat{\boldsymbol{\theta}}),\boldsymbol{\xi}_t^o)=\boldsymbol{\xi}_t^o-\boldsymbol{\xi}_t(\hat{\boldsymbol{\theta}})$. The optimal cumulative loss is zero, i.e.  $ L(\boldsymbol{\xi}(\boldsymbol{\theta}^*))=0$, with full knowledge of the parameter. Three other methods are used for comparison: (i) Pytorch Adam solver \citep{pillonetto2025deep}, (ii) DMDc \citep{proctor2016dynamic}, and (iii) PDP \citep{jin2020pontryagin}.
No measurement noise are injected into observed data for existing methods due to their inherent limitations. For OCIL, $\sigma =0.05$ for all of the systems.

\begin{figure}[ht]
     \centering
     \begin{subfigure}[b]{0.3\textwidth}
         \centering
         \includegraphics[width=\textwidth]{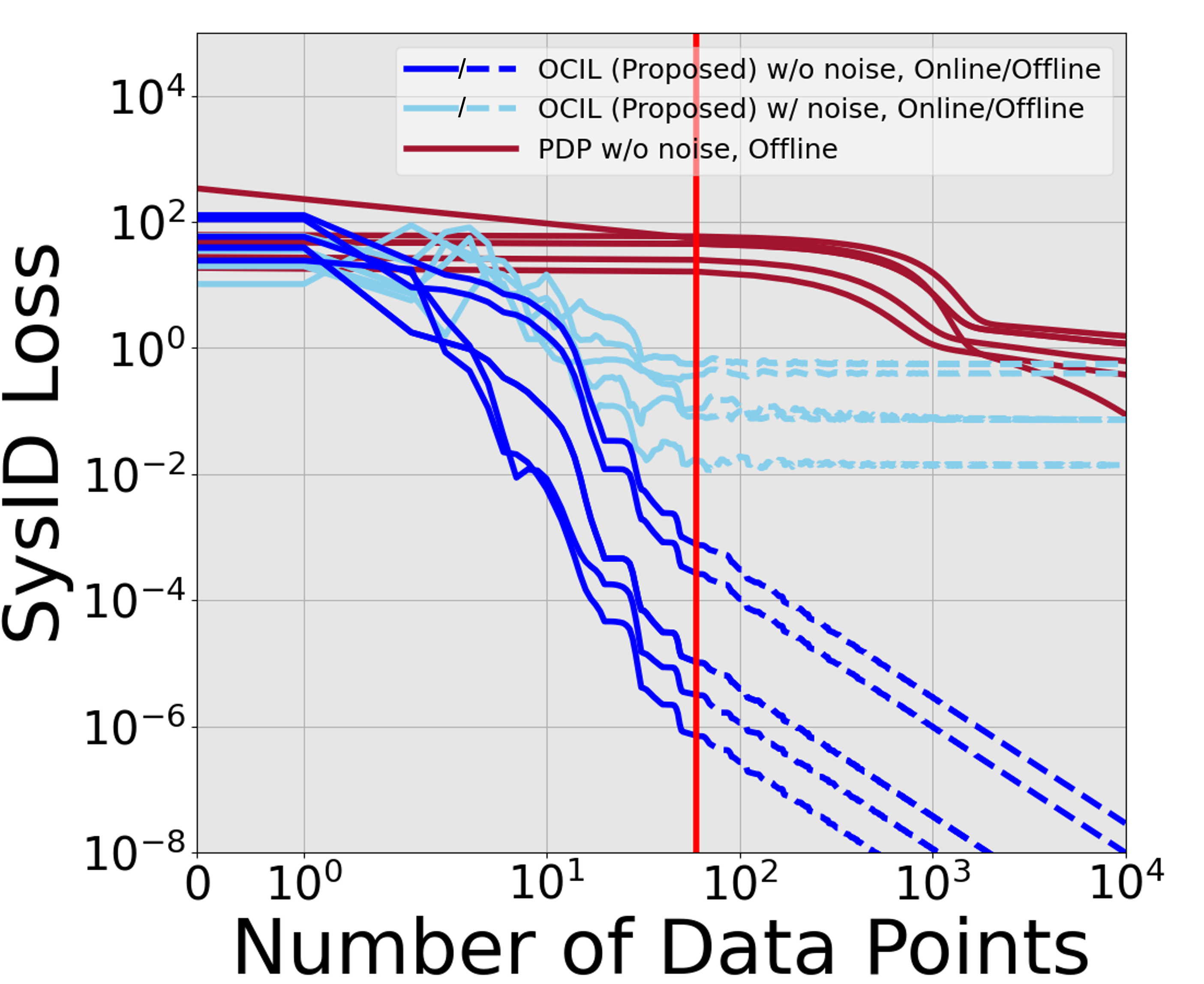}
         \caption{Cartpole}
         \label{fig:SysId_Cartpole}
     \end{subfigure}
     \begin{subfigure}[b]{0.3\textwidth}
         \centering
         \includegraphics[width=\textwidth]{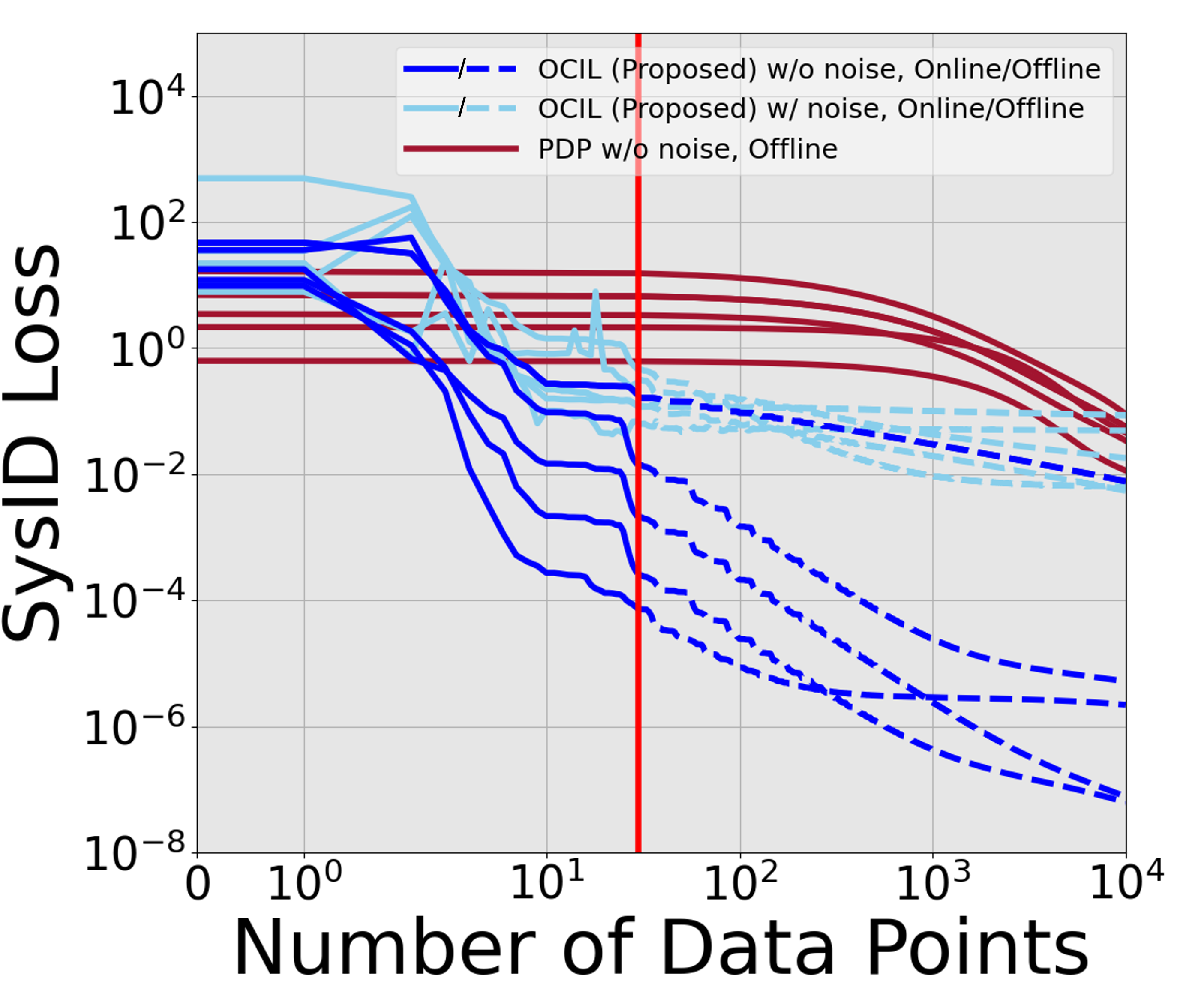}
         \caption{Quadrotor}
         \label{fig:SysId_Quadrotor}
     \end{subfigure}
     \begin{subfigure}[b]{0.3\textwidth}
         \centering
         \includegraphics[width=\textwidth]{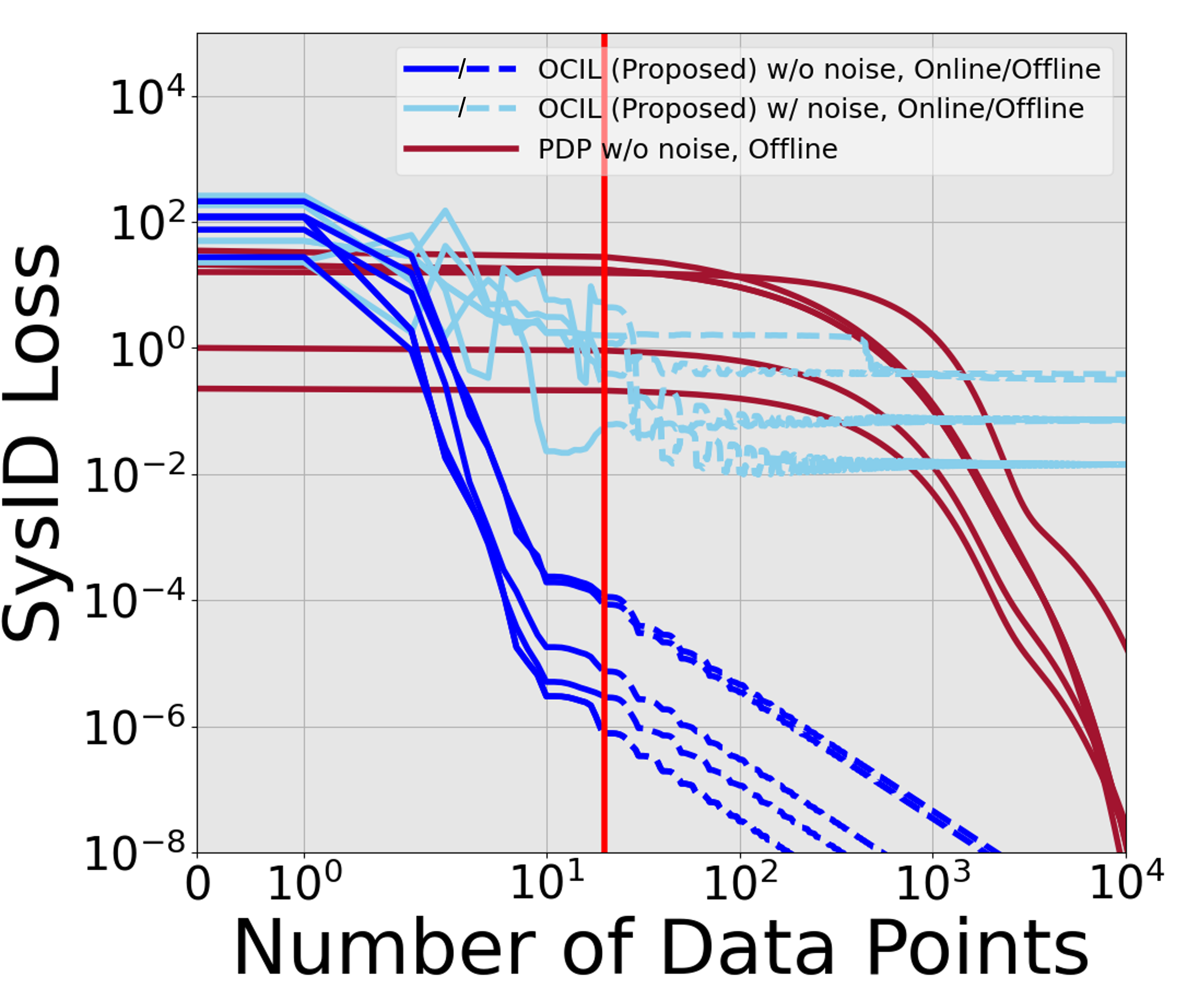}
         \caption{Rocket}
         \label{fig:SysId_Rocket}
     \end{subfigure}
     \vskip 1pt%\baselineskip
     \centering
     \begin{subfigure}[b]{0.3\textwidth}
         \centering
         \includegraphics[width=\textwidth]{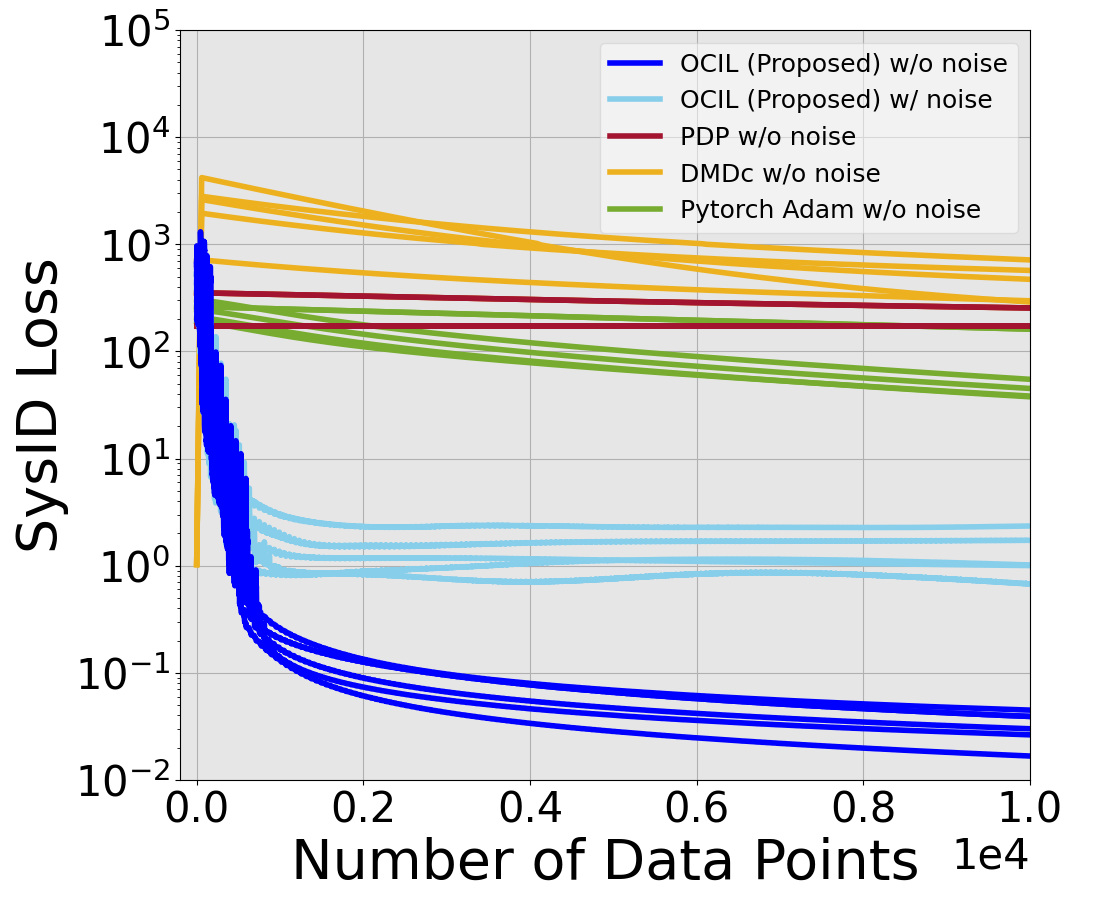}
         \caption{Cartpole, offline, NN dynamics}
         \label{fig:SysID_Cartpole_NN}
     \end{subfigure}
     \begin{subfigure}[b]{0.3\textwidth}
         \centering
         \includegraphics[width=\textwidth]{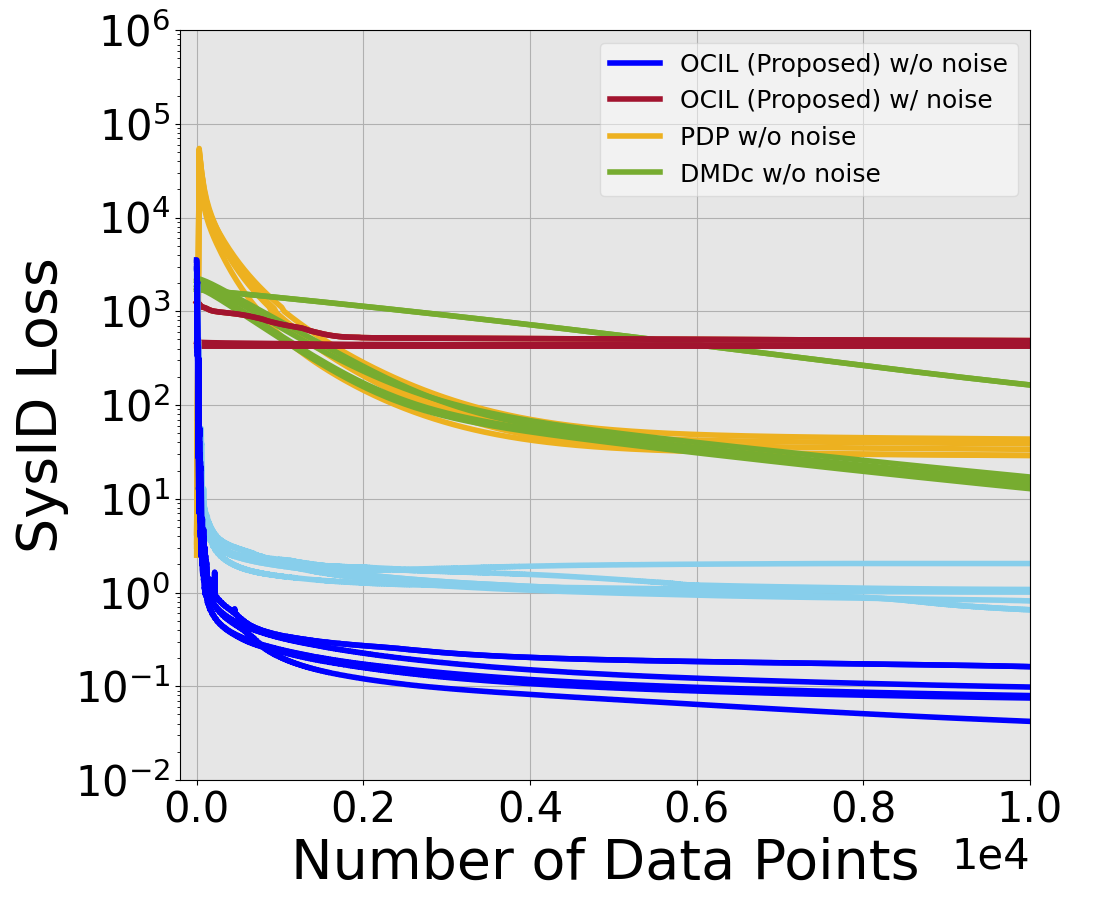}
         \caption{Quadrotor, offline, NN dynamics}
         \label{fig:SysID_Quadrotor_NN}
     \end{subfigure}
     \begin{subfigure}[b]{0.3\textwidth}
         \centering
         \includegraphics[width=\textwidth]{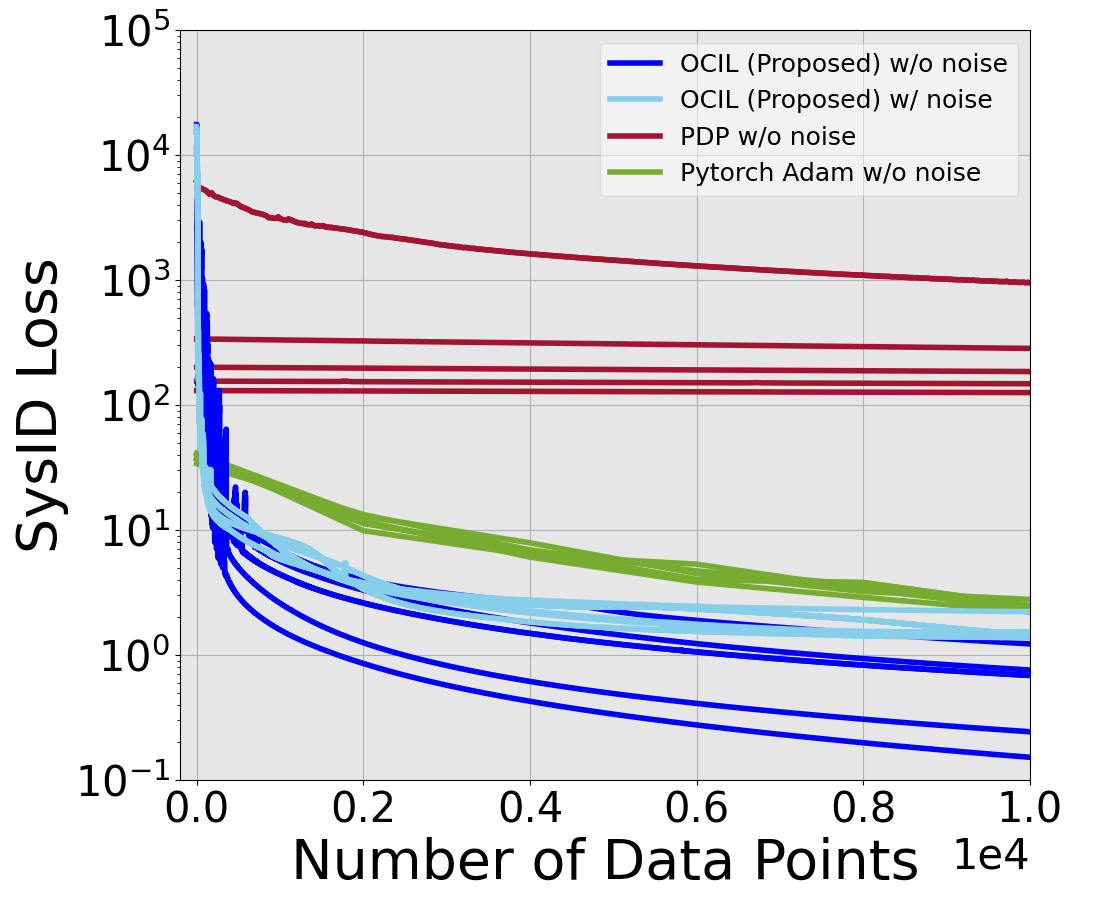}
         \caption{Rocket, offline, NN dynamics}
         \label{fig:SysID_Rocket_NN}
         \end{subfigure}
     \caption{SysID loss v.s. number of data points}
        \label{fig:SysId}
\end{figure}

\begin{figure}[ht] \centering
\begin{subfigure}[b]{0.3\textwidth}
     \centering
     \includegraphics[width=\textwidth]{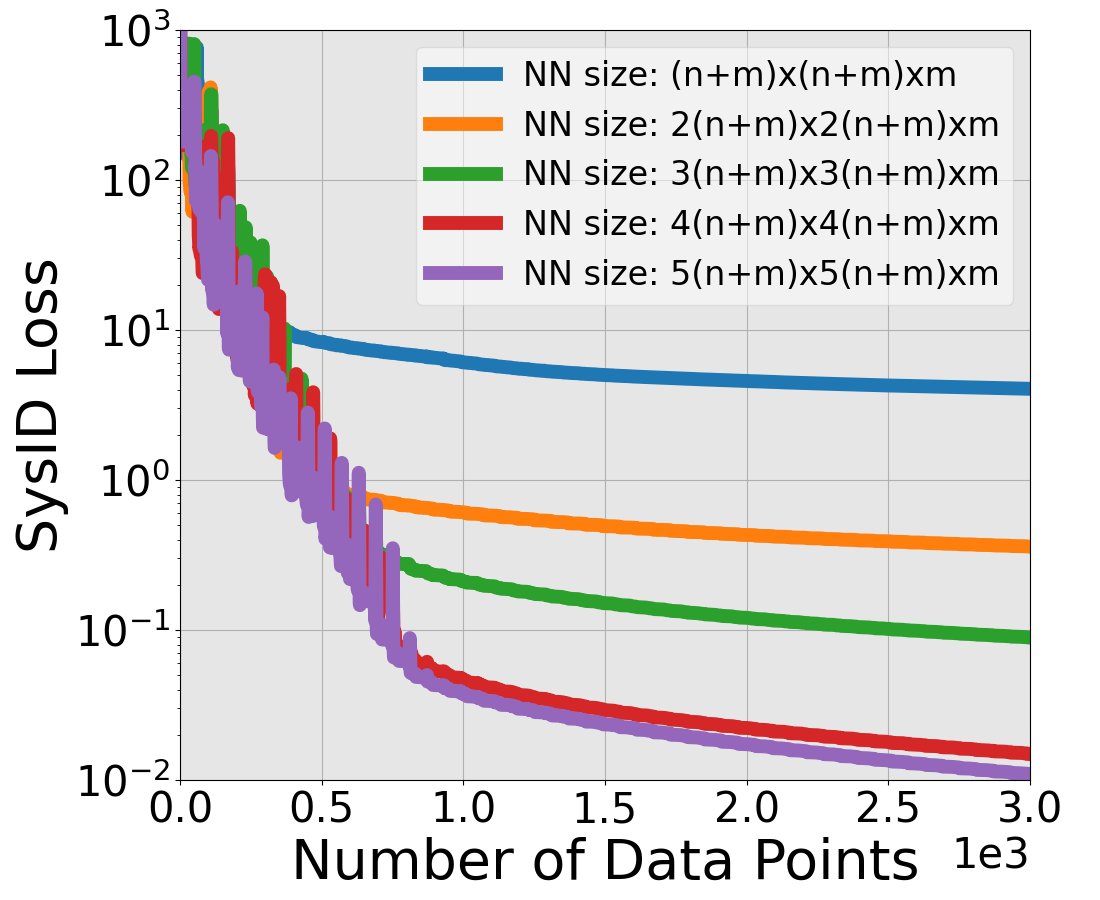}
     \label{fig:SysID_Cartpole_nn}
 \caption{Cartpole w/o noise}
\end{subfigure}
\begin{subfigure}[b]{0.3\textwidth}
     \centering
     \includegraphics[width=\textwidth]{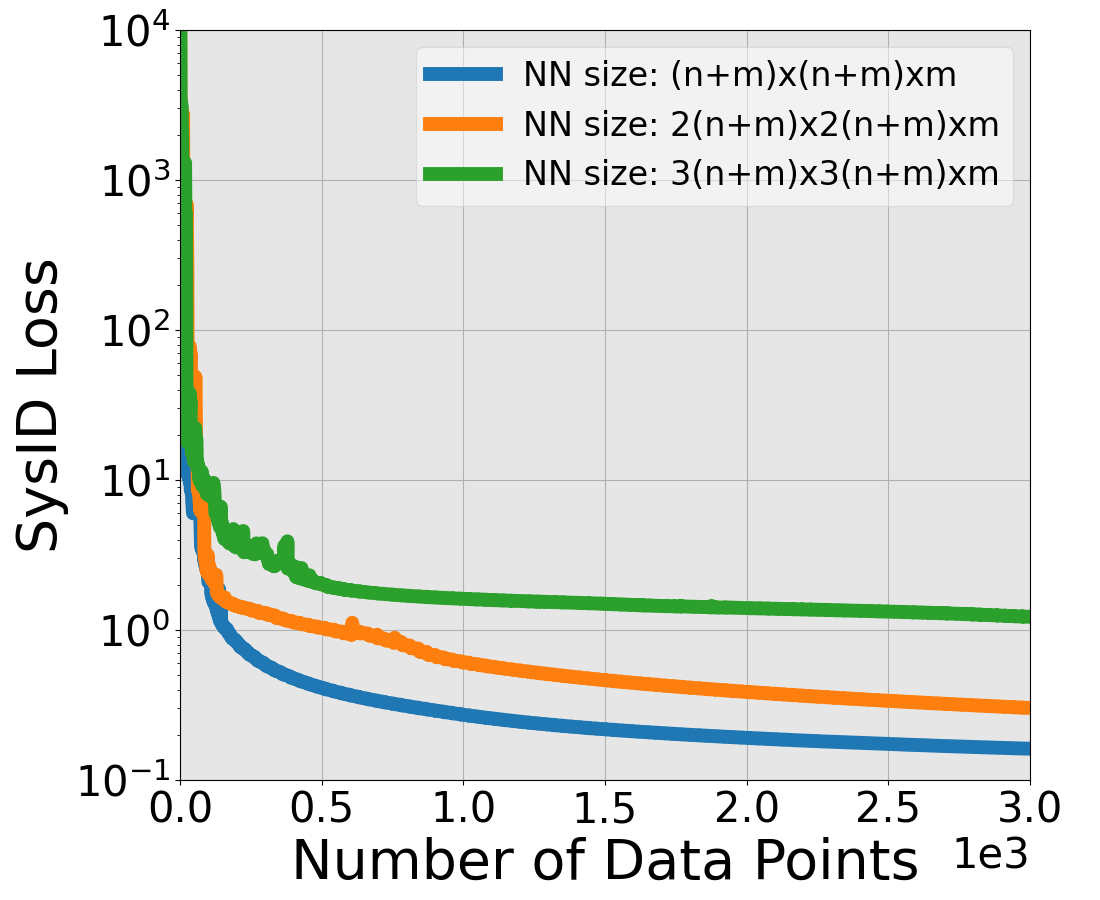}
     \label{fig:SysID_uav_nn}
 \caption{Quadrotor w/o noise}
\end{subfigure}
\begin{subfigure}[b]{0.3\textwidth}
     \centering
     \includegraphics[width=\textwidth]{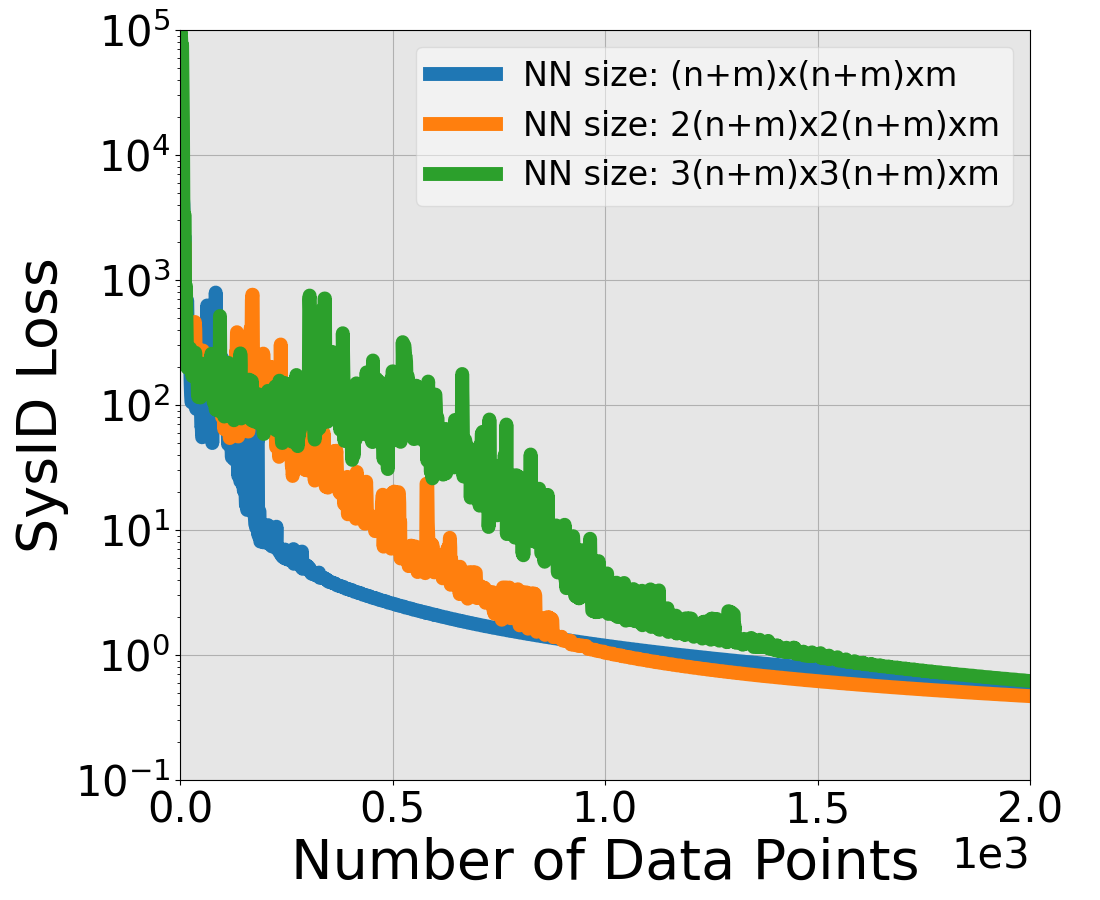}
     \label{fig:SysID_rocket_input}
 \caption{Rocket w/o noise}
\end{subfigure}
\caption{SysID Loss v.s. number of data points, given different sizes of neural dynamics}
\label{fig:SysIdwithDifferentNN}
\end{figure}

Fig.~\ref{fig:SysId_Cartpole}-\ref{fig:SysId_Rocket} summarize the result, where OCIL outperforms PDP for faster convergence and lower loss, in both online and offline phases.
Different than Online Imitation Learning, OCIL does not decrease its SysID loss significantly at first because the number of data points is not sufficient for online learning. Once the number of data points becomes sufficient, the SysID loss starts decreasing significantly.
This phenomenon can also be observed in the other methods, but their critical number of data points is significantly larger than OCIL's.
In Fig.~\ref{fig:SysID_Cartpole_NN}-\ref{fig:SysID_Rocket_NN}, OCIL and other methods are applied to learn the neural dynamics using the same observed trajectory. It can be seen that OCIL outperforms other methods for lower loss. 
Fig.~\ref{fig:SysIdwithDifferentNN} demonstrates the capability of OCIL to deal with neural dynamics that have different sizes of NN.

% \vspace{-1em}
% \vspace{-1em}

\textbf{Policy Tuning On-the-fly.} 
% OCIL is used to do Policy Tuning On-the-fly.
The parameterized OC system in \ref{sys:Control} is used here, where the policy is in a state-feedback form and parameterized by the tunable parameter $\hat{\boldsymbol{\theta}}$. The signed residual function is set to be $l(\boldsymbol{\xi}_t(\hat{\boldsymbol{\theta}}), \boldsymbol{\xi}_t^*)=\boldsymbol{\xi}_t^*-\boldsymbol{\xi}_t(\hat{\boldsymbol{\theta}})$, where $\boldsymbol{\xi}_t^*$ is the trajectory that needs to be tracked. The optimal cumulative loss is zero, i.e.  $ L(\boldsymbol{\xi}(\boldsymbol{\theta}^*))=0$, with full knowledge of the parameter. Other methods are used for comparison (i) iLQR \citep{li2004iterative} (ii) GPS \citep{levine2014learning}, and (iii) PDP \citep{jin2020pontryagin}.
No measurement noise is included for existing methods due to their limitations. For OCIL, $\sigma =0.1$ for cartpole and quadrotor; $\sigma =0.25$ for rocket.

\begin{figure}[ht]
     \centering
     \begin{subfigure}[b]{0.3\textwidth}
         \centering
          \includegraphics[width=\textwidth]{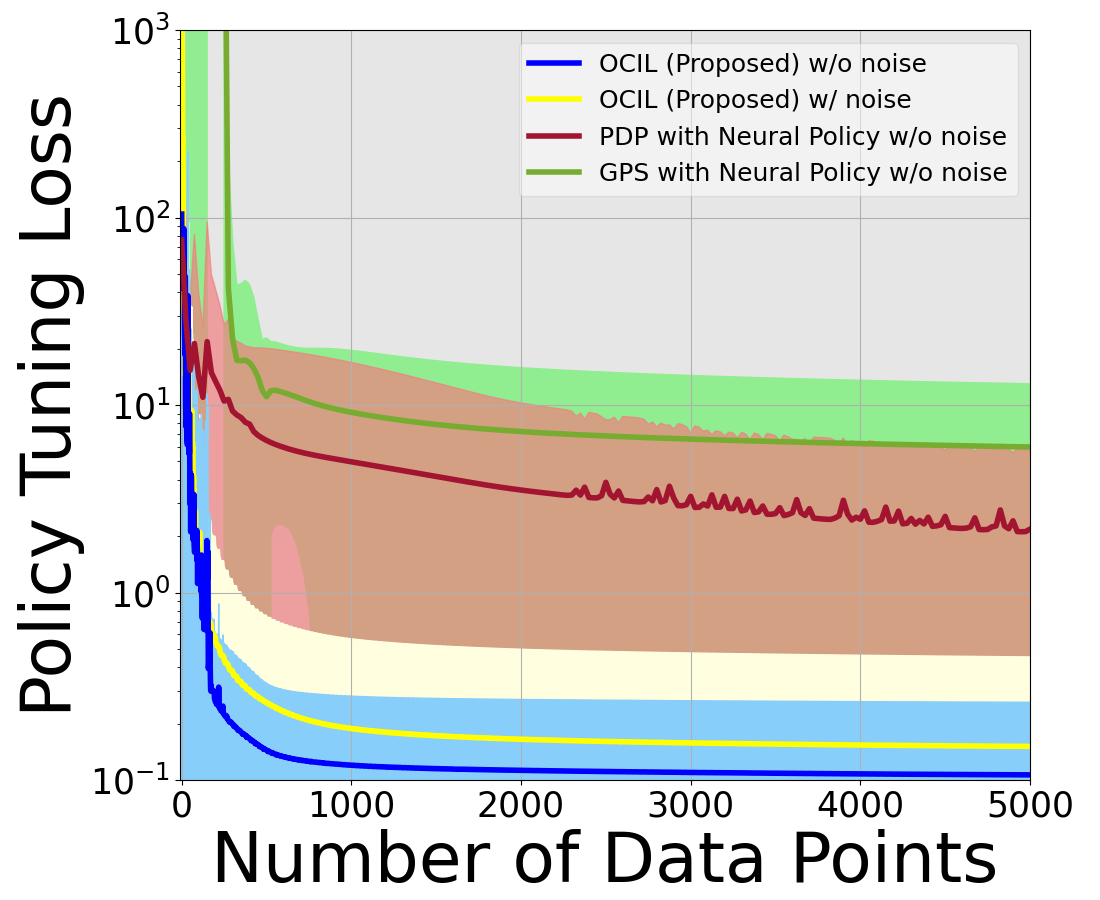}
         \caption{Cartpole}
         \label{fig:PT_Cartpole}
     \end{subfigure}
     \begin{subfigure}[b]{0.3\textwidth}
         \centering
         \includegraphics[width=\textwidth]{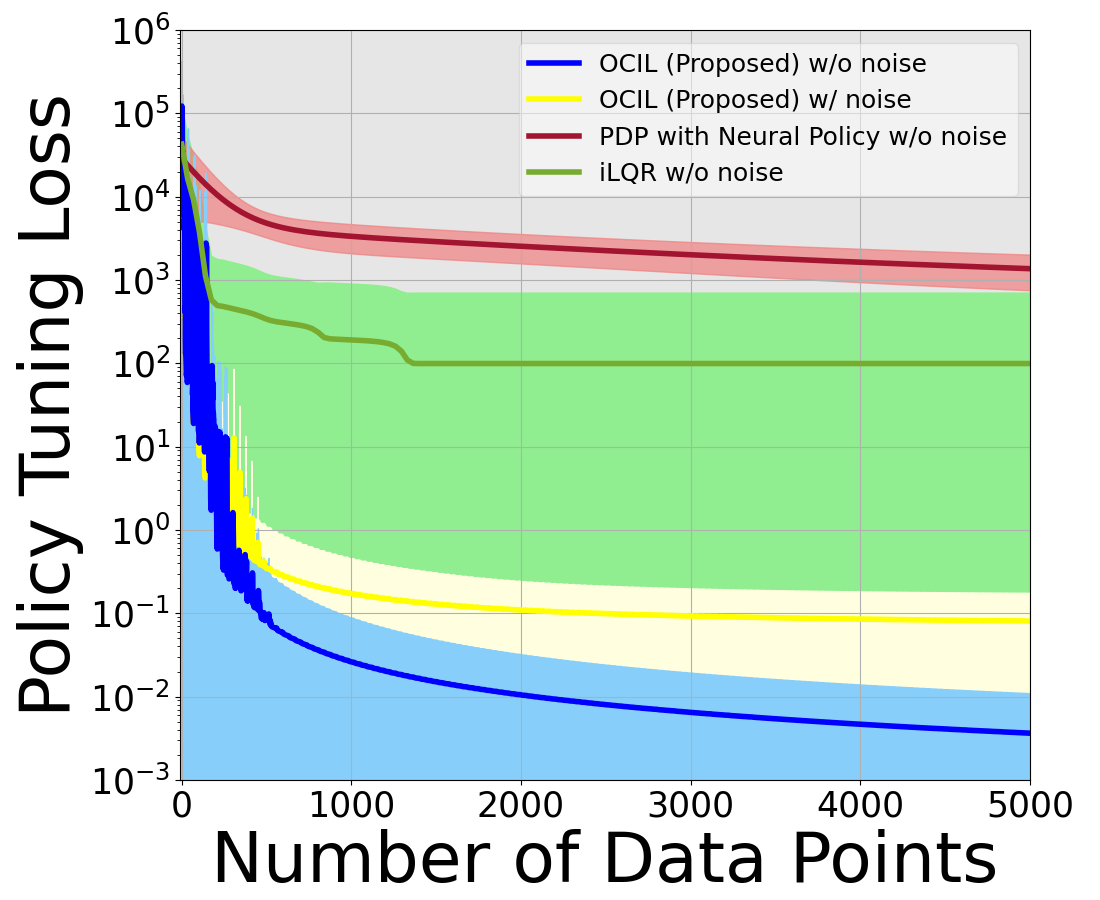}
         \caption{Quadrotor}
         \label{fig:PT_Quadrotor}
     \end{subfigure}
     \begin{subfigure}[b]{0.3\textwidth}
         \centering
         \includegraphics[width=\textwidth]{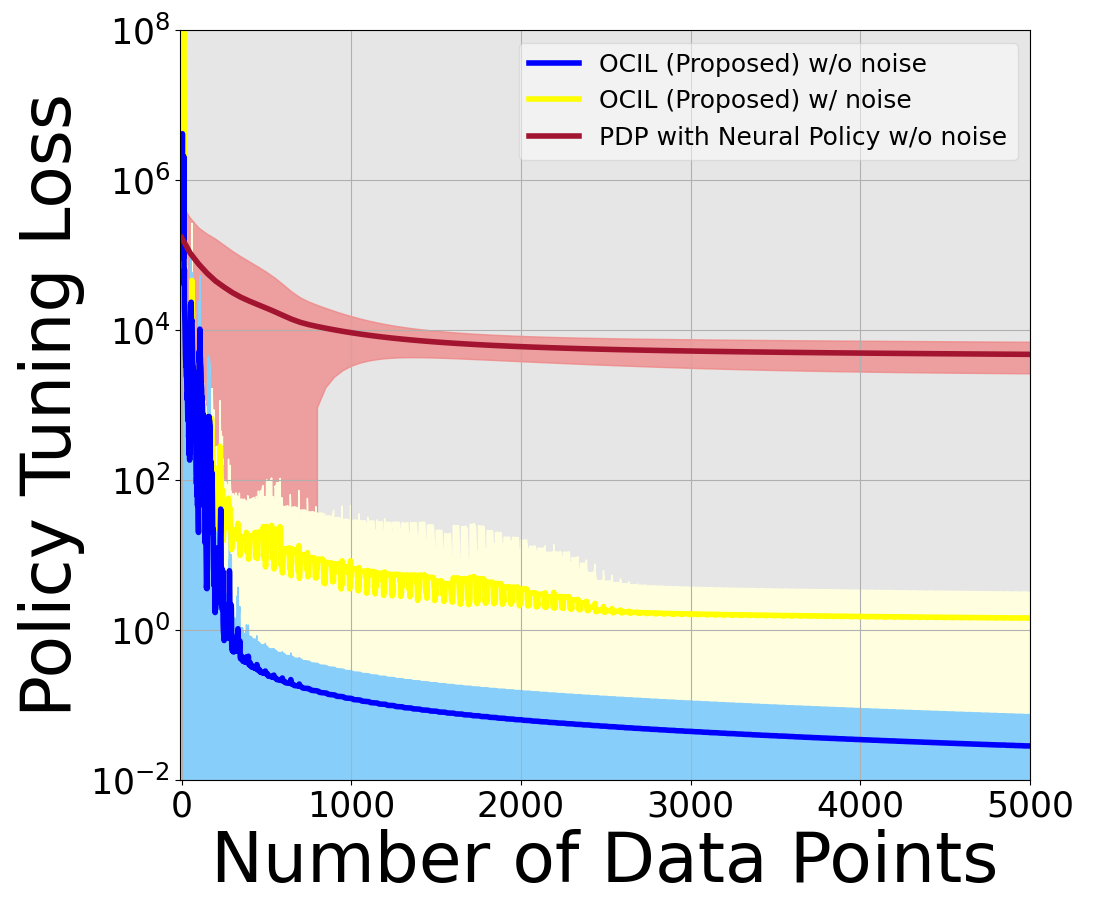}
         \caption{Rocket}
         \label{fig:PT_Rocket}
     \end{subfigure}
     \vskip 1pt%\baselineskip
     \centering
     \begin{subfigure}[b]{0.3\textwidth}
         \centering
         \includegraphics[width=\textwidth]{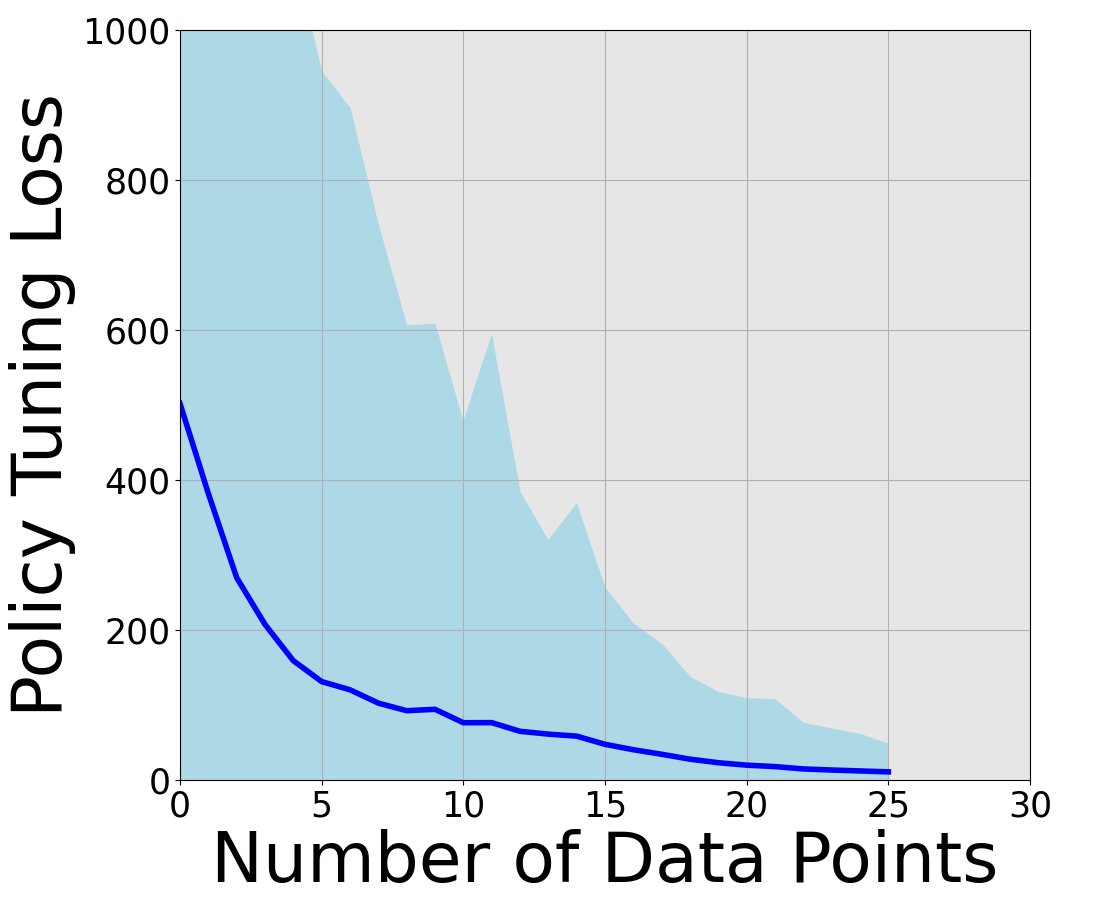}
         \caption{Cartpole, online phase, with noise}
         \label{fig:PT_Cartpole_online}
     \end{subfigure}
     \begin{subfigure}[b]{0.3\textwidth}
         \centering
         \includegraphics[width=\textwidth]{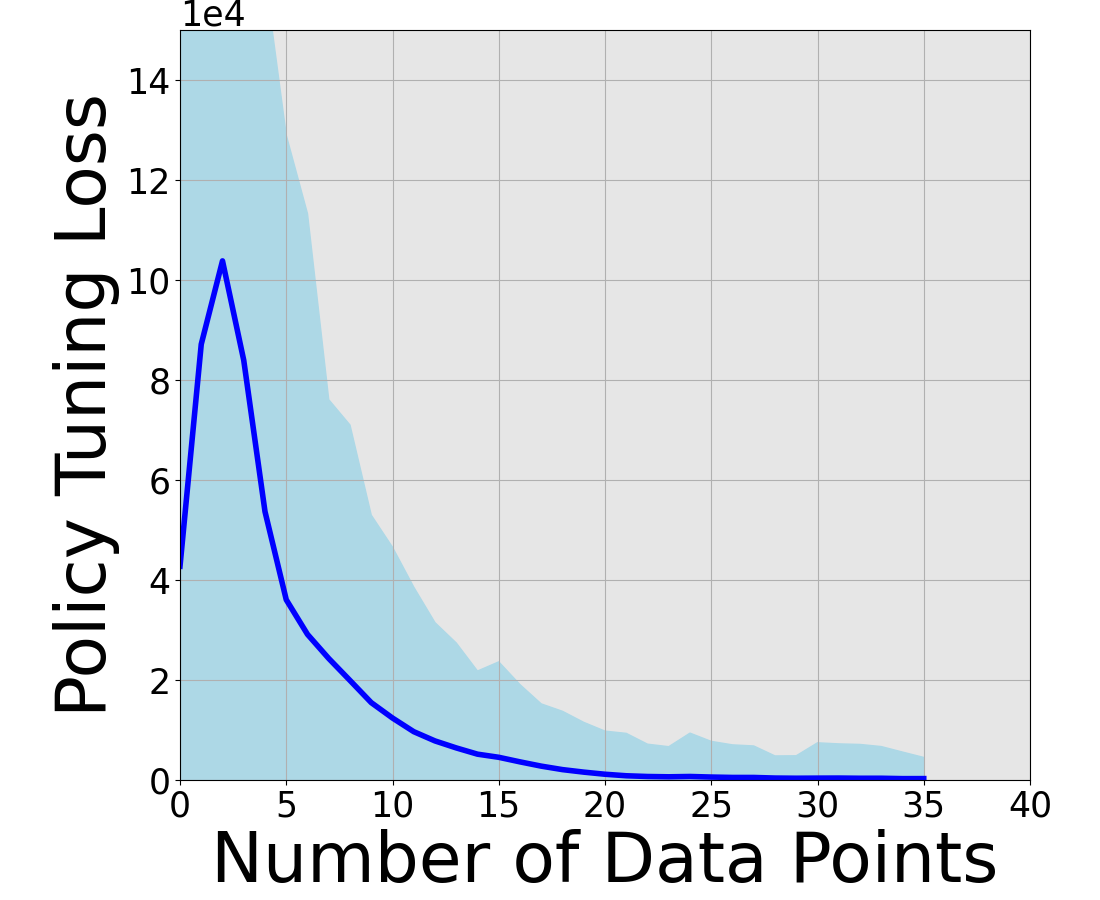}
         \caption{Quadrotor, online phase, with noise}
         \label{fig:PT_Quadrotor_online}
     \end{subfigure}
     \begin{subfigure}[b]{0.3\textwidth}
         \centering
         \includegraphics[width=\textwidth]{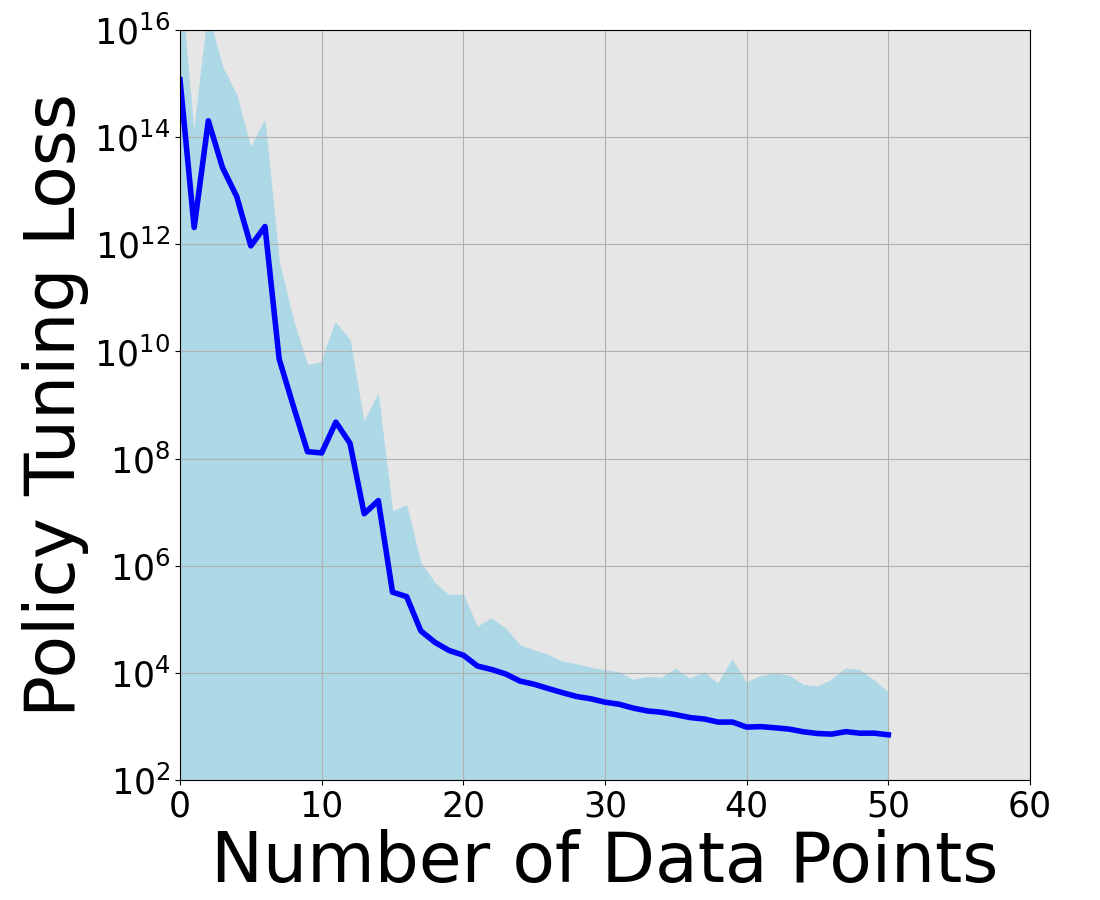}
         \caption{Rocket, online phase, with noise}
         \label{fig:PT_Rocket_online}
     \end{subfigure}
     \caption{Policy Tuning Loss v.s. number of data points. Buffers represent loss variation under 3$\sigma$ with random initial conditions.}
        \label{fig:PT}
\end{figure}

Fig.~\ref{fig:PT_Cartpole}-\ref{fig:PT_Rocket} summarize the result, where the loss and its variation of OCIL converge very quickly.
The buffers in Fig.~\ref{fig:PT_Cartpole}-\ref{fig:PT_Rocket_online} indicate 3 times of standard deviation.
Fig.~\ref{fig:PT_Cartpole_online}-\ref{fig:PT_Rocket_online} presents the online phase of OCIL given 1000 random trials, which further validates the effectiveness and robustness of OCIL given measurement noise. 

\begin{figure*}[ht]
\centering\includegraphics[width=0.32\linewidth]{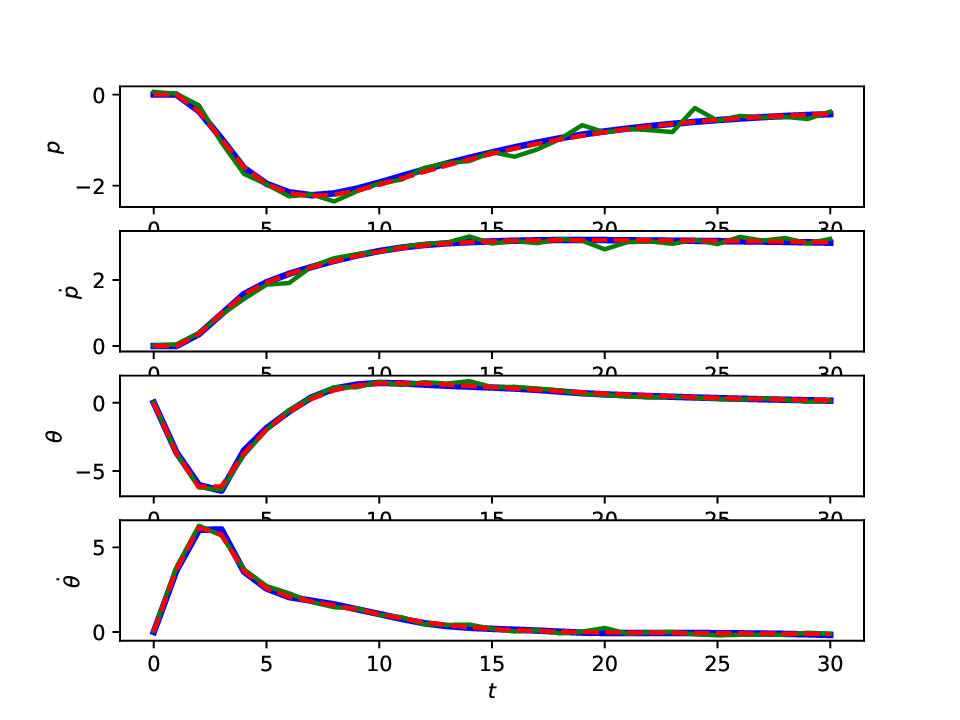}\label{fig:PT_cartpole_state}
\centering\includegraphics[width=0.32\linewidth]{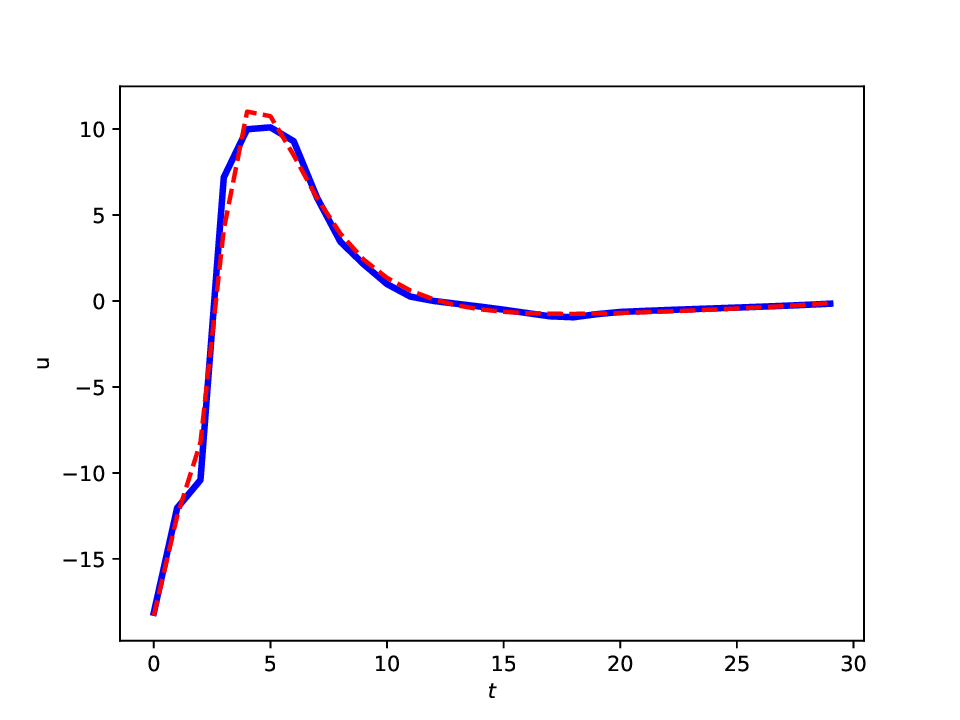}
% \hfill
\caption{Trajectories of the cartpole in policy tuning on-the-fly. Blue solid lines: learned trajectory. Green solid lines: observed noisy trajectory. Red dashed lines: ground truth.}\label{fig:PT_cartpole}
\end{figure*} 

\begin{figure*}[ht]
\centering\includegraphics[width=0.32\linewidth]{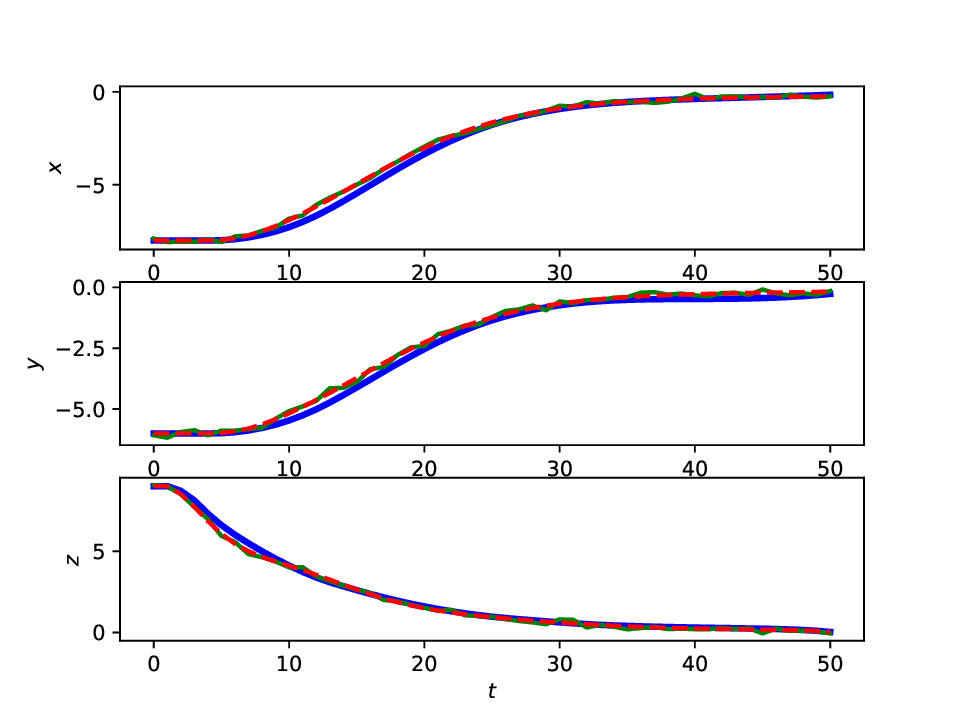}
\centering\includegraphics[width=0.32\linewidth]{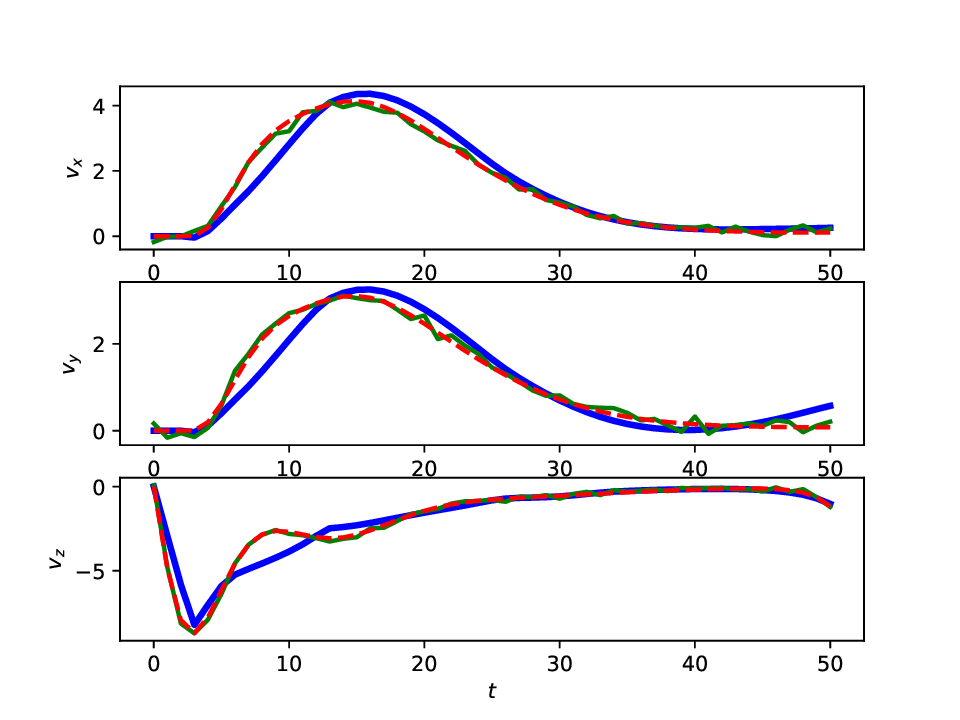}
% \hfill
\centering\includegraphics[width=0.32\linewidth]{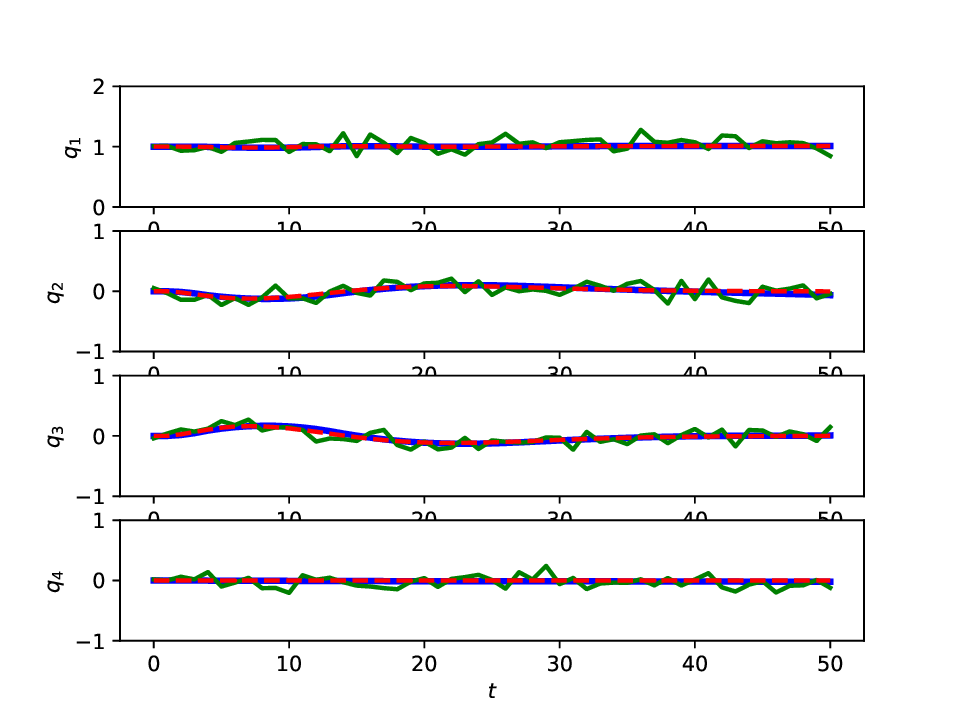}
\centering\includegraphics[width=0.32\linewidth]{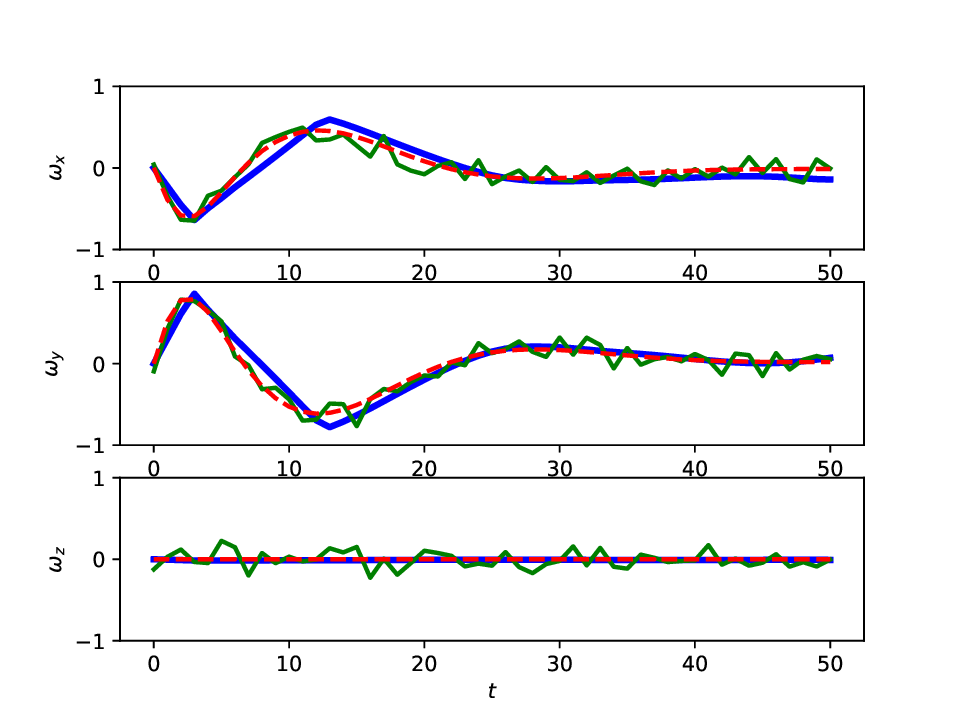}
\centering\includegraphics[width=0.32\linewidth]{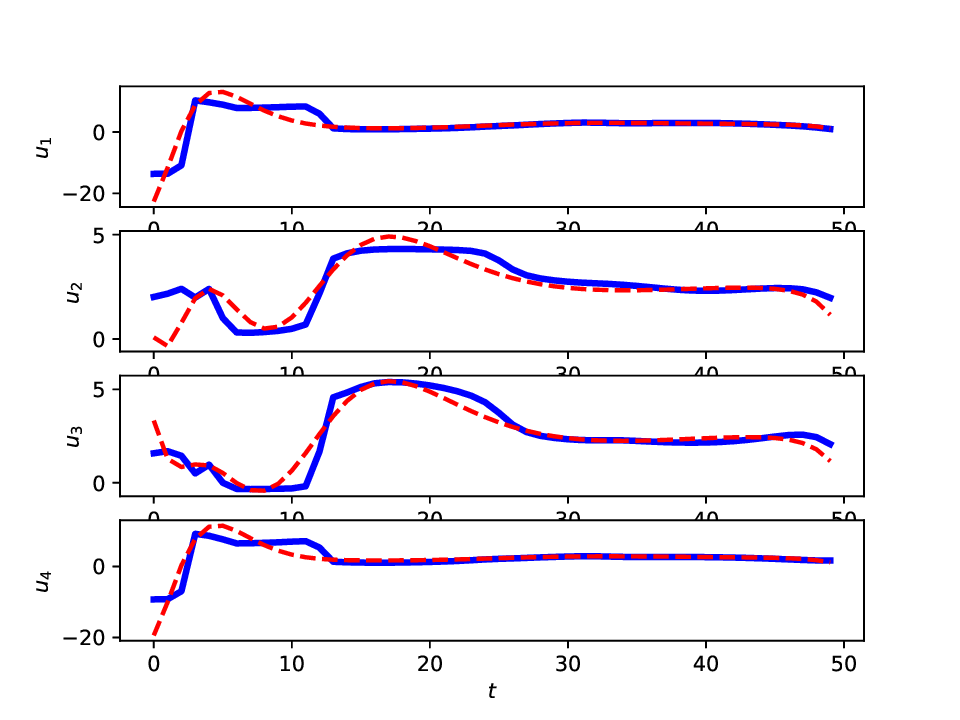}
\caption{Trajectories of the quadrotor in policy tuning on-the-fly. Blue solid lines: learned trajectory. Green solid lines: observed noisy trajectory. Red dashed lines: ground truth.} \label{fig:PT_uav}
\end{figure*} 

In general, OCIL from all figures does not have a smooth loss trajectory as the other offline methods. This is because at the online phase, an optimal gain matrix $\boldsymbol{K}_t$ from (\ref{eq:EKFpredict}) is computed to update $\hat{\boldsymbol{\theta}}_t$, whereas the other methods either use a constant or iteration-dependent step size. The optimal gain is conceptually similar to searching an optimal step size in the line-search optimization algorithms.
Thus, it is observed that the loss variation, as represented by blue buffers, is relatively high initially but starts decreasing significantly as new data comes in because $\boldsymbol{K}_t$ is continually updated. In contrast, the loss variation barely changes for the other offline methods after some data points.

% \vspace{-1em}
% \vspace{-1em}

\section{Online Computational Performance} \label{Appendix:computationPerformacnce}

The experiments with OCIL were performed on a desktop with one Intel Core i7-8700k CPU with 8GB RAM. No GPU was used.
The experiments with other methods were performed on a desktop with one AMD Ryzen 9 5900X CPU, one Nvidia Geforce RTX 4070ti, and 32 GB RAM.
A more powerful PC was selected for the other methods because of their high computational cost. As noted at the beginning of Section \ref{sec:applications}, only 5 trials were conducted due to the computational expense.

To demonstrate that the computational performance of OCIL is enough to be used in an online fashion, we recorded the computational time for OCIL in different modes for 100 trials. The code is implemented in Python, utilizing the CasADi library with the IPOPT solver to solve the OC problem.
% During the online phase, for each time step $t$, the average and standard deviation for the total elapsed time of Algorithm \ref{algorithm:OPDP} (include GG and the online parameter estimator) is $17.18\pm5.15$ ms. For the online parameter estimator, the average and standard deviation of elapsed time for each time step is $1.93\pm0.85$ ms. For the gradient generator, the average and standard deviation of elapsed time for each time step is $6.05\pm2.59$. $\Delta = 50ms$
Table~\ref{table:computationSysID} summarizes the OCIL's computational performance for the system identification task for three environments, where OCIL Time indicates the computational time of running OCIL at each time $t$, i.e. the iteration within the for-loop of Algorithm~\ref{algorithm:OPDP}; GG Time indicates the computational time of running gradient generator (GG) at each time $t$, i.e. Algorithm~\ref{algorithm:GG}; Estimator Time indicates the computational time of updating $\hat{\boldsymbol{\theta}}$, i.e. Line 7 of Algorithm~\ref{algorithm:OPDP}; $\Delta$ indicates the time step of each environment, i.e., the time duration between two consecutive data measurements or the maximum allowed time duration of online algorithms to perform computation; Percentage indicates the percentage of the average OCIL time with respect to $\Delta$. The header of Table~\ref{table:computationPT} and Table~\ref{table:computationIL} are the same.  Roughly speaking, OCIL time = GG time + Estimator time + Optimal Control computation time.

Table~\ref{table:computationSysID} illustrates that OCIL can estimate the dynamical system with neural network representation in an online fashion, within the system frequency of getting new data.
Table~\ref{table:computationPT} illustrates that OCIL can tune the neural policy online.
As indicated in Line 2 of Algorithm~\ref{algorithm:OPDP}, the most computationally heavy part is solving optimal control trajectory in an online fashion, instead of GG and the parameter estimator.
As demonstrated at the beginning of this section, OCIL does not require huge computational resources, such as GPU.
Therefore, OCIL has the capability to run in an online fashion.

% As can be seen, the core module of OCIL, i.e. GG and the online parameter estimator, does not induce a large computational cost. Most of the computational overhead is caused by generating the optimal trajectory using CASADI.

\begin{table}[ht]
\centering
\begin{threeparttable}
\caption{Computational Performance for Online Imitation Learning} \label{table:computationIL}
\begin{tabular}{c c c c c c}
\toprule
Env. & OCIL Time [ms] & GG Time [ms] & Estimator Time [ms] & $\Delta$ [ms] & Percentage \\
\midrule
Cartpole & $62.10\pm6.63$ & $7.47\pm0.25$ & $0.031\pm0.0023$ & $100$ & 62.10 \% \\
\midrule
Quadrotor & $81.70\pm2.51$ & $21.72\pm0.84$ & $0.058\pm0.039$ & $100$ & 81.70 \% \\
\midrule
Rocket & $72.25\pm13.91$ & $19.77\pm6.26$ & $0.060\pm0.012$ & $100$ &  72.25\% \\
\bottomrule
\end{tabular}
% \begin{tablenotes}
% \small
% \item $\dagger$ Percentage of OCIL time to time step $\Delta$
% \end{tablenotes}
\end{threeparttable}
\centering
\end{table}

\begin{table}[ht]
\centering
\begin{threeparttable}
\caption{Computational Performance for SysID with Neural System} \label{table:computationSysID}
\begin{tabular}{c c c c c c}
\toprule
Env. & OCIL Time [ms] & GG Time [ms] & Estimator Time [ms] & $\Delta$ [ms] & Percentage \\
\midrule
Cartpole & $17.18\pm5.15$ & $6.05\pm2.59$ & $1.93\pm0.85$ & $50$ & 34.36 \% \\
\midrule
Quadrotor & $35.53\pm8.98$ & $12.18\pm4.77$ & $16.18\pm6.11$ & $100$ & 35.53 \% \\
\midrule
Rocket & $29.54\pm8.29$ & $11.25\pm4.67$ & $12.89\pm5.29$ & $200$ & 14.77 \% \\
\bottomrule
\end{tabular}
% \begin{tablenotes}
% \small
% \item $\dagger$ Percentage of OCIL time to time step $\Delta$
% \end{tablenotes}
\end{threeparttable}
\centering
\end{table}

\begin{table}[ht]
\centering
\begin{threeparttable}
\caption{Computational Performance for Policy Tuning with Neural Policy} \label{table:computationPT}
\begin{tabular}{c c c c c c}
\toprule
Env. & OCIL Time [ms] & GG Time [ms] & Estimator Time [ms] & $\Delta$ [ms] & Percentage \\
\midrule
Cartpole & $16.41\pm5.35$ & $7.72\pm3.07$ & $3.96\pm1.91$ & $50$ & 32.82 \%  \\
\midrule
Quadrotor & $62.67\pm9.51$ & $30.86\pm2.25$ & $22.47\pm8.33$ & $100$ & 62.67 \%  \\
\midrule
Rocket & $59.02\pm7.25$ & $33.99\pm2.98$ & $12.94\pm5.62$ & $100$ & 59.02 \%  \\
\bottomrule
\end{tabular}
% \begin{tablenotes}
% \small
% \item $\dagger$ Percentage of the average OCIL time to time step $\Delta$
% \end{tablenotes}
\end{threeparttable}
\centering
\end{table}

% \section{Discussion} \label{sec:discussion}

% \textbf{Limitations.} \underline{i) Local convergence: }Since OCIL is based on first-order gradient, it can only achieve local minima for general non-convex optimization problems in \ref{eq:generalOpti0}. Global convergence of the bi-level optimization can be established via assumptions on the convexity and smoothness on the OC system. \underline{ii) Parameterization matters for global convergence: } When performing simulations, we find that how models are parameterized matters for good convergence performance. For example, in online SysID mode, we observe that using a neural network dynamics (in Fig. \ref{fig:SysID_Cartpole_NN}-\ref{fig:SysID_Rocket_NN}) is more likely to get trapped in local minima than using the true dynamics with unknown parameters (in Fig. \ref{fig:SysId_Cartpole}-\ref{fig:SysId_Rocket}). \underline{iii) Initialization matters: } As shown in Remark \ref{remark:initialGuess}, a bad initial guess might cause the estimator to diverge. For more detail of the discussion on limitations of the proposed method, please refer to Section \ref{Appendix:Limitations} in the Appendix.

\section{Limitations.} \label{Appendix:Limitations}

This section discusses the major limitations of the proposed framework from three perspectives.

\textbf{Local convergence: }Since OCIL is based on first-order gradients, it can only achieve local minima for general non-convex optimal control problems in (\ref{eq:est_auto_sys}). Furthermore, the general problem proposed in this paper belongs to a bi-level optimization framework. Under certain assumptions such as convexity and smoothness on models (e.g., dynamics model, policy, loss function, and control objective function), global convergence of the bi-level optimization can be established.
However, such conditions are too restrictive in the context of dynamic control systems. 
Therefore, the local convergence analysis based on general nonlinear optimization is enough.

\textbf{Parameterization matters for global convergence: } When performing experiments, we find that how models are parameterized matters for good convergence performance. For example, in online SysID mode, we observe that using a neural network dynamics (in Fig. \ref{fig:SysID_Cartpole_NN}-\ref{fig:SysID_Rocket_NN}) is more likely to get trapped in local minima than using the true dynamics with unknown parameters (in Fig. \ref{fig:SysId_Cartpole}-\ref{fig:SysId_Rocket})).
In general, more complex parameterization will bring extreme non-convexity to the optimization problem, making the algorithm more easily trapped in local minima. Determining the parameterization of an object to be learned requires prior or expert knowledge, which is common in ML.

\textbf{Initialization matters: } As OCIL borrows how optimal gain updates from EKF, they share the same drawback that convergence depends on the selection of initialization. As shown in Remark \ref{remark:initialGuess}, a bad initial guess might cause the estimator to diverge according to Lemma \ref{lemma:convergenceKalman}.
Therefore, if a relatively good initial guess is hard to retrieve, one might need to use other methods to cold start OCIL. 

\section{Conclusions}
This paper proposes Online Control-Informed Learning (OCIL), an online learning method tailored for diverse learning tasks. By considering an optimal control system with a tunable parameter, OCIL is a unified learning framework that effectively addresses tasks such as online imitation learning, online system identification, and tuning policy on-the-fly. By designing a signed residual function specific to each task and treating the parameter as a state of a new system, we employ the online parameter estimator to estimate the parameter online and minimize the signed residual at each time step. Theoretical analysis establishes the convergence conditions for OCIL, while experiments on various environments, tasks, and existing methods are done to validate its data efficiency, versatility, and robustness against measurement noise.

% Only add this information once your submission is accepted and deanonymized.

\subsubsection*{Acknowledgments}
This material is based upon work supported by the Office of Naval Research (ONR) and Saab, Inc. under the Threat and Situational Understanding of Networked Online Machine Intelligence (TSUNOMI) program (grant no. N00014-23-C-1016). Any opinions, findings and conclusions or recommendations expressed in this material are those of the author(s) and do not necessarily reflect the views of the ONR, the U.S. Government, or Saab, Inc.

% \bibliography{main}
\bibliographystyle{tmlr}
\bibliography{zihao, ref_zehui}

\begin{thebibliography}{96}
\providecommand{\natexlab}[1]{#1}
\providecommand{\url}[1]{\texttt{#1}}
\expandafter\ifx\csname urlstyle\endcsname\relax
  \providecommand{\doi}[1]{doi: #1}\else
  \providecommand{\doi}{doi: \begingroup \urlstyle{rm}\Url}\fi

\bibitem[Abadi et~al.(2016)Abadi, Agarwal, Barham, Brevdo, Chen, Citro,
  Corrado, Davis, Dean, Devin, et~al.]{abadi2016tensorflow}
Mart{\'\i}n Abadi, Ashish Agarwal, Paul Barham, Eugene Brevdo, Zhifeng Chen,
  Craig Citro, Greg~S Corrado, Andy Davis, Jeffrey Dean, Matthieu Devin, et~al.
\newblock Tensorflow: Large-scale machine learning on heterogeneous distributed
  systems.
\newblock \emph{arXiv preprint arXiv:1603.04467}, 2016.

\bibitem[Abbeel \& Ng(2004)Abbeel and Ng]{abbeel2004apprenticeship}
Pieter Abbeel and Andrew~Y Ng.
\newblock Apprenticeship learning via inverse reinforcement learning.
\newblock In \emph{International Conference on Machine Learning}, pp.\  1--8,
  2004.

\bibitem[Abbeel et~al.(2006)Abbeel, Quigley, and Ng]{abbeel2006using}
Pieter Abbeel, Morgan Quigley, and Andrew~Y Ng.
\newblock Using inaccurate models in reinforcement learning.
\newblock In \emph{International Conference on Machine Learning}, pp.\  1--8,
  2006.

\bibitem[Abraham \& Murphey(2019)Abraham and Murphey]{abraham2019active}
Ian Abraham and Todd~D Murphey.
\newblock Active learning of dynamics for data-driven control using koopman
  operators.
\newblock \emph{IEEE Transactions on Robotics}, 35\penalty0 (5):\penalty0
  1071--1083, 2019.

\bibitem[Amos et~al.(2018)Amos, Jimenez, Sacks, Boots, and
  Kolter]{amos2018differentiable}
Brandon Amos, Ivan Jimenez, Jacob Sacks, Byron Boots, and J~Zico Kolter.
\newblock Differentiable mpc for end-to-end planning and control.
\newblock In \emph{Advances in Neural Information Processing Systems}, pp.\
  8289--8300, 2018.

\bibitem[Arora \& Doshi(2021)Arora and Doshi]{arora2021survey}
Saurabh Arora and Prashant Doshi.
\newblock A survey of inverse reinforcement learning: Challenges, methods and
  progress.
\newblock \emph{Artificial Intelligence}, 297:\penalty0 103500, 2021.

\bibitem[Athans(1971)]{athans1971role}
Michael Athans.
\newblock The role and use of the stochastic linear-quadratic-gaussian problem
  in control system design.
\newblock \emph{IEEE Transactions on Automatic Control}, 16\penalty0
  (6):\penalty0 529--552, 1971.

\bibitem[Bastianello et~al.(2024)Bastianello, Carli, and Zampieri]{10189107}
Nicola Bastianello, Ruggero Carli, and Sandro Zampieri.
\newblock Internal model-based online optimization.
\newblock \emph{IEEE Transactions on Automatic Control}, 69\penalty0
  (1):\penalty0 689--696, 2024.

\bibitem[Beintema et~al.(2023)Beintema, Schoukens, and
  T{\'o}th]{beintema2023deep}
Gerben~I Beintema, Maarten Schoukens, and Roland T{\'o}th.
\newblock Deep subspace encoders for nonlinear system identification.
\newblock \emph{Automatica}, 156:\penalty0 111210, 2023.

\bibitem[Benning et~al.(2019)Benning, Celledoni, Ehrhardt, Owren, and
  Sch{\"o}nlieb]{benning2019deep}
Martin Benning, Elena Celledoni, Matthias~J Ehrhardt, Brynjulf Owren, and
  Carola-Bibiane Sch{\"o}nlieb.
\newblock Deep learning as optimal control problems: models and numerical
  methods.
\newblock \emph{arXiv preprint arXiv:1904.05657}, 2019.

\bibitem[Berkenkamp et~al.(2023)Berkenkamp, Krause, and
  Schoellig]{berkenkamp2023bayesian}
Felix Berkenkamp, Andreas Krause, and Angela~P Schoellig.
\newblock Bayesian optimization with safety constraints: safe and automatic
  parameter tuning in robotics.
\newblock \emph{Machine Learning}, 112\penalty0 (10):\penalty0 3713--3747,
  2023.

\bibitem[Bertsekas(2022)]{bertsekas2022lessons}
Dimitri Bertsekas.
\newblock \emph{Lessons from AlphaZero for optimal, model predictive, and
  adaptive control}.
\newblock Athena Scientific, 2022.

\bibitem[Bertsekas(1996)]{1996Bertsekas}
Dimitri~P. Bertsekas.
\newblock Incremental least squares methods and the extended kalman filter.
\newblock \emph{SIAM Journal on Optimization}, 6\penalty0 (3):\penalty0
  807--822, 1996.

\bibitem[Bock \& Plitt(1984)Bock and Plitt]{bock1984multiple}
Hans~Georg Bock and Karl-Josef Plitt.
\newblock A multiple shooting algorithm for direct solution of optimal control
  problems.
\newblock \emph{IFAC Proceedings Volumes}, 17\penalty0 (2):\penalty0
  1603--1608, 1984.

\bibitem[Bojarski et~al.(2016)Bojarski, Del~Testa, Dworakowski, Firner, Flepp,
  Goyal, Jackel, Monfort, Muller, Zhang, et~al.]{bojarski2016end}
Mariusz Bojarski, Davide Del~Testa, Daniel Dworakowski, Bernhard Firner, Beat
  Flepp, Prasoon Goyal, Lawrence~D Jackel, Mathew Monfort, Urs Muller, Jiakai
  Zhang, et~al.
\newblock End to end learning for self-driving cars.
\newblock \emph{arXiv preprint arXiv:1604.07316}, 2016.

\bibitem[B{\"o}ttcher et~al.(2022)B{\"o}ttcher, Antulov-Fantulin, and
  Asikis]{bottcher2022ai}
Lucas B{\"o}ttcher, Nino Antulov-Fantulin, and Thomas Asikis.
\newblock Ai pontryagin or how artificial neural networks learn to control
  dynamical systems.
\newblock \emph{Nature Communications}, 13\penalty0 (1):\penalty0 333, 2022.

\bibitem[Boutayeb \& Aubry(1999)Boutayeb and Aubry]{Boutayeb1999}
M.~Boutayeb and D.~Aubry.
\newblock A strong tracking extended kalman observer for nonlinear
  discrete-time systems.
\newblock \emph{IEEE Transactions on Automatic Control}, 44\penalty0
  (8):\penalty0 1550--1556, 1999.

\bibitem[Boutayeb et~al.(1997)Boutayeb, Rafaralahy, and Darouach]{Boutayeb1997}
M.~Boutayeb, H.~Rafaralahy, and M.~Darouach.
\newblock Convergence analysis of the extended kalman filter used as an
  observer for nonlinear deterministic discrete-time systems.
\newblock \emph{IEEE Transactions on Automatic Control}, 42\penalty0
  (4):\penalty0 581--586, 1997.

\bibitem[Brown et~al.(2019)Brown, Goo, Nagarajan, and
  Niekum]{brown2019extrapolating}
Daniel Brown, Wonjoon Goo, Prabhat Nagarajan, and Scott Niekum.
\newblock Extrapolating beyond suboptimal demonstrations via inverse
  reinforcement learning from observations.
\newblock In \emph{International Conference on Machine Learning}, pp.\
  783--792. PMLR, 2019.

\bibitem[Bryson(1975)]{bryson1975applied}
Arthur~Earl Bryson.
\newblock \emph{Applied optimal control: optimization, estimation and control}.
\newblock CRC Press, 1975.

\bibitem[Casti et~al.(2023)Casti, Bastianello, Carli, and
  Zampieri]{casti2023control}
Umberto Casti, Nicola Bastianello, Ruggero Carli, and Sandro Zampieri.
\newblock A control theoretical approach to online constrained optimization.
\newblock \emph{arXiv preprint arXiv:2309.15498}, 2023.

\bibitem[Chan \& van~der Schaar(2021)Chan and van~der Schaar]{chan2021scalable}
Alex~J Chan and Mihaela van~der Schaar.
\newblock Scalable bayesian inverse reinforcement learning.
\newblock \emph{arXiv preprint arXiv:2102.06483}, 2021.

\bibitem[Czarnecki et~al.(2019)Czarnecki, Pascanu, Osindero, Jayakumar,
  Swirszcz, and Jaderberg]{czarnecki2019distilling}
Wojciech~M Czarnecki, Razvan Pascanu, Simon Osindero, Siddhant Jayakumar,
  Grzegorz Swirszcz, and Max Jaderberg.
\newblock Distilling policy distillation.
\newblock In \emph{The 22nd International Conference on Artificial Intelligence
  and Statistics}, pp.\  1331--1340. PMLR, 2019.

\bibitem[Deisenroth \& Rasmussen(2011)Deisenroth and
  Rasmussen]{deisenroth2011pilco}
Marc Deisenroth and Carl~E Rasmussen.
\newblock Pilco: A model-based and data-efficient approach to policy search.
\newblock In \emph{Proceedings of the 28th International Conference on Machine
  Learning (ICML-11)}, pp.\  465--472, 2011.

\bibitem[Englert et~al.(2017)Englert, Vien, and Toussaint]{englert2017inverse}
Peter Englert, Ngo~Anh Vien, and Marc Toussaint.
\newblock Inverse kkt: Learning cost functions of manipulation tasks from
  demonstrations.
\newblock \emph{The International Journal of Robotics Research}, 36\penalty0
  (13-14):\penalty0 1474--1488, 2017.

\bibitem[Finn et~al.(2016)Finn, Goodfellow, and Levine]{finn2016unsupervised}
Chelsea Finn, Ian Goodfellow, and Sergey Levine.
\newblock Unsupervised learning for physical interaction through video
  prediction.
\newblock In \emph{Advances in Neural Information Processing Systems}, pp.\
  64--72, 2016.

\bibitem[Ghaemi et~al.(2009)Ghaemi, Sun, and
  Kolmanovsky]{ghaemi2009neighboring}
Reza Ghaemi, Jing Sun, and Ilya~V Kolmanovsky.
\newblock Neighboring extremal solution for nonlinear discrete-time optimal
  control problems with state inequality constraints.
\newblock \emph{IEEE Transactions on Automatic Control}, 54\penalty0
  (11):\penalty0 2674--2679, 2009.

\bibitem[Gu et~al.(2016)Gu, Lillicrap, Sutskever, and Levine]{gu2016continuous}
Shixiang Gu, Timothy Lillicrap, Ilya Sutskever, and Sergey Levine.
\newblock Continuous deep q-learning with model-based acceleration.
\newblock In \emph{International Conference on Machine Learning}, pp.\
  2829--2838, 2016.

\bibitem[Guo \& Pan(2023)Guo and Pan]{guo2023composite}
Kai Guo and Yongping Pan.
\newblock Composite adaptation and learning for robot control: A survey.
\newblock \emph{Annual Reviews in Control}, 55:\penalty0 279--290, 2023.

\bibitem[Han et~al.(2019)Han, Li, et~al.]{han2019mean}
Jiequn Han, Qianxiao Li, et~al.
\newblock A mean-field optimal control formulation of deep learning.
\newblock \emph{Research in the Mathematical Sciences}, 6\penalty0
  (1):\penalty0 1--41, 2019.

\bibitem[Hao et~al.(2023)Hao, Heredia, Huang, Lu, Liang, and
  Mou]{hao2023policy}
Wenjian Hao, Paulo~C Heredia, Bowen Huang, Zehui Lu, Zihao Liang, and Shaoshuai
  Mou.
\newblock Policy learning based on deep koopman representation.
\newblock \emph{arXiv preprint arXiv:2305.15188}, 2023.

\bibitem[Hao et~al.(2024)Hao, Lu, Upadhyay, and Mou]{hao2024distributed}
Wenjian Hao, Zehui Lu, Devesh Upadhyay, and Shaoshuai Mou.
\newblock A distributed deep koopman learning algorithm for control.
\newblock \emph{arXiv preprint arXiv:2412.07212}, 2024.

\bibitem[Haruno et~al.(2001)Haruno, Wolpert, and Kawato]{haruno2001mosaic}
Masahiko Haruno, Daniel~M Wolpert, and Mitsuo Kawato.
\newblock Mosaic model for sensorimotor learning and control.
\newblock \emph{Neural Computation}, 13\penalty0 (10):\penalty0 2201--2220,
  2001.

\bibitem[Heess et~al.(2015)Heess, Wayne, Silver, Lillicrap, Erez, and
  Tassa]{heess2015learning}
Nicolas Heess, Gregory Wayne, David Silver, Timothy Lillicrap, Tom Erez, and
  Yuval Tassa.
\newblock Learning continuous control policies by stochastic value gradients.
\newblock In \emph{Advances in Neural Information Processing Systems}, pp.\
  2944--2952, 2015.

\bibitem[Ioannou \& Sun(2012)Ioannou and Sun]{ioannou2012robust}
Petros~A Ioannou and Jing Sun.
\newblock \emph{Robust adaptive control}.
\newblock Courier Corporation, 2012.

\bibitem[Jin \& Mou(2021)Jin and Mou]{jin2021distributed}
Wanxin Jin and Shaoshuai Mou.
\newblock Distributed inverse optimal control.
\newblock \emph{Automatica}, 129:\penalty0 109658, 2021.
\newblock ISSN 0005-1098.

\bibitem[Jin et~al.(2019)Jin, Kuli{\'c}, Lin, Mou, and Hirche]{jin2019inverse}
Wanxin Jin, Dana Kuli{\'c}, Jonathan Feng-Shun Lin, Shaoshuai Mou, and Sandra
  Hirche.
\newblock Inverse optimal control for multiphase cost functions.
\newblock \emph{IEEE Transactions on Robotics}, 35\penalty0 (6):\penalty0
  1387--1398, 2019.

\bibitem[Jin et~al.(2020)Jin, Wang, Yang, and Mou]{jin2020pontryagin}
Wanxin Jin, Zhaoran Wang, Zhuoran Yang, and Shaoshuai Mou.
\newblock Pontryagin differentiable programming: An end-to-end learning and
  control framework.
\newblock \emph{Advances in Neural Information Processing Systems},
  33:\penalty0 7979--7992, 2020.

\bibitem[Jin et~al.(2021{\natexlab{a}})Jin, Kuli{\'c}, Mou, and
  Hirche]{jin2021inverse}
Wanxin Jin, Dana Kuli{\'c}, Shaoshuai Mou, and Sandra Hirche.
\newblock Inverse optimal control from incomplete trajectory observations.
\newblock \emph{The International Journal of Robotics Research}, 40\penalty0
  (6-7):\penalty0 848--865, 2021{\natexlab{a}}.

\bibitem[Jin et~al.(2021{\natexlab{b}})Jin, Mou, and Pappas]{jin2021safe}
Wanxin Jin, Shaoshuai Mou, and George~J Pappas.
\newblock Safe pontryagin differentiable programming.
\newblock \emph{Advances in Neural Information Processing Systems},
  34:\penalty0 16034--16050, 2021{\natexlab{b}}.

\bibitem[Kappen(2005)]{kappen2005path}
Hilbert~J Kappen.
\newblock Path integrals and symmetry breaking for optimal control theory.
\newblock \emph{Journal of Statistical Mechanics: Theory and Experiment}, 2005.

\bibitem[Karniadakis et~al.(2021)Karniadakis, Kevrekidis, Lu, Perdikaris, Wang,
  and Yang]{karniadakis2021physics}
George~Em Karniadakis, Ioannis~G Kevrekidis, Lu~Lu, Paris Perdikaris, Sifan
  Wang, and Liu Yang.
\newblock Physics-informed machine learning.
\newblock \emph{Nature Reviews Physics}, 3\penalty0 (6):\penalty0 422--440,
  2021.

\bibitem[Kashinath et~al.(2021)Kashinath, Mustafa, Albert, Wu, Jiang,
  Esmaeilzadeh, Azizzadenesheli, Wang, Chattopadhyay, Singh,
  et~al.]{kashinath2021physics}
Karthik Kashinath, M~Mustafa, Adrian Albert, JL~Wu, C~Jiang, Soheil
  Esmaeilzadeh, Kamyar Azizzadenesheli, R~Wang, Ashesh Chattopadhyay, A~Singh,
  et~al.
\newblock Physics-informed machine learning: case studies for weather and
  climate modelling.
\newblock \emph{Philosophical Transactions of the Royal Society A},
  379\penalty0 (2194):\penalty0 20200093, 2021.

\bibitem[Keesman(2011)]{keesman2011system}
Karel~J Keesman.
\newblock \emph{System identification: an introduction}.
\newblock Springer Science \& Business Media, 2011.

\bibitem[Keshavarz et~al.(2011)Keshavarz, Wang, and
  Boyd]{keshavarz2011imputing}
Arezou Keshavarz, Yang Wang, and Stephen Boyd.
\newblock Imputing a convex objective function.
\newblock In \emph{IEEE International Symposium on Intelligent Control}, pp.\
  613--619, 2011.

\bibitem[Khosravi et~al.(2021)Khosravi, Behrunani, Myszkorowski, Smith,
  Rupenyan, and Lygeros]{khosravi2021performance}
Mohammad Khosravi, Varsha~N Behrunani, Piotr Myszkorowski, Roy~S Smith, Alisa
  Rupenyan, and John Lygeros.
\newblock Performance-driven cascade controller tuning with bayesian
  optimization.
\newblock \emph{IEEE Transactions on Industrial Electronics}, 69\penalty0
  (1):\penalty0 1032--1042, 2021.

\bibitem[Kuipers(1999)]{kuipers1999quaternions}
Jack~B Kuipers.
\newblock \emph{Quaternions and rotation sequences}, volume~66.
\newblock Princeton University Press, 1999.

\bibitem[Kumpati et~al.(1990)Kumpati, Kannan,
  et~al.]{kumpati1990identification}
S~Narendra Kumpati, Parthasarathy Kannan, et~al.
\newblock Identification and control of dynamical systems using neural
  networks.
\newblock \emph{IEEE Transactions on Neural Networks}, 1\penalty0 (1):\penalty0
  4--27, 1990.

\bibitem[LeCun et~al.(2015)LeCun, Bengio, and Hinton]{lecun2015deep}
Yann LeCun, Yoshua Bengio, and Geoffrey Hinton.
\newblock Deep learning.
\newblock \emph{Nature}, 521\penalty0 (7553):\penalty0 436--444, 2015.

\bibitem[Levine \& Abbeel(2014)Levine and Abbeel]{levine2014learning}
Sergey Levine and Pieter Abbeel.
\newblock Learning neural network policies with guided policy search under
  unknown dynamics.
\newblock In \emph{Advances in Neural Information Processing Systems}, pp.\
  1071--1079, 2014.

\bibitem[Lewis et~al.(1998)Lewis, Jagannathan, and Yesildirak]{lewis1998neural}
FW~Lewis, Suresh Jagannathan, and Aydin Yesildirak.
\newblock \emph{Neural network control of robot manipulators and non-linear
  systems}.
\newblock CRC press, 1998.

\bibitem[Li \& Hao(2018)Li and Hao]{li2018optimal}
Qianxiao Li and Shuji Hao.
\newblock An optimal control approach to deep learning and applications to
  discrete-weight neural networks.
\newblock \emph{arXiv preprint arXiv:1803.01299}, 2018.

\bibitem[Li et~al.(2018)Li, Chen, Tai, and Weinan]{li2018maximum}
Qianxiao Li, Long Chen, Cheng Tai, and E~Weinan.
\newblock Maximum principle based algorithms for deep learning.
\newblock \emph{Journal of Machine Learning Research}, 18\penalty0
  (165):\penalty0 1--29, 2018.

\bibitem[Li \& Todorov(2004)Li and Todorov]{li2004iterative}
Weiwei Li and Emanuel Todorov.
\newblock Iterative linear quadratic regulator design for nonlinear biological
  movement systems.
\newblock In \emph{International Conference on Informatics in Control,
  Automation and Robotics}, pp.\  222--229, 2004.

\bibitem[Liang et~al.(2022)Liang, Jin, and Mou]{liang2022iterative}
Zihao Liang, Wanxin Jin, and Shaoshuai Mou.
\newblock An iterative method for inverse optimal control.
\newblock In \emph{2022 13th Asian Control Conference (ASCC)}, pp.\  959--964,
  2022.

\bibitem[Liang et~al.(2023)Liang, Hao, and Mou]{liang2023data}
Zihao Liang, Wenjian Hao, and Shaoshuai Mou.
\newblock A data-driven approach for inverse optimal control.
\newblock In \emph{2023 62nd IEEE Conference on Decision and Control (CDC)},
  pp.\  3632--3637, 2023.

\bibitem[Liu \& Markowich(2020)Liu and Markowich]{liu2020selection}
Hailiang Liu and Peter Markowich.
\newblock Selection dynamics for deep neural networks.
\newblock \emph{Journal of Differential Equations}, 269\penalty0 (12):\penalty0
  11540--11574, 2020.

\bibitem[Liu et~al.(2024)Liu, Li, Wang, and Long]{liu2024koopa}
Yong Liu, Chenyu Li, Jianmin Wang, and Mingsheng Long.
\newblock Koopa: Learning non-stationary time series dynamics with koopman
  predictors.
\newblock \emph{Advances in Neural Information Processing Systems}, 36, 2024.

\bibitem[Lu et~al.(2024)Lu, Zhang, and Wang]{10771962}
Zehui Lu, Tianpeng Zhang, and Yebin Wang.
\newblock Torque constraint modeling and reference shaping for servo systems.
\newblock \emph{IEEE Control Systems Letters}, 8:\penalty0 2637--2642, 2024.

\bibitem[Luo et~al.(2023)Luo, Peng, Hu, and Ghosh]{luo2023adaptive}
Rui Luo, Zhinan Peng, Jiangping Hu, and Bijoy~Kumar Ghosh.
\newblock Adaptive optimal control of affine nonlinear systems via
  identifier--critic neural network approximation with relaxed pe conditions.
\newblock \emph{Neural Networks}, 167:\penalty0 588--600, 2023.

\bibitem[Lutter et~al.(2019)Lutter, Ritter, and Peters]{lutter2019deep}
Michael Lutter, Christian Ritter, and Jan Peters.
\newblock Deep lagrangian networks: Using physics as model prior for deep
  learning.
\newblock \emph{arXiv preprint arXiv:1907.04490}, 2019.

\bibitem[Mauroy et~al.(2020)Mauroy, Susuki, and Mezi{\'c}]{mauroy2020koopman}
Alexandre Mauroy, Y~Susuki, and Igor Mezi{\'c}.
\newblock \emph{Koopman operator in systems and control}.
\newblock Springer, 2020.

\bibitem[Mnih et~al.(2013)Mnih, Kavukcuoglu, Silver, Graves, Antonoglou,
  Wierstra, and Riedmiller]{mnih2013playing}
Volodymyr Mnih, Koray Kavukcuoglu, David Silver, Alex Graves, Ioannis
  Antonoglou, Daan Wierstra, and Martin Riedmiller.
\newblock Playing atari with deep reinforcement learning.
\newblock \emph{arXiv preprint arXiv:1312.5602}, 2013.

\bibitem[Mnih et~al.(2015)Mnih, Kavukcuoglu, Silver, Rusu, Veness, Bellemare,
  Graves, Riedmiller, Fidjeland, Ostrovski, et~al.]{mnih2015human}
Volodymyr Mnih, Koray Kavukcuoglu, David Silver, Andrei~A Rusu, Joel Veness,
  Marc~G Bellemare, Alex Graves, Martin Riedmiller, Andreas~K Fidjeland, Georg
  Ostrovski, et~al.
\newblock Human-level control through deep reinforcement learning.
\newblock \emph{Nature}, 518\penalty0 (7540):\penalty0 529--533, 2015.

\bibitem[Mombaur et~al.(2010)Mombaur, Truong, and Laumond]{mombaur2010human}
Katja Mombaur, Anh Truong, and Jean-Paul Laumond.
\newblock From human to humanoid locomotion—an inverse optimal control
  approach.
\newblock \emph{Autonomous Robots}, 28\penalty0 (3):\penalty0 369--383, 2010.

\bibitem[Nelles \& Nelles(2020)Nelles and Nelles]{nelles2020nonlinear}
Oliver Nelles and Oliver Nelles.
\newblock \emph{Nonlinear dynamic system identification}.
\newblock Springer, 2020.

\bibitem[Oh et~al.(2016)Oh, Chockalingam, Singh, and Lee]{oh2016control}
Junhyuk Oh, Valliappa Chockalingam, Satinder Singh, and Honglak Lee.
\newblock Control of memory, active perception, and action in minecraft.
\newblock \emph{arXiv preprint arXiv:1605.09128}, 2016.

\bibitem[Okada et~al.(2017)Okada, Rigazio, and Aoshima]{okada2017path}
Masashi Okada, Luca Rigazio, and Takenobu Aoshima.
\newblock Path integral networks: End-to-end differentiable optimal control.
\newblock \emph{arXiv preprint arXiv:1706.09597}, 2017.

\bibitem[Patterson \& Rao(2014)Patterson and Rao]{patterson2014gpops}
Michael~A Patterson and Anil~V Rao.
\newblock Gpops-ii: A matlab software for solving multiple-phase optimal
  control problems using hp-adaptive gaussian quadrature collocation methods
  and sparse nonlinear programming.
\newblock \emph{ACM Transactions on Mathematical Software}, 41\penalty0
  (1):\penalty0 1, 2014.

\bibitem[Pereira et~al.(2018)Pereira, Fan, An, and Theodorou]{pereira2018mpc}
Marcus Pereira, David~D Fan, Gabriel~Nakajima An, and Evangelos Theodorou.
\newblock Mpc-inspired neural network policies for sequential decision making.
\newblock \emph{arXiv preprint arXiv:1802.05803}, 2018.

\bibitem[Pillonetto et~al.(2025)Pillonetto, Aravkin, Gedon, Ljung, Ribeiro, and
  Sch{\"o}n]{pillonetto2025deep}
Gianluigi Pillonetto, Aleksandr Aravkin, Daniel Gedon, Lennart Ljung,
  Ant{\^o}nio~H Ribeiro, and Thomas~B Sch{\"o}n.
\newblock Deep networks for system identification: a survey.
\newblock \emph{Automatica}, 171:\penalty0 111907, 2025.

\bibitem[Pontryagin(2018)]{pontryagin2018mathematical}
Lev~Semenovich Pontryagin.
\newblock \emph{Mathematical theory of optimal processes}.
\newblock Routledge, 2018.

\bibitem[Proctor et~al.(2016)Proctor, Brunton, and Kutz]{proctor2016dynamic}
Joshua~L Proctor, Steven~L Brunton, and J~Nathan Kutz.
\newblock Dynamic mode decomposition with control.
\newblock \emph{SIAM Journal on Applied Dynamical Systems}, 15\penalty0
  (1):\penalty0 142--161, 2016.

\bibitem[Proctor et~al.(2018)Proctor, Brunton, and
  Kutz]{proctor2018generalizing}
Joshua~L Proctor, Steven~L Brunton, and J~Nathan Kutz.
\newblock Generalizing koopman theory to allow for inputs and control.
\newblock \emph{SIAM Journal on Applied Dynamical Systems}, 17\penalty0
  (1):\penalty0 909--930, 2018.

\bibitem[Ratliff et~al.(2006)Ratliff, Bagnell, and
  Zinkevich]{ratliff2006maximum}
Nathan~D Ratliff, J~Andrew Bagnell, and Martin~A Zinkevich.
\newblock Maximum margin planning.
\newblock In \emph{International Conference on Machine Learning}, pp.\
  729--736, 2006.

\bibitem[Revach et~al.(2022)Revach, Shlezinger, Ni, Escoriza, van Sloun, and
  Eldar]{revach2022kalman}
Guy Revach, Nir Shlezinger, Xiaoyong Ni, Adrià~López Escoriza, Ruud J.~G. van
  Sloun, and Yonina~C. Eldar.
\newblock Kalmannet: Neural network aided kalman filtering for partially known
  dynamics.
\newblock \emph{IEEE Transactions on Signal Processing}, 70:\penalty0
  1532--1547, 2022.

\bibitem[Ribeiro(2004)]{ribeiro2004kalman}
Maria~Isabel Ribeiro.
\newblock Kalman and extended kalman filters: Concept, derivation and
  properties.
\newblock \emph{Institute for Systems and Robotics}, 43\penalty0 (46):\penalty0
  3736--3741, 2004.

\bibitem[Saemundsson et~al.(2020)Saemundsson, Terenin, Hofmann, and
  Deisenroth]{saemundsson2020variational}
Steindor Saemundsson, Alexander Terenin, Katja Hofmann, and Marc Deisenroth.
\newblock Variational integrator networks for physically structured embeddings.
\newblock In \emph{International Conference on Artificial Intelligence and
  Statistics}, pp.\  3078--3087, 2020.

\bibitem[Sasaki \& Yamashina(2021)Sasaki and Yamashina]{sasaki2021behavioral}
Fumihiro Sasaki and Ryota Yamashina.
\newblock Behavioral cloning from noisy demonstrations.
\newblock In \emph{International Conference on Learning Representations}, 2021.

\bibitem[Schneider(1997)]{schneider1997exploiting}
Jeff~G Schneider.
\newblock Exploiting model uncertainty estimates for safe dynamic control
  learning.
\newblock In \emph{Advances in Neural Information Processing Systems}, pp.\
  1047--1053, 1997.

\bibitem[Schwenzer et~al.(2021)Schwenzer, Ay, Bergs, and
  Abel]{schwenzer2021review}
Max Schwenzer, Muzaffer Ay, Thomas Bergs, and Dirk Abel.
\newblock Review on model predictive control: An engineering perspective.
\newblock \emph{The International Journal of Advanced Manufacturing
  Technology}, 117\penalty0 (5):\penalty0 1327--1349, 2021.

\bibitem[Scokaert \& Rawlings(1998)Scokaert and
  Rawlings]{scokaert1998constrained}
Pierre~OM Scokaert and James~B Rawlings.
\newblock Constrained linear quadratic regulation.
\newblock \emph{IEEE Transactions on Automatic Control}, 43\penalty0
  (8):\penalty0 1163--1169, 1998.

\bibitem[Sharma et~al.(2023)Sharma, Chung, Akoush, and Ihme]{sharma2023review}
Pushan Sharma, Wai~Tong Chung, Bassem Akoush, and Matthias Ihme.
\newblock A review of physics-informed machine learning in fluid mechanics.
\newblock \emph{Energies}, 16\penalty0 (5):\penalty0 2343, 2023.

\bibitem[Song \& Grizzle(1992)Song and Grizzle]{Song1992}
Yongkyu Song and Jessy~W. Grizzle.
\newblock The extended kalman filter as a local asymptotic observer for
  nonlinear discrete-time systems.
\newblock In \emph{1992 American Control Conference}, pp.\  3365--3369, 1992.

\bibitem[Sorourifar et~al.(2021)Sorourifar, Makrygirgos, Mesbah, and
  Paulson]{sorourifar2021data}
Farshud Sorourifar, Georgios Makrygirgos, Ali Mesbah, and Joel~A Paulson.
\newblock A data-driven automatic tuning method for mpc under uncertainty using
  constrained bayesian optimization.
\newblock \emph{IFAC-PapersOnLine}, 54\penalty0 (3):\penalty0 243--250, 2021.

\bibitem[Srinivas et~al.(2018)Srinivas, Jabri, Abbeel, Levine, and
  Finn]{srinivas2018universal}
Aravind Srinivas, Allan Jabri, Pieter Abbeel, Sergey Levine, and Chelsea Finn.
\newblock Universal planning networks.
\newblock \emph{arXiv preprint arXiv:1804.00645}, 2018.

\bibitem[Taylor \& Stone(2009)Taylor and Stone]{taylor2009transfer}
Matthew~E Taylor and Peter Stone.
\newblock Transfer learning for reinforcement learning domains: A survey.
\newblock \emph{Journal of Machine Learning Research}, 10\penalty0 (7), 2009.

\bibitem[Torabi et~al.(2018)Torabi, Warnell, and Stone]{torabi2018behavioral}
Faraz Torabi, Garrett Warnell, and Peter Stone.
\newblock Behavioral cloning from observation.
\newblock \emph{arXiv preprint arXiv:1805.01954}, 2018.

\bibitem[Von~Rueden et~al.(2021)Von~Rueden, Mayer, Beckh, Georgiev,
  Giesselbach, Heese, Kirsch, Pfrommer, Pick, Ramamurthy,
  et~al.]{von2021informed}
Laura Von~Rueden, Sebastian Mayer, Katharina Beckh, Bogdan Georgiev, Sven
  Giesselbach, Raoul Heese, Birgit Kirsch, Julius Pfrommer, Annika Pick,
  Rajkumar Ramamurthy, et~al.
\newblock Informed machine learning--a taxonomy and survey of integrating prior
  knowledge into learning systems.
\newblock \emph{IEEE Transactions on Knowledge and Data Engineering},
  35\penalty0 (1):\penalty0 614--633, 2021.

\bibitem[Wang et~al.(2014)Wang, Man, Shen, Cao, Zheng, Jin,
  et~al.]{wang2014robust}
Hai Wang, Zhihong Man, Weixiang Shen, Zhenwei Cao, Jinchuan Zheng, Jiong Jin,
  et~al.
\newblock Robust control for steer-by-wire systems with partially known
  dynamics.
\newblock \emph{IEEE Transactions on Industrial Informatics}, 10\penalty0
  (4):\penalty0 2003--2015, 2014.

\bibitem[Watter et~al.(2015)Watter, Springenberg, Boedecker, and
  Riedmiller]{watter2015embed}
Manuel Watter, Jost Springenberg, Joschka Boedecker, and Martin Riedmiller.
\newblock Embed to control: A locally linear latent dynamics model for control
  from raw images.
\newblock In \emph{Advances in Neural Information Processing Systems}, pp.\
  2746--2754, 2015.

\bibitem[Wu et~al.(2018)Wu, Xiao, and Paterson]{wu2018physics}
Jin-Long Wu, Heng Xiao, and Eric Paterson.
\newblock Physics-informed machine learning approach for augmenting turbulence
  models: A comprehensive framework.
\newblock \emph{Physical Review Fluids}, 3\penalty0 (7):\penalty0 074602, 2018.

\bibitem[Xu et~al.(2020)Xu, Roosta, and Mahoney]{xu2020second}
Peng Xu, Fred Roosta, and Michael~W Mahoney.
\newblock Second-order optimization for non-convex machine learning: An
  empirical study.
\newblock In \emph{SIAM International Conference on Data Mining}, pp.\
  199--207, 2020.

\bibitem[Zhang et~al.(2018)Zhang, Sukhbaatar, Lerer, Szlam, and
  Fergus]{zhang2018composable}
Amy Zhang, Sainbayar Sukhbaatar, Adam Lerer, Arthur Szlam, and Rob Fergus.
\newblock Composable planning with attributes.
\newblock In \emph{International Conference on Machine Learning}, pp.\
  5842--5851, 2018.

\bibitem[Zhang et~al.(2019)Zhang, Zhang, Lu, Zhu, and Dong]{zhang2019you}
Dinghuai Zhang, Tianyuan Zhang, Yiping Lu, Zhanxing Zhu, and Bin Dong.
\newblock You only propagate once: Painless adversarial training using maximal
  principle.
\newblock \emph{arXiv preprint arXiv:1905.00877}, 2019.

\bibitem[Ziebart et~al.(2008)Ziebart, Maas, Bagnell, and
  Dey]{ziebart2008maximum}
Brian~D Ziebart, Andrew Maas, J~Andrew Bagnell, and Anind~K Dey.
\newblock Maximum entropy inverse reinforcement learning.
\newblock In \emph{AAAI Conference on Artificial Intelligence}, pp.\
  1433--1438, 2008.

\end{thebibliography}

\newpage
\appendix
\begin{appendix}

{\centering\textbf{\Large Appendix}}
\section{Proof of Lemma \ref{lemma:convergenceKalman}} \label{Appendix:ProofLemma2}
   % We will start with building relationship between the estimation error and prediction error:
   %  \begin{align}
   %      \tilde{\boldsymbol{\theta}}^-_t&=\boldsymbol{\theta}_t-\hat{\boldsymbol{\theta}}^-_t\\
   %      &=\boldsymbol{\theta}_{t-1}-\hat{\boldsymbol{\theta}}_{t-1}\\
   %      &=\tilde{\boldsymbol{\theta}}_{t-1}. \label{proof:tildetheta}
   %  \end{align}
    Since the matrices $\boldsymbol{P}_t$ and $\boldsymbol{L}_t$ are bounded according to Assumption \ref{assumption:observable} and \ref{assumption:Lbounded}, from (\ref{eq:EKFupdate}), one will have:
    \begin{align}
        \boldsymbol{K}_t&=\boldsymbol{P}_t\boldsymbol{L}_t^\prime \boldsymbol{R}_t^{-1} \label{proof:Kt1}\\
&=\boldsymbol{P}^{-}_t\boldsymbol{L}_{t}^\prime(\boldsymbol{L}_t\boldsymbol{P}^{-}_{t}\boldsymbol{L}^\prime_t+\boldsymbol{R}_t)^{-1}. \label{proof:Kt2}
    \end{align}
    Then, by taking the inverse of (\ref{proof:Kt1}) and (\ref{proof:Kt2}), one will get:
    \begin{equation} \label{proof:Ptinv}
        \boldsymbol{P}^{-1}_t=(\boldsymbol{P}_t^-)^{-1}+\boldsymbol{L}_t^\prime \boldsymbol{R}_t^{-1}\boldsymbol{L}_t.
    \end{equation}

    Substituting (\ref{proof:Kt1}) into (\ref{eq:EKFupdate}) and subtracting both sides from $\boldsymbol{\theta}_t$, one will have:
    \begin{equation} \label{eq:EKFupdate2}
        \tilde{\boldsymbol{\theta}}_{t}=\tilde{\boldsymbol{\theta}}^{-}_t - \boldsymbol{P}_t\boldsymbol{L}_t^\prime \boldsymbol{R}_t^{-1}\boldsymbol{e}_t.
    \end{equation}
    Then, plug (\ref{eq:EKFupdate2}) into the Lyapunov function (\ref{eq:Lyapunov}):
    \begin{align}\label{eq:Lyapunov2}
    V_t&=\tilde{\boldsymbol{\theta}}^\prime_t\boldsymbol{P}^{-1}_{t}\tilde{\boldsymbol{\theta}}_t\\
        &=(\tilde{\boldsymbol{\theta}}^{-}_t - \boldsymbol{P}_t\boldsymbol{L}_t^\prime \boldsymbol{R}_t^{-1}\boldsymbol{e}_t)^\prime\boldsymbol{P}^{-1}_{t}(\tilde{\boldsymbol{\theta}}^{-}_t - \boldsymbol{P}_t\boldsymbol{L}_t^\prime \boldsymbol{R}_t^{-1}\boldsymbol{e}_t)\\
        &=(\tilde{\boldsymbol{\theta}}^{-}_t)^\prime\boldsymbol{P}^{-1}_{t}\tilde{\boldsymbol{\theta}}^{-}_t-(\tilde{\boldsymbol{\theta}}^{-}_t)^\prime\boldsymbol{L}_t^\prime \boldsymbol{R}_t^{-1}\boldsymbol{e}_t-\boldsymbol{e}_t^\prime \boldsymbol{R}_t^{-1}\boldsymbol{L}_t\tilde{\boldsymbol{\theta}}^{-}_t+\boldsymbol{e}_t ^\prime\boldsymbol{R}_t^{-1}\boldsymbol{L}_t\boldsymbol{P}_t\boldsymbol{L}_t^\prime \boldsymbol{R}_t^{-1}\boldsymbol{e}_t \label{proof:Lyaounov3}
    \end{align}
    Next, we plug (\ref{proof:Ptinv}) into (\ref{proof:Lyaounov3}):
    \begin{align}
        V_t&=(\tilde{\boldsymbol{\theta}}^{-}_t)^\prime((\boldsymbol{P}_t^-)^{-1}+\boldsymbol{L}_t^\prime \boldsymbol{R}_t^{-1}\boldsymbol{L}_t)\tilde{\boldsymbol{\theta}}^{-}_t-(\tilde{\boldsymbol{\theta}}^{-}_t)^\prime\boldsymbol{L}_t^\prime \boldsymbol{R}_t^{-1}\boldsymbol{e}_t-\boldsymbol{e}_t^\prime \boldsymbol{R}_t^{-1}\boldsymbol{L}_t\tilde{\boldsymbol{\theta}}^{-}_t+\boldsymbol{e}_t^\prime \boldsymbol{R}_t^{-1}\boldsymbol{L}_t\boldsymbol{P}_t\boldsymbol{L}_t^\prime \boldsymbol{R}_t^{-1}\boldsymbol{e}_t\\      
        &=V_t^-+(\tilde{\boldsymbol{\theta}}^{-}_t)^\prime\boldsymbol{L}_t^\prime \boldsymbol{R}_t^{-1}\boldsymbol{L}_t\tilde{\boldsymbol{\theta}}^{-}_t-(\tilde{\boldsymbol{\theta}}^{-}_t)^\prime\boldsymbol{L}_t^\prime \boldsymbol{R}_t^{-1}\boldsymbol{e}_t-\boldsymbol{e}_t^\prime\boldsymbol{R}_t^{-1}\boldsymbol{L}_t\tilde{\boldsymbol{\theta}}^{-}_t+\boldsymbol{e}_t^\prime\boldsymbol{R}_t^{-1}\boldsymbol{L}_t\boldsymbol{P}_t\boldsymbol{L}_t^\prime \boldsymbol{R}_t^{-1}\boldsymbol{e}_t, \label{proof:Lyapunov4}
    \end{align}
    where
    \begin{align}
        V_t^-&=(\tilde{\boldsymbol{\theta}}^{-}_t)^\prime(\boldsymbol{P}_t^-)^{-1}\tilde{\boldsymbol{\theta}}^{-}_t\\
        &=(\tilde{\boldsymbol{\theta}}_{t-1})^\prime\mathbfcal{G}_t^\prime\boldsymbol{P}_{t-1}^{-1}\mathbfcal{G}_t\tilde{\boldsymbol{\theta}}_{t-1}.
    \end{align}
    Using (\ref{eq:errorapproximate}), (\ref{proof:Lyapunov4}) becomes:
    \begin{align}
        V_t &= V_t^-+(\tilde{\boldsymbol{\theta}}^{-}_t)^\prime\boldsymbol{L}_t^\prime \boldsymbol{R}_t^{-1}\boldsymbol{L}_t\tilde{\boldsymbol{\theta}}^{-}_t-(\tilde{\boldsymbol{\theta}}^{-}_t)^\prime\boldsymbol{L}_t^\prime \boldsymbol{R}_t^{-1}\boldsymbol{e}_t-\boldsymbol{e}_t^\prime \boldsymbol{R}_t^{-1}\boldsymbol{L}_t\tilde{\boldsymbol{\theta}}^{-}_t+\boldsymbol{e}_t ^\prime\boldsymbol{R}_t^{-1}\boldsymbol{L}_t\boldsymbol{P}_t\boldsymbol{L}_t^\prime \boldsymbol{R}_t^{-1}\boldsymbol{e}_t\\
        &= V_t^-+\boldsymbol{e}_t^\prime\mathbfcal{F}_t \boldsymbol{R}_t^{-1}\mathbfcal{F}_t\boldsymbol{e}_t-\boldsymbol{e}_t^\prime\mathbfcal{F}_t \boldsymbol{R}_t^{-1}\boldsymbol{e}_t-\boldsymbol{e}_t^\prime\boldsymbol{R}_t^{-1}\mathbfcal{F}_t\boldsymbol{e}_t+\boldsymbol{e}_t^\prime \boldsymbol{R}_t^{-1}\boldsymbol{L}_t\boldsymbol{P}_t\boldsymbol{L}_t^\prime \boldsymbol{R}_t^{-1}\boldsymbol{e}_t\\
        &= V_t^-+\boldsymbol{e}_t^\prime(\mathbfcal{F}_t \boldsymbol{R}_t^{-1}\mathbfcal{F}_t-\mathbfcal{F}_t \boldsymbol{R}_t^{-1}-\boldsymbol{R}_t^{-1}\mathbfcal{F}_t+ \boldsymbol{R}_t^{-1}\boldsymbol{L}_t\boldsymbol{P}_t\boldsymbol{L}_t^\prime \boldsymbol{R}_t^{-1})\boldsymbol{e}_t.
    \end{align}
    To ensure that the Lyapunov function $\{V_t\}_{t=1,2,\hdots}$ is a decreasing sequence, $V_{t}-V_{t-1}\le0$.
    \begin{align}
        V_{t}-&V_{t-1}\\
        &=\boldsymbol{e}_t^\prime(\mathbfcal{F}_t \boldsymbol{R}_t^{-1}\mathbfcal{F}_t-\mathbfcal{F}_t \boldsymbol{R}_t^{-1}-\boldsymbol{R}_t^{-1}\mathbfcal{F}_t+ \boldsymbol{R}_t^{-1}\boldsymbol{L}_t\boldsymbol{P}_t\boldsymbol{L}_t^\prime \boldsymbol{R}_t^{-1})\boldsymbol{e}_t\\
        &+(\tilde{\boldsymbol{\theta}}_{t-1})^\prime(\mathbfcal{G}_t^\prime\boldsymbol{P}_{t-1}^{-1}\mathbfcal{G}_t-\boldsymbol{P}_{t-1}^{-1})\tilde{\boldsymbol{\theta}}_{t-1}\le0.
    \end{align}
    Therefore, to ensure the Lyapunov function is a decreasing sequence,
    \begin{equation} \label{proof:suffcond1}
        \mathbfcal{F}_t \boldsymbol{R}_t^{-1}\mathbfcal{F}_t-\mathbfcal{F}_t \boldsymbol{R}_t^{-1}-\boldsymbol{R}_t^{-1}\mathbfcal{F}_t+ \boldsymbol{R}_t^{-1}\boldsymbol{L}_t\boldsymbol{P}_t\boldsymbol{L}_t^\prime \boldsymbol{R}_t^{-1}\le0,
    \end{equation}
    and
    \begin{equation} \label{proof:suffcond2}
        \mathbfcal{G}_t^\prime\boldsymbol{P}_{t-1}^{-1}\mathbfcal{G}_t-\boldsymbol{P}_{t-1}^{-1}\le0.
    \end{equation}
With some manipulations:
\begin{align} \label{proof:suffcond11}
(\mathbfcal{F}_t-\boldsymbol{I}_{s})\boldsymbol{R}_t^{-1}(\mathbfcal{F}_t-\boldsymbol{I}_{s})-\boldsymbol{R}_t^{-1}+\boldsymbol{R}_t^{-1}\boldsymbol{L}_t\boldsymbol{P}_t\boldsymbol{L}_t^\prime \boldsymbol{R}_t^{-1}&\le0,\\
(\mathbfcal{F}_t-\boldsymbol{I}_{s})\boldsymbol{R}_t^{-1}(\mathbfcal{F}_t-\boldsymbol{I}_{s})-\boldsymbol{R}_t^{-1}(\boldsymbol{I}_s-\boldsymbol{L}_t\boldsymbol{P}^{-}_t\boldsymbol{L}_{t}^\prime(\boldsymbol{L}_t\boldsymbol{P}^{-}_{t}\boldsymbol{L}^\prime_t+\boldsymbol{R}_t)^{-1})&\le0.
\end{align}
By letting $\boldsymbol{I}_s=(\boldsymbol{L}_t\boldsymbol{P}^{-}_{t}\boldsymbol{L}^\prime_t+\boldsymbol{R}_t)(\boldsymbol{L}_t\boldsymbol{P}^{-}_{t}\boldsymbol{L}^\prime_t+\boldsymbol{R}_t)^{-1}$, we have
\begin{align}
(\mathbfcal{F}_t-\boldsymbol{I}_{s})\boldsymbol{R}_t^{-1}(\mathbfcal{F}_t-\boldsymbol{I}_{s})-(\boldsymbol{L}_t\boldsymbol{P}^{-}_{t}\boldsymbol{L}^\prime_t+\boldsymbol{R}_t)^{-1}&\le0.
    \end{align}
Since $\mathbfcal{F}_t$ and $\boldsymbol{R}_t$ are diagonal matrices, we will have 
\begin{align}
\boldsymbol{R}_t^{-1}(\mathbfcal{F}_t-\boldsymbol{I}_{s})^2-(\boldsymbol{L}_t\boldsymbol{P}^{-}_{t}\boldsymbol{L}^\prime_t+\boldsymbol{R}_t)^{-1}&\le0,
    \end{align}
which at the end yields:
\begin{equation}
    (\mathbfcal{F}_t-\boldsymbol{I}_{s})^2\le\boldsymbol{R}_t(\boldsymbol{L}_t\boldsymbol{P}^{-}_{t}\boldsymbol{L}^\prime_t+\boldsymbol{R}_t)^{-1},
\end{equation}
therefore the proof is completed.

\section{Proof of Theorem \ref{theorem}} \label{Appendix:ProofTheorem}

This proof is straightforward once Lemma~\ref{lemma:convergenceKalman} is provided.
Consider the assumptions \ref{assumption:observable} and \ref{assumption:Lbounded} are met, according to Lemma \ref{lemma:convergenceKalman}, with the exact gradient generated by the gradient generator in (\ref{eq:backwardPass})-(\ref{eq:recursiveGG}), $ \lim_{t\rightarrow\infty}\tilde{\boldsymbol{\theta}}_t=0$. As the estimated $\hat{\boldsymbol{\theta}}_t$ converges to the true $\boldsymbol{\theta}^*$, where the true parameter gives zero cumulative loss, the cumulative loss $L(\boldsymbol{\xi}(\hat{\boldsymbol{\theta}}))$ goes to $0$.

\section{Experiment Details} \label{Appendix:experimentalDetails}

\subsection{System/Environment Setups}

\textbf{Cartpole.} We consider the following continuous dynamics of the cartpole 
\begin{equation} \label{dyn_cartpole_cont}
\begin{bmatrix}
\dot{p}\\
\ddot{p}\\
\dot{\theta}\\
\ddot{\theta}
\end{bmatrix}
=
\begin{bmatrix}
\dot{p}\\
(F+\frac{m_p l\dot{\theta}^2 \sin(\theta)}{m_t}) - \frac{m_p l \ddot{\theta} \cos(\theta)}{m_t}\\
\dot{\theta}\\
\frac{g\sin(\theta) - \cos(\theta)(F+\frac{m_p l\dot{\theta}^2 \sin(\theta)}{m_t})}{l(\frac{4}{3} - \frac{m_{p}\cos(\theta)^2}{m_{t}})}
\end{bmatrix}
,
\end{equation}
where $p\in\mathbb{R}$ is the horizontal displacement of the cart; $\theta\in\mathbb{R}$ is the pole angle; $F\in\mathbb{R}$ denotes the horizontal force applied to the cart which is between $-1$ and $+1$; $l\in\mathbb{R}$ is the length of the pole; $m_p, m_t\in\mathbb{R}$ are the masses of the pole and total cartpole, respectively. By defining  the states and control inputs of the cartpole 
\begin{equation}\label{cart_dyn}
\boldsymbol{x} \triangleq
\begin{bmatrix}
p&\dot{p}&\theta&\dot{\theta}
\end{bmatrix}^{\prime}
\quad \text{and} \quad 
\boldsymbol u\triangleq F
\end{equation}
respectively.

% \textbf{Two-link Robot Arm.} We consider that a two-link robot arm moves in vertical plane with  continuous dynamics given by \cite[p. 209]{spong2008robot}

% \begin{equation} \label{armeDyn}
% M(\boldsymbol\theta)\ddot{\boldsymbol\theta}+C(\boldsymbol\theta,\dot{\boldsymbol\theta})\dot{\boldsymbol\theta}+\boldsymbol g(\boldsymbol\theta)=\boldsymbol\tau,
% \end{equation}
% where $\boldsymbol\theta=[\theta_1, \theta_2]^{\prime} \in \mathbb{R}^2$ is the angle vector of the joints; $M(\boldsymbol \theta)\in \mathbb{R}^{2\times2}$ is the inertia matrix; $C(\boldsymbol\theta,\dot{\boldsymbol\theta}) \in \mathbb{R}^{2\times2}$ is the Coriolis matrix; $\boldsymbol g(\boldsymbol\theta) \in \mathbb{R}^2$ is the gravity vector; and $\boldsymbol\tau=[\tau_1, \tau_2]^{\prime}\in \mathbb{R}^2$ are the torques applied to each joint. We adopt the parameters in \cite[p. 209]{spong2008robot}: the link mass $m_1=m_2=1 \mathrm{kg}$, the link length $l_1=l_2=1 \mathrm{m}$;  the distance from joint to center of mass (COM) $r_{1}=r_{2}=0.5 \mathrm{m}$, and the moment of inertia with respect to COM $I_1=I_2=1/12 \mathrm{kgm^2}$.	
% By defining  the  states and control inputs of the robot arm system
% \begin{equation}
% \boldsymbol x\triangleq\begin{bmatrix}
% \theta_1&\dot{\theta_1}&\theta_2&\dot{\theta_2}
% \end{bmatrix}^{\prime}\quad \text{and} \quad \boldsymbol u\triangleq\boldsymbol\tau=
% \begin{bmatrix}
% \tau_1&\tau_2
% \end{bmatrix}^{\prime},
% \end{equation}
% respectively.

\textbf{Quadrotor UAV.} We consider a quadrotor UAV with the following dynamics

\begin{equation} \label{eq:quadrotor_dyn}
    \begin{aligned} 
         \dot{\boldsymbol{p}}_I &= \boldsymbol{v}_I, \\
         m\dot{\boldsymbol{v}}_I &= m\boldsymbol{g}_I+\mathbf{F}_I, \\
         \dot{\boldsymbol{q}}_{B/I} &= \frac{1}{2}\boldsymbol{\Omega}(\boldsymbol{\omega}_B)\boldsymbol{q}_{B/I}, \\
         J_B\dot{\boldsymbol{\omega}}_B &= \mathbf{M}_B - \boldsymbol{\omega} \times J_B\boldsymbol{\omega}_B.
    \end{aligned}
\end{equation}
Here, the subscription $_B$ and $_I$ denote a quantity expressed in the body frame and inertial (world) frame, respectively; $m$ and $J_B\in\mathbb{R}^{3\times 3}$ are the mass and moment of inertia with respect to the body frame of the UAV, respectively. $g$ is the gravitational constant ($g=10\text{ m}/\text{s}^2$), $\boldsymbol{g}_I=[0,0,g]^\prime$. $\boldsymbol{p}\in\mathbb{R}^3$  and $\boldsymbol{v}\in\mathbb{R}^3$ are the position and velocity vector of the UAV; $\boldsymbol\omega_B\in\mathbb{R}^3$ is the angular velocity vector of the UAV; $\boldsymbol{q}_{B/I}\in\mathbb{R}^4$ is the unit quaternion \cite{kuipers1999quaternions} that describes the attitude of the UAV with respect to the inertial frame; $\boldsymbol{\Omega}(\boldsymbol\omega_B)$ is defined as:
\begin{equation}
    \boldsymbol{\Omega}(\boldsymbol\omega_B) =
    \begin{bmatrix}
    0&-\omega_x&-\omega_y&-\omega_z \\
    \omega_x&0&\omega_z &-\omega_y \\
    \omega_y&-\omega_z &0&\omega_x \\
    \omega_z &\omega_y&-\omega_x&0
    \end{bmatrix},
\end{equation}
$\mathbf{M}_B\in\mathbb{R}^3$ is the torque applied to the UAV; $\mathbf{F}_I\in\mathbb{R}^3$ is the force vector applied to the UAV center of mass. The total force magnitude $f=\|\mathbf{F}_I\|\in\mathbb{R}$ (along z-axis of the body frame) and torque $\mathbf{M}_B=[M_x,M_y,M_z]^\prime$ are generated by thrust from four rotating propellers $[T_1, T_2, T_3, T_4]^\prime$, their relationship can be expressed as:

\begin{equation}
    \begin{bmatrix}
    f\\
    M_x\\
    M_y\\
    M_z
    \end{bmatrix}
    =
    \begin{bmatrix}
    1&1&1&1 \\
    0&-l_w/2&0&l_w/2 \\
    -l_w/2&0&l_w/2&0 \\
    c&-c&c&-c
    \end{bmatrix}
    \begin{bmatrix}
    T_1\\
    T_2\\
    T_3\\
    T_4
    \end{bmatrix},
\end{equation}
where $l_w$ is the wing length of the UAV and $c$ is a fixed constant.
The state and input vectors of the UAV are defined as:
\begin{equation}
    \begin{aligned}
        \boldsymbol{x} &\triangleq
        \begin{bmatrix}
        \boldsymbol{p}'& \boldsymbol{v}'& \boldsymbol{q}' & \boldsymbol\omega'
        \end{bmatrix}
        ' \in \mathbb{R}^{13}, \\
        \boldsymbol{u} &\triangleq
        \begin{bmatrix}
        T_1&T_2&T_3&T_4
        \end{bmatrix}
        '\in\mathbb{R}^4.
    \end{aligned}
\end{equation}

\textbf{Rocket.} The rocket is treated as a rigid body subject to constant gravitational acceleration, $g_I\in\mathbb{R}^3$, and neglects aerodynamic forces. The vehicle is assumed to actuate a single gimbaled rocket engine to generate a thrust vector within a feasible range of magnitudes and gimbal angles. We assume that at the landing phase, the depletion of fuel is insignificant. Therefore, we omit the dynamics of rocket mass. The rocket has the following dynamics:
\begin{equation}
    \begin{aligned} 
         \dot{\boldsymbol{p}}_{\mathcal{I}} &= \boldsymbol{v}_{\mathcal{I}}, \\
         \dot{\boldsymbol{v}}_{\mathcal{I}} &= \frac{1}{m}C_{\mathcal{I}/\mathcal{B}}\mathbf{T}_{\mathcal{B}}+\mathbf{g}_{\mathcal{I}}, \\
         \dot{\boldsymbol{q}}_{\mathcal{B}/\mathcal{I}} &= \frac{1}{2}\Omega(\boldsymbol{\omega}_{\mathcal{B}})\boldsymbol{q}_{\mathcal{B}/\mathcal{I}}, \\
         J_{\mathcal{B}}\dot{\boldsymbol{\omega}}_{\mathcal{B}} &= \mathbf{M}_{\mathcal{B}}- [\boldsymbol{\omega}_{\mathcal{B}} \times] J_{\mathcal{B}}\boldsymbol{\omega}_{\mathcal{B}}.
    \end{aligned}
\end{equation}

Here, the subscription $_{\mathcal{B}}$ and $_{\mathcal{I}}$ denote a quantity expressed in the body frame and inertial (world) frame, respectively; $m$ and $J_{\mathcal{B}}\in\mathbb{R}^{3\times3}$ are the mass and moment of inertia with respect to body frame of the rocket, respectively. $\boldsymbol{p}\in\mathbb{R}^3$  and $\boldsymbol{v}\in\mathbb{R}^3$ are the position and velocity vector of the rocket; $\boldsymbol\omega_{\mathcal{B}}\in\mathbb{R}^3$ is the angular velocity vector of the rocket; $\boldsymbol{q}_{\mathcal{B}/\mathcal{I}}=[q_0,q_1,q_2,q_3]$ is the unit quaternion that describes the attitude of rocket with respect to the inertial frame; $\mathbf{T}_{\mathcal{B}}\in\mathbb{R}^3$ is the commanded thrust vector; $\mathbf{M}_{\mathcal{B}}\in\mathbb{R}^3$ is the torque applied to the rocket; $C_{\mathcal{B}/\mathcal{I}}$ is the direction cosine matrix that encodes the attitude transformation from body frame to inertia frame and related to $\boldsymbol{q}_{\mathcal{B}/\mathcal{I}}$ by the following relationship:

\begin{equation*}
    C_{\mathcal{B}/\mathcal{I}} =
    \begin{bmatrix}
    1-2(q_2^2+q_3^2)&2(q_1q_2+q_0q_3)&2(q_1q_3-q_0q_2) \\
    2(q_1q_2-q_0q_3)&1-2(q_1^2+q_3^2)&2(q_2q_3+q_0q_1)\\
    2(q_1q_3+q_0q_2)&2(q_2q_3-q_0q_1)&1-2(q_1^2+q_1^2)\\
    \end{bmatrix},
\end{equation*}
The inverse transformation is denoted as $C_{\mathcal{I}/\mathcal{B}}=C_{\mathcal{B}/\mathcal{I}}^T$;

The skew-symmetric matrices $[\boldsymbol{\omega}_{\mathcal{B}} \times]$ and $\Omega(\boldsymbol\omega_{\mathcal{B}})$ are defined as follow:

\begin{equation*}
    [\boldsymbol{\omega}_{\mathcal{B}} \times] \triangleq
    \begin{bmatrix}
    0&-\omega_z&\omega_y\\
    \omega_z&0&-\omega_x\\
    -\omega_y&\omega_x &0\\
    \end{bmatrix},
\quad
    \Omega(\boldsymbol\omega_{\mathcal{B}}) \triangleq
    \begin{bmatrix}
    0&-\omega_x&-\omega_y&-\omega_z \\
    \omega_x&0&\omega_z &-\omega_y \\
    \omega_y&-\omega_z &0&\omega_x \\
    \omega_z &\omega_y&-\omega_x&0
    \end{bmatrix},
\end{equation*}

The state and input vectors of the rocket are defined as:
\begin{equation}
    \begin{aligned}
        \boldsymbol{x} &=
        \begin{bmatrix}
        \boldsymbol{p}_{\mathcal{I}}'& \boldsymbol{v}_{\mathcal{I}}'& \boldsymbol{q}_{\mathcal{B}/\mathcal{I}}' & \boldsymbol\omega_{\mathcal{B}}'
        \end{bmatrix}
        ' \in \mathbb{R}^{13}, \\
        \boldsymbol{u} &=\mathbf{T}_{\mathcal{B}}=
        \begin{bmatrix}
        T_x&T_y&T_z
        \end{bmatrix}
        '\in\mathbb{R}^3,
    \end{aligned}
\end{equation}

\textbf{Discretization.} Discretization is done by the following discrete-time  form
\begin{equation}\label{equ_discretization}
\boldsymbol x_{t+1}\approx\boldsymbol x_t+\Delta\cdot\boldsymbol g(\boldsymbol x_t, \boldsymbol u_{t})\triangleq\boldsymbol{f}(\boldsymbol x_t, \boldsymbol u_{t}),
\end{equation}
where $\Delta$ is the discretization interval.

\subsection{Online Imitation Learning}
\textbf{Data acquisition.} The dataset of expert demonstrations is generated by solving an 
optimal control system with the true dynamics and control objective parameter $\boldsymbol{\theta}^*=\{\boldsymbol{\theta}_{dyn},\boldsymbol{\theta}_{obj}\}$
given. We generate five trajectories with different initial conditions $\boldsymbol{x}_0$ and time horizons $T$.

\textbf{PDP.} We employed the PDP in \cite{jin2020pontryagin} to solve this problem. The learning rate is $\eta=10^{-4}$. Five trials were run given random initial $\boldsymbol{\theta}_0$.

\textbf{Inverse KKT method.} We choose the inverse KKT method \cite{englert2017inverse} for comparison because it is
suitable for learning objective functions for high-dimensional continuous-space systems. We adopt the inverse KKT method and define the KKT loss as the norm-2 violation of the KKT condition by the demonstration data:
\begin{equation*}
\min_{\boldsymbol{\theta},\boldsymbol{\lambda}_{1:T}} \left( \left\| \frac{\partial L}{\partial \boldsymbol{x}_{0:T}} (\boldsymbol{x}^*_{0:T}, \boldsymbol{u}^*_{0:T-1}) \right\|^2 + \left\| \frac{\partial L}{\partial \boldsymbol{u}_{0:T-1}} (\boldsymbol{x}^*_{0:T}, \boldsymbol{u}^*_{0:T-1}) \right\|^2 \right).
\end{equation*}

\textbf{Neural policy cloning.} For the neural policy cloning, we directly learn a neural network policy $\boldsymbol{u} = \boldsymbol{\mu}(\boldsymbol{x},\boldsymbol{\theta})$ from the dataset using supervised learning, that is
\begin{equation}
\min_{\boldsymbol{\theta}}\sum_{t=0}^{T-1}\|\boldsymbol{u}^*_t-\boldsymbol{\mu}(\boldsymbol{x}^*_t,\boldsymbol{\theta})\|^2
\end{equation}

\subsection{Online System Identification}
\textbf{Data acquisition.} In the system identification experiment, we collect a total number of five
trajectories from systems with dynamics known, wherein different trajectories $\boldsymbol{\xi}^o=\{\boldsymbol{x}^o_{0:T},\boldsymbol{u}_{0:T-1}\}$ have different initial conditions $\boldsymbol{x}_0$ and horizons $T$ ($T$ ranges from 10 to 20 depending on different environment and task), with
random inputs $\boldsymbol{u}_{0:T-1}$ drawn from uniform distribution.

\textbf{PDP.} We employed the PDP in \cite{jin2020pontryagin} to solve this problem. The learning rate is $\eta=10^{-4}$. Five trials were run given random initial $\boldsymbol{\theta}_0$. For the neural dynamics case, the learning rate is $\eta=10^{-5}$.

\textbf{Pytorch Adam to learn neural dynamics.} We consider the dynamics of each system (cartpole, quadrotor, and rocket) are represented by a fully-connected feed-forward neural network $\hat{\boldsymbol{f}}(\boldsymbol{x}_t,\boldsymbol{u}_t,\boldsymbol{\theta})$. The neural network has a layer structure of $(n+m)$-$2(n+m)$-$n$ with $tanh$ activation functions, i.e., there is an
input layer with $(n+m)$ neurons equal to the dimension of state, one hidden layer with $2(n+m)$ neurons and one output layer with $n$ neurons. The $\boldsymbol{\xi}^o=\{\boldsymbol{x}^o_{0:T},\boldsymbol{u}_{0:T-1}\}$ obtained previously are used in stage loss. We conducted five trials for each method with different initial $\boldsymbol{\theta}$. We use Pytorch Adam to train the neural network by minimizing the following residual
\begin{equation}
\min_{\boldsymbol{\theta}}\sum_{t=0}^{T-1}\|\boldsymbol{x}^o_{t+1}-\hat{\boldsymbol{f}}(\boldsymbol{x}^o_t,\boldsymbol{u}_t,\boldsymbol{\theta})\|^2.
\end{equation}

\textbf{DMDc.} The DMDc method \cite{proctor2016dynamic} is a method that is based on Koopman theory to represent nonlinear dynamics with linear dynamics of observables. Observables $\boldsymbol{\psi}(\boldsymbol{x}_t)$ are some basis functions of states. The observable space has a much higher dimension compared to state space. The estimation of the dynamics is achieved by the following optimization:
\begin{equation}
\min_{\boldsymbol{A},\boldsymbol{B}}\sum_{t=0}^{T-1}\|\boldsymbol{\psi}(\boldsymbol{x}^o_{t+1})-\boldsymbol{A}\boldsymbol{\psi}(\boldsymbol{x}^o_t)-\boldsymbol{B}\boldsymbol{u}_t\|^2.
\end{equation}

\subsection{Policy Tuning On-the-fly}

\textbf{Neural State Feedback Policy.} In this application, we learn the parameters of a neural state feedback policy by minimizing given control objective functions. Specifically, we use a fully connected feed-forward neural network that has a layer structure of $3n$-$3n$-$m$ with $tanh$ activation functions, i.e., there is an input layer with $3n$ neurons equal to the dimension of state, one hidden layer with $3n$ neurons and one output layer with $m$ neurons. The policy parameter $\boldsymbol{\theta}$ is the neural network parameter. For comparison, we apply the guided policy search (GPS) method \cite{levine2014learning} and iLQR \cite{li2004iterative} to solve the same problem.

\textbf{PDP.} We employed the PDP in \cite{jin2020pontryagin} to solve this problem. The learning rate is set to be $\eta=10^{-4}$ or $=10^{-6}$. Five trials were run given random initial $\boldsymbol{\theta}_0$. For the neural objective function case, the learning rate is $\eta=10^{-5}$.

\end{appendix} 
\end{document}